\documentclass[reqno,10pt]{amsart}
\usepackage{xcolor}
\usepackage{hyperref}

\usepackage{enumerate}
\usepackage{amssymb}
\usepackage{comment}
\usepackage[all,cmtip]{xy}

\newtheorem{Proposition}{Proposition}[section]
\newtheorem{Definition}[Proposition]{Definition}
\newtheorem{Lemma}[Proposition]{Lemma}
\newtheorem{Theorem}[Proposition]{Theorem}
\newtheorem{Corollary}[Proposition]{Corollary}
\newtheorem{Remark}[Proposition]{Remark}

\newtheorem{Example}[Proposition]{Example}
\DeclareMathOperator{\Val}{Val}

\DeclareMathOperator{\Gr}{Gr}
\DeclareMathOperator{\AGr}{AGr}
\DeclareMathOperator{\Res}{Res}

\renewcommand{\Re}{\mathrm {Re}}
\renewcommand{\Im}{\mathrm {Im}}
\DeclareMathOperator{\vol}{vol}

\DeclareMathOperator{\Dens}{Dens}
\DeclareMathOperator{\Span}{Span}

\DeclareMathOperator{\Ker}{Ker}
\DeclareMathOperator{\Hom}{Hom}
\DeclareMathOperator{\Stab}{Stab}

\DeclareMathOperator{\Sym}{Sym}
\DeclareMathOperator{\Sp}{Sp}
\DeclareMathOperator{\CA}{CA}
\DeclareMathOperator{\Pfaff}{Pf}
\DeclareMathOperator{\tr}{Tr}
\DeclareMathOperator{\OO}{O}
\DeclareMathOperator{\SO}{SO}

\DeclareMathOperator{\Alt}{Alt}

\DeclareMathOperator{\Kl}{Kl}

\DeclareMathOperator{\Cr}{Cr}
\DeclareMathOperator{\WF}{WF}

\DeclareMathOperator{\Id}{Id}
\DeclareMathOperator{\ind}{ind}

\DeclareMathOperator{\GL}{GL}
\DeclareMathOperator{\SL}{SL}
\DeclareMathOperator{\sign}{sign}

\newcommand{\R}{\mathbb{R}}
\newcommand{\C}{\mathbb{C}}

\title{Contact integral geometry and the Heisenberg algebra}

\author{Dmitry Faifman}

\email{dmitry.faifman@umontreal.ca}

\address{Centre de recherches math\'ematiques, Universit\'e de Montr\'eal, Pavillon Andr\'e-Aisenstadt,
	2920, Chemin de la tour, bur. 5357 Montréal (Québec) H3T 1J4 Canada}

\date{}

\begin{document}
\maketitle

\begin{abstract}
Generalizing Weyl's tube formula and building on Chern's work, Alesker reinterpreted the Lipschitz-Killing curvature integrals as a family of valuations (finitely-additive measures with good analytic properties), attached canonically to any Riemannian manifold, which is universal with respect to isometric embeddings. In this note, we uncover a similar structure for contact manifolds. Namely, we show that a contact manifold admits a canonical family of generalized valuations, which are universal under contact embeddings. Those valuations assign numerical invariants to even-dimensional submanifolds, which in a certain sense measure the curvature at points of tangency to the contact structure. 
Moreover, these valuations generalize to the class of manifolds equipped with the structure of a Heisenberg algebra on their cotangent bundle. Pursuing the analogy with Euclidean integral geometry, we construct symplectic-invariant distributions on Grassmannians to produce Crofton formulas on the contact sphere. Using closely related distributions, we obtain Crofton formulas also in the linear symplectic space.

\end{abstract}

\section{Introduction}

\subsection{Background and motivation}

Intrinsic volumes first appeared in convex geometry through Steiner's formula: given a compact convex body $K\subset \R^n$, $\vol(K+\epsilon B^n)=\sum_{k=0}^n \omega_{n-k}\mu_k(K) \epsilon^{n-k}$ where $B^j$ is unit Euclidean ball in $\R^j$ and $\omega_j$ is its volume. The coefficient $\mu_k(K)$, the $k$-th intrinsic volume, can be written explicitly for smooth $K$ as $\mu_k(K)=c_{n,k}\int_{\partial K} \sigma_{n-1-k}(\kappa_1,\dots,\kappa_{n-1})d\textrm{Area}_K$, where $\kappa_j$ are the principal curvatures of $\partial K$.
Alternatively, $\mu_k(K)$ can be given by integral-geometric Crofton formulas: $\mu_{k}(K)=c'_{n,k}\int_{\AGr_{n-k}(\R^n)}\chi(K \cap E)dE$ where $dE$ is the rigid motion invariant measure on the affine Grassmannian, and $\chi$ is the Euler characteristic. A third, axiomatic definition, was given by Hadwiger, who described the intrinsic volumes as the unique rigid-motion invariant continuous finitely additive measures on compact convex sets.

A closely related famous result is Weyl's tube formula \cite{weyl}. It asserts that the volume of an $\epsilon$-tube around a Riemannian manifold $M$ embedded isometrically in Euclidean space $\R^N$ is a polynomial in $\epsilon\ll1$, whose coefficients are, remarkably, intrinsic invariants of the Riemannian manifold $M$, independent of the isometric embedding. These coefficients, now known as the intrinsic volumes of $M$, are intimately linked with the asymptotic expansion of the heat kernel, see \cite{donnelly}.

These results fall naturally in the domain of valuations on manifolds, a fairly young branch of valuation theory introduced by Alesker et al. in a sequence of works \cite{alesker1, alesker2,alesker3,alesker4}, see also \cite{alesker_fu_survey} for a survey. Valuation theory itself is a mixture of convex and integral geometry, originating in the early 20th century in works of Steiner, Blaschke, Chern and Santalo, as well as in Dehn's solution to Hilbert's third problem. Generally speaking, valuations are finitely additive measures on some family of nice subsets. In this note, there is typically some analytic restriction on the nature of the valuation, such as smoothness or smoothness with singularities, and the subsets are manifolds with corners or differentiable polyhedra.

Building on results of Chern, Alesker noticed a natural extension of Weyl's theorem for valuations: restricting to $M$ the intrinsic volumes of $\R^N$, (considered as valuations), yields an intrinsically defined family of valuations on $M$, now known as the Lipschitz-Killing valuations. Weyl's intrinsic volumes of $M$ are then the integrals of the corresponding valuations. In a recent work, Fu and Wannerer \cite{fu_wannerer} characterized the Lipschitz-Killing valuations as the unique family of valuations attached canonically (in a sense made precise therein) to arbitrary Riemannian manifolds that are universal to isometric embeddings.

Other spaces whose (smooth) valuation theories were considered in recent years include complex space forms \cite{bernig_fu_hermitian, BFS}, the quaternionic plane \cite{bernig_solanes_quaternions1, bernig_solanes_quaternions2}, the octonionic plane \cite{bernig_voide} and exceptional spheres \cite{solanes_wannerer}.

Numerous intriguing connections between convex and symplectic geometries are known to exist. To name a few: Viterbo's conjectured isoperimetric inequality for capacities of convex bodies \cite{viterbo}, was later shown by Artstein-Avidan, Karasev and Ostrover to imply Mahler's famous conjecture \cite{AAKO}. Capacities have been successfully studied up to a bounded factor using convex techniques, see \cite{AAMO} and \cite{gluskin_ostrover}.
In \cite{APB, APBTZ}, links are established between systolic geometry, contact geometry, Mahler's conjecture and the geometry of numbers. Sch\"affer's dual girth conjecture for normed spaces has been proved using symplectic techniques \cite{alvarez-paiva}, and generalized further in \cite{faifman_quotient} using hamiltonian group actions. In very recent works of Abbondandolo et al. \cite{abbondandolo}, some links are established between the geometry of the group of symplectomorphisms and systolic geometry of the 2-sphere. For an exposition of some of those connections, see \cite{ostrover}. 

The main objective for this work is to further explore the convex-symplectic link by studying the valuation theory of contact manifolds, using the Riemannian case and Weyl's principle as guides. 
\subsection{Informal summary}
 We find that like in the Riemannian setting, contact manifolds possess a canonical family of valuations associated to them. We describe those valuations in two ways: geometrically through a curvature-type formula, and also dynamically through the invariants of a certain vector field at its singular points.

The contact valuations satisfy Weyl's principle of universality under embeddings, similarly to the valuation extension of the Weyl principle. Let us emphasize the role played by valuations in this phenomenon: It so happens that contact valuations only assume non-zero values on even dimensional submanifolds, while they live on odd-dimensional contact manifolds. Thus unlike the Riemannian case, Weyl's principle in the contact setting is only manifested in its extended to valuations form, as the statement in the original form becomes vacuous: all integrals of the contact valuations vanish.

Contact valuations are in fact an instance of a natural collection of valuations associated to the larger class of manifolds whose cotangent spaces admit a smoothly varying structure of Heisenberg algebras. The Heisenberg algebra provides a unifying link between contact, symplectic and metric geometries.

This path leads us to consider the valuation theory of the dual Heisenberg algebra, invariant under the group of the automorphisms of the Heisenberg algebra, which is closely related to the symplectic group. From this perspective, this is another step in the study of the valuation theory of non-compact Lie groups, which up to now has only been considered for the indefinite orthogonal group in \cite{alesker_faifman, bernig_faifman, faifman_opq}, and in a somewhat different framework for the special linear group, see \cite{ludwig_reitzner_annals,ludwig_reitzner_99}.

Further similarity to the metric setting is exhibited by the contact sphere, where we prove Crofton formulas and a Hadwiger-type theorem, thus establishing an integral-geometric and an axiomatic description of the contact valuations. 

Finally in the last part we explore the valuation theory of linear symplectic spaces. We show that there are no non-trivial invariant valuations, but nevertheless one can write oriented Crofton formulas for the symplectic volume of submanifolds.

\subsection{Main results}

Let us very briefly recall or indicate the relevant notions. For precise definitions see sections \ref{sec:basics} and \ref{sec:valuations}.

A contact manifold $M^{2n+1}$ is given by a maximally non-integrable hyperplane distribution, namely a smooth field of tangent hyperplanes $H\subset TM$ s.t. locally one can find $\alpha\in\Omega^1(M)$ with $H=\Ker(\alpha)$ and $d\alpha|_H$ a non-degenerate 2-form.

A smooth valuation $\phi$ on an orientable manifold $M^n$, denoted $\phi\in\mathcal V^\infty(M)$, is a finitely additive measure on the compact differentiable polyhedra of $M$, denoted $\mathcal P(M)$, which has the form $\phi(X)=\int_X\mu+\int_{N^*X}\omega$ for some forms $\mu\in\Omega^n(M)$ and $\omega\in\Omega^{n-1}(S^*M)$. Here $S^*M$ is the cosphere bundle, and $N^*X$ is the conormal cycle of $X$, which is just the conormal bundle when $X$ is a manifold. Orientability is not essential, and is only assumed to simplify the exposition. 

There is a natural filtration $\mathcal W^\infty_n(M)\subset \dots\subset \mathcal W^\infty_0(M)=\mathcal V^\infty(M)$. Very roughly speaking, $\mathcal W^\infty_j(M)$ consists of valuations which are locally homogeneous of degree at least $j$, for example $\mathcal W^\infty_n(M)$ are the smooth measures on $M$.

The generalized valuations $\mathcal V^{-\infty}(M)$ are, roughly speaking, distributional valuations: we allow $\omega$ and $\mu$ to be currents rather than smooth forms. Generalized valuations can be naturally evaluated on sufficiently nice subsets $X \in \mathcal P(M)$. The filtration $\mathcal W_j^\infty(M)$ on $\mathcal V^{\infty}(M)$ extends to a filtration  $\mathcal W_j^{-\infty}(M)$ on $\mathcal V^{-\infty}(M)$.

\begin{Theorem}\label{mainthm:contact}
	To any contact manifold $M^{2n+1}$ with contact distribution $H\subset TM$ there are canonically associated, linearly independent generalized valuations $\phi^M_{2k}\in\mathcal V^{-\infty}(M)$, $0\leq k\leq n$. They have the following properties:
	\begin{enumerate}
		\item $\phi^M_0$ is the Euler characteristic, and $\phi^M_{2k}\in \mathcal W_{2k}^{-\infty}(M)\setminus  \mathcal W_{2k+1}^{-\infty}(M)$.
		\item $\phi_{2k}^M$ can be naturally evaluated on submanifolds in generic position relative to the contact structure (see Definition \ref{def:normally transversal}). For a generic closed hypersurface $F$, 
		
		\[ \phi^M_{2k}(F)=\sum_{T_pF=H_p}  \phi^M_{2k}(F,p)  \]
		where the \emph{local contact area} $\phi^M_{2k}(F,p)$ only depends on the germ of $F$ at $p$. It is described explicitly below in equations \eqref{eq:local_dynamics} and \eqref{eq:local_geometry}.
		
		\item Universality to restriction under embedding: if $i:N^{2m+1}\to M^{2n+1}$ is a contact embedding, then for $k\leq m$ one has $i^*\phi^M_{2k}=\phi^N_{2k}$.
		\item  For a $2k$-dimensional submanifold in general position $F$, $\phi^M_{2k}(F)\geq 0$, with equality if and only if there are no contact tangent points.
		\item The space of generalized valuations on $M$ invariant under all contactomorphisms of $M$ is spanned by $(\phi^M_{2k})_{k=0}^n$.
	\end{enumerate}
\end{Theorem}
\begin{Remark}
The universality with respect to embeddings is sometimes referred to as the Weyl principle. Thus we recover a Weyl principle in the contact setting. 
\end{Remark}	
	
	The local contact areas $\phi^M_{2k}(F,p)$ can be given explicitly in two different ways, through a geometric or a dynamical approach.
	\begin{itemize}
		\item From a dynamical point of view, $\phi^M_{2k}(F,p)$ encodes the invariants of the linearized vector field $B\in\mathcal X(F)$ representing the characteristic foliation. While there are many such vector fields, there is a distinguished choice, to linear order, at the critical points - choose an arbitrary contact form $\alpha$ near $p$ and let $B$ be given by $d\alpha|_F(B,\bullet)=\alpha|_F$.  One then has 
		\begin{equation}\label{eq:local_dynamics}  \phi^M_{2k}(F,p)= \frac{\tr \wedge^{2n-2k} d_pB}{|\det d_pB|} \end{equation}
		
		\item 	From a geometric point of view, 
		\begin{equation}\label{eq:local_geometry} \phi^M_{k}(F,p)= {2n \choose k} |\det (S-h)|^{-1}D(S-h[2n-k], J[k]).  
		\end{equation}
		Here $D$ denotes the mixed discriminant, 
		$J=\left(\begin{array}{cc}0 & -I_n\\ I_n & 0\end{array}\right)$, $S$ is the second fundamental form of $F$ at $p$, and $h$ the second fundamental form of the contact distribution at $p$ (see Definition \ref{def:second_fund_form}), both written with respect to a frame compatible with the contact structure in a natural way at $p$, see Definition \ref{def:compatible_riemannian}.
		Thus they are reminiscent of (certain symmetric functions of) the principal radii of an embedded hypersurface in a Riemannian manifolds. This point of view is applicable in the wider setting of DH manifolds, as described below.
	\end{itemize}
	
We may extend $\phi_{2k}$ to a non-negative lower semicontinuos functional on all $2k$-dimensional submanifolds with boundary, denoted $\CA_{2k}(F)$, the contact area of $F$. We observe that $\CA_{2k}(F)=0$ if and only if $F$ can be made nowhere tangent to the contact distribution by an arbitrarily small perturbation.
\newline\newline
In Riemannian or Hermitian manifolds, a fair amount of the valuation theory appears already in the corresponding flat space, which can be thought of as the tangent space to the given manifold. In the contact setting, it is no longer true: the tangent space of a contact manifold does not itself inherit a contact structure. 

The main observation guiding this paper is that every cotangent space of a contact manifold is canonically the Heisenberg Lie algebra. We are thus led to study the valuation theory of general manifolds with such structure.
\begin{Definition}
	A manifold $X$ equipped with a hyperplane distribution $H$ (called horizontal) and a smooth field of nowhere-degenerate forms $\omega\in \Gamma^\infty(X, \wedge ^2 H_x^*\otimes (T_xX/H_x))$ will be called a dual Heisenberg (DH) manifold. 

\end{Definition}

The space of valuations naturally associated to such manifolds turns out to resemble somewhat the Lipschitz-Killing space of valuations in Riemannian geometry, in particular, they exhibit universality with respect to embeddings.

Theorem \ref{mainthm:contact} is then the contact instance of the following general result.

\begin{Theorem} \label{mainthm:DH}
	To any DH manifold $M^{2n+1}$ with horizontal distribution $H\subset TM$ there are canonically associated, generalized valuations $\phi^M_{k}\in\mathcal W_k^{-\infty}(M)$, $0\leq k\leq 2n$. They have the following properties:
\begin{enumerate}
	\item $\phi^M_0$ is the Euler characteristic, and $\phi^M_{2k}\in \mathcal W_{2k}^{-\infty}(M)\setminus \mathcal W_{2k+1}^{-\infty}(M)$.
	
	\item $\phi_{k}^M$ can be naturally evaluated on submanifolds in generic position with respect to the horizontal distribution. For a generic closed hypersurface $F$, 
	
	\[ \phi^M_{k}(F)=\sum_{T_pF=H_p}  \phi^M_{k}(F,p)  \]
	where $\phi^M_{k}(F,p)$ only depends on the germ of $F$ at $p$. 
	\item If $i:N^{2m+1}\to M^{2n+1}$ is a DH embedding, then for $k\leq 2m$ one has $i^*\phi^M_{k}=\phi^N_{k}$.
\end{enumerate}
\end{Theorem}
The local contact areas $\phi_k^M(F,p)$ are given by the same curvature-type formula \eqref{eq:local_geometry} as in the contact case.

As an intermediate step of independent interest, we obtain a Hadwiger-type theorem for the dual Heisenberg algebra itself. We denote by $U=\R^{2n+1}$ the dual of the Heisenberg Lie algebra $\mathfrak h_{2n+1}$, and by $\Sp_H(U)$ its automorphism group.

\begin{Theorem}\label{mainthm:DH_Hadwiger}
It holds that $\Val^{-\infty}(U)^{\Sp_H(U)}$ consists of even valuations. For $0\leq k\leq n$, $\Val^{-\infty}_{2k+1}(U)^{\Sp_H(U)}=\{0\}$ while $\dim \Val^{-\infty}_{2k}(U)^{\Sp_H(U)}=1$. The corresponding Klain sections are zero-order distributions (that is, regular Borel measures).
\end{Theorem}
We also consider $\Sp_H^+(U)$, the connected component of the identity. We prove

\begin{Theorem}\label{mainthm:DH+_Hadwiger}
For $0\leq k\leq n$, $\Val^{-\infty}_{2k+1}(U)^{\Sp_H^+(U)}=\{0\}$ while $\dim \Val^{-\infty}_{2k}(U)^{\Sp^+_H(U)}=2$.
In the latter space, the $\Sp_H(U)$-invariant valuations are complemented by a one-dimensional space of odd valuations.
\end{Theorem}

We then consider the standard contact sphere $S^{2n+1}$, which we identify with the oriented projectivization $\mathbb P_+(V)$ of a symplectic space $V=(\R^{2n+2},\omega)$.  We construct a canonical distribution $\mu_\omega\in \mathcal M^{-\infty}(\Gr_{2k}(V))^{\Sp(V)}$. We obtain the following Crofton formulas, establishing further common grounds with the Riemannian setting. 

\begin{Theorem}\label{mainthm:crofton_sphere}
	Define for $0\leq i\leq n$ the generalized valuations $\psi_{2i}$ on $S^{2n+1}$ given by the Crofton formula 
	\[\psi_{2i}:=\int_{\Gr_{2n+2-2i}(V)}\chi(\bullet\cap E)d\mu_{\omega}(E).\] 
	Then for certain explicit constans $c^n_{ij}$ one has
	
	\begin{equation*} \psi_{2i}=\sum_{j=i}^n c^n_{ij}\phi_{2j}.\end{equation*}
\end{Theorem}

We also establish a Hadwiger-type theorem for the contact sphere.

\begin{Theorem}\label{mainthm:sphere_Hadwiger}
Both $(\phi_{2k})_{k=0}^n$ and $(\psi_{2k})_{k=0}^n$ are bases of $ \mathcal V^{-\infty}(S^{2n+1})^{\Sp(2n+2)}$. 
\end{Theorem}

Finally, we find that while symplectic space and manifolds do not possess interesting invariant valuations, one can nevertheless write certain integral geometric formulas for symplectic volumes of manifolds.
We construct a canonical distribution on the affine oriented Grassmannian $\overline \mu_\omega\in \mathcal M^{-\infty}(\AGr^+_{2k}(\R^{2n}))$ which is translation- and $\Sp(2n)$-invariant, and odd to orientation reversal. We prove the following Crofton formula on symplectic linear space.

\begin{Theorem}\label{mainthm:symplectic} Let $F^{2k}\subset \R^{2n}$ be a $C^1$ compact, oriented submanifold with boundary. Then
	\[  \int_F \omega^k=(-1)^\kappa{n\choose k}{2n \choose 2k}^{-1}\frac{(2n-1)^\kappa}{2^{\kappa+1} } \int_{\AGr^+_{2n-2k}(\R^{2n})} I(E, F)d\overline\mu_\omega(E) \]
	where $\kappa=\min(k, n-k)$ and $I$ is the oriented intersection index.
	
\end{Theorem}

\subsection{Plan of the paper}

In section \ref{sec:basics} we introduce notation and present the basic geometric facts we will use. In section \ref{sec:valuations} we recall the basics of valuation theory and prove some lemmas we will need. In section \ref{sec:DH_manifolds} we construct the canonical valuations on general DH manifolds and establish their universality to embeddings, proving Theorem \ref{mainthm:DH}. We also explore some geometric properties of those valuations.
In section \ref{sec:flat_DH} we classify the translation-invariant valuations of the dual Heisenberg algebra, proving Theorems \ref{mainthm:DH_Hadwiger} and \ref{mainthm:DH+_Hadwiger}, and note their relation to gaussian curvature. Apart from its intrinsic interest, the linear classification is needed for the uniqueness statement in Theorem \ref{mainthm:contact}, as well as for Theorem \ref{mainthm:sphere_Hadwiger}. In section \ref{sec:contact} we specialize the DH valuations to contact manifolds and prove Theorem \ref{mainthm:contact}. In particular, we give the dynamical description of the contact valuations. In section \ref{sec:grassmannians} we construct symplectic-invariant distributions on linear and affine Grassmannians in symplectic space, which are used in the subsequent two sections. In section \ref{sec:contact_sphere} we consider the standard contact sphere. We produce Crofton formulas for $\phi_{2k}$ that are invariant under $\Sp(V)$, proving theorems \ref{mainthm:crofton_sphere} and \ref{mainthm:sphere_Hadwiger}, and compute some examples explicitly. We also bound from below the contact valuations of a convex set. Finally in section \ref{sec:symplectic}, we study the integral geometry of linear symplectic space, proving Theorem \ref{mainthm:symplectic}.
\subsection{Acknowledgements}
I am grateful to Semyon Alesker, Andreas Bernig, Yael Karshon, Leonid Polterovich, Egor Shelukhin and Gil Solanes for fruitful discussions and insightful suggestions. I am also indebted to Andreas and Semyon for helpful comments on the first draft.

This research was partially supported by an NSERC Discovery Grant, and conducted while the author was employed as Coxeter Assistant Professor at the University of Toronto.

\section{Preliminaries}\label{sec:basics}
\subsection{Notation}
We use $\sigma_1$ to denote the unit $\SO(n)$-invariant measure on various homogeneous spaces of the special orthogonal group.
We write $\Gr_k(V)$ for the $k$-Grassmannian in $V$, $\Gr_k^+$ for the oriented Grassmannian, and $\AGr^{(+)}$ for the affine (oriented) Grassmannian. $\mathcal K(V)$ is the set of compact convex subsets of $V$.

The one-dimensional space of real-valued Lebesgue measures over $V$ is denoted $\Dens(V)$. For a manifold $M$, $|\omega_M|$ denotes the line bundle of densities, whose fiber over $x\in M$ is $\Dens(T_xM)$. 
We will write $M^{tr}$ for the translation-invariant elements of a module $M$ over $V$.
For a group $G\subset \textrm{GL}(V)$, $\overline G$ is the group generated by $G$ and all translations in $V$.

We will write $\Omega_{-\infty}(M)$ for the space of currents on $M$, since we typically consider them as generalized differential forms.

Throughout the note, $J=\left(\begin{array}{cc}0 & -I_n\\ I_n & 0\end{array}\right)$ is the standard real form of $\sqrt{-1}$.

\subsection{The symplectic group action on the Grassmannian}

This subsection is the symplectic version of the corresponding section in \cite{bernig_faifman} where $\OO(p,q)$ is considered.

Let $V=(\R^{2n}, \omega)$ be a symplectic space. Let $X^{k}_r(V)$, $0\leq r\leq \kappa:=\lfloor \min(k, 2n-k)/2\rfloor$ be the orbits of $\Gr_{k}(V)$ under the real symplectic group $\Sp(V)$, where \begin{equation}
\label{eq:orbits}X^{k}_r(V)=\{E\in \Gr_{k}(V): \dim\Ker \omega|_E=\min(k, 2n-k)-2r\}.
\end{equation} 
When no confusion can arise we write simply $X^{k}_r$.
In the oriented Grassmannian $\Gr_{k}^+(\R^{2n})$, the double covers $X^{k, +}_r$ of $X^{k}_r$ are orbits of $\Sp(V)$ when $0\leq r\leq \kappa-1$. For $r=\kappa$, there are two possibilities: for even $k$, the double cover of $X^{k, +}_\kappa$  splits into two open orbits, denoted $X^{k}_\pm$, corresponding to the different orientations induced by $\omega$ on the subspace; for odd $k$, $X^{k, +}_\kappa$ is a single orbit.

We use the same notation for the corresponding $\overline{\Sp(V)}$-orbits in the (oriented) affine Grassmannian. 

We will need a simple lemma from linear algebra, which is a linearized version of Witt's theorem.
\begin{Lemma}\label{lem:isotropic_tangent} Take $E\in X^k_r$. Write $E_0=E\cap E^\omega$, $\dim E_0=\min(k,2n-k)-2r$. Let $\pi_0: E^*\otimes V/E\to E_0^*\otimes E_0^*$ denote the map $\pi_0(T)(x,y)=\omega(Tx,y)$. Then
	\[T_EX^{k}_{r}=\{T:E\to V/E : \pi_0(T)\in \Sym^2 E_0^*\}\]
\end{Lemma} 
\proof
Clearly $\pi_0$ is onto. Denote by $\pi_E: V^*\otimes V\to E^*\otimes (V/E)$ the natural projection. Recall that $\mathfrak{sp}(V)=\{T\in\mathfrak {gl}(V): \omega(Tx, y)=\omega(Ty, x)\quad\forall x,y\in V\}$.
We have to show that $\pi_{E} \mathfrak {sp}(V)=\pi_0^{-1}\Sym^2E_0^*$.

The inclusion $\pi_{E} \mathfrak {sp}(V)\subset\pi_0^{-1}\Sym^2E_0^*$ is immediate. In the other direction, take $T:E\to V/E$ s.t. $\pi_0T\in \Sym^2E_0^*$. We should lift $T$ to $\tilde T\in \mathfrak {sp}(V)$.

Choose any subspaces $E' \subset E$, $E'' \subset E^\omega$ with $E=E_0 \oplus E'$, $E^\omega=E_0 \oplus E''$.  

Denote $\kappa=\min(k,2n-k)$. Then $\dim E'=k-\kappa+2r$ and $\dim E''=2n-k-\kappa+2r$, and they are both non-degenerate. Moreover, $E' \oplus E'' \subset V$ is a non-degenerate subspace of dimension $2n-2\kappa+4r$, and $E_0$ is an isotropic subspace of the non-degenerate space $W:=(E' \oplus E'')^\omega$, and $\dim W=2\kappa-4r$. That is, $E_0\subset W$ is a Lagrangian subspace. Fix a Lagrangian complementing space $F\subset W$ s.t. $W=E_0 \oplus F$. Then $V=E' \oplus E'' \oplus E_0 \oplus F$.  

Let $T_1 \in \Hom(E,V)$ be a lift of $T$ such that $T_1(E') \subset E'' \oplus F$, $T_1(E_0) \subset E'' \oplus F \oplus E'$ and 
\begin{displaymath}
\omega( \pi'(T_1x),e):=\omega(Te,x), \quad \forall x \in E_0, e \in E',
\end{displaymath}
where $ \pi'(T_1x)$ is the $E'$-component of $T_1x$.

Note that $\omega$ gives identifications $E^*=E' \oplus F$, as well as $(E \oplus E'')^*=E' \oplus E'' \oplus F$. We extend $T_1$ to a map $T_2 \in \Hom(E \oplus E'',V)$ by requiring $T_2(E'') \subset E' \oplus F$ and 
\begin{displaymath}
\omega(T_2e'',e):=\omega(T_1e,e''), \quad \forall e'' \in E'', e \in E.
\end{displaymath}

Note that $\omega(T_2e_1,e_2)=0$ for all $e_1,e_2 \in  E''$. 

Finally, we extend $T_2$ to a map $T' \in \Hom(V,V)$ by requiring $T'(F) \subset E' \oplus E'' \oplus F$ and  
\begin{displaymath}
\omega(T'f,x)=\omega(T_2x,f), \quad \forall x \in E \oplus E'', f \in F.
\end{displaymath}
Again we have $\omega(T'f,f')=0$ for $f,f' \in F$. Then $T' \in \mathfrak{sp}(V)$ lifts $T$ as required.
\endproof

\begin{Corollary}\label{cor:normal_space}
	There is a natural identification of $N_EX^{k}_r$ with $\wedge^2 E_0^*$.
\end{Corollary}

\proof\begin{align*}N_EX^{k}_r&=T_E\Gr_{k}(V)/T_EX^{k}_r=E^*\otimes (V/E)/ \pi_0^{-1}\Sym^2E_0^* \\&=(E_0^*\otimes E_0^*)/\Sym^2E_0^*=\wedge^2 E_0^*\end{align*}
\endproof

\begin{Lemma}\label{lem:transversal_intersection}
	Fix $L\in \mathbb P_+(V)$, and denote $Z_L=\{E\in\Gr_{k}(V): L\subset E\}$. Then $Z_L$ intersects $X^{k}_{r}$ transversally for all $r$.
\end{Lemma}
\proof
Take $E\in Z_L\cap X^{k}_{r}$, denote $E_0=E\cap E^\omega$. Let $A:E\to V/E$ be an arbitrary linear map. We will find a representation $A=A_1+A_2$ with $A_1\in T_EZ_L$ (that is, $A_1|_L=0$) and $A_2\in T_E X^{k}_r$.  Fix $l\in L$ and write $Al=v+E$.
\\\textit{Case 1}: $L\not\subset E_0$. We set $A_2|_L:=A|_L$, $A_2|_{E_0}=0$.
\\ \textit{Case 2}: $L\subset E_0$ and $\omega(v, E_0)= 0$. Decompose $E_0=L\oplus F$ and set $A_2(l)=v+E$, $A_2(f)=0$ for all $f\in F$.
\\\textit{Case 3}: $L\subset E_0$ and $\omega(L, AL)\neq 0$. Since $\omega(l, v)\neq 0$, we may decompose $E_0=E_0\cap v^\omega\oplus L$. We then set $A_2(l)=v+E$ and $A_2(x)=0$ for $x\in E_0\cap v^\omega$.
\\\textit{Case 4}: $L\subset E_0$, $\omega(v, E_0)\neq 0$ and $\omega(L, AL)= 0$. Then we decompose $E_0=L\oplus F\oplus \Span(u)$ where $L\oplus F=v^\omega\cap E_0$ and $\omega(u, v)=1$. Choose $w\in V$ such that $\omega(w, l)=-1$ and $w\in F^\omega$. This is possible since $ l^\omega\neq F^\omega$. Then set $A_2l=v+E$, $A_2f=0$ for $f\in F$, and $A_2u=w+E$.
\\In all cases, extend $A_2$ arbitrarily to $E$. Thus in all cases, by Lemma \ref{lem:isotropic_tangent}, $A_2\in T_E X^{k}_{r}$, and $A_1:=A-A_2\in T_EZ_L$.
\endproof

\section{Valuation theory}\label{sec:valuations}

\subsection{Valuations on manifolds}
For a manifold $X$, we let $\mathbb P_X:=\mathbb P_+(T^*X)$ denote the oriented projectivization of its cotangent bundle and $\pi:\mathbb P_X\to X$ the projection. $\mathbb P_X$ has a canonical contact structure. A form $\omega\in\Omega(\mathbb P_X)$ that vanishes when restricted to the contact distribution is usually called vertical. However, we will have several different notions of verticality, so we will call such forms Legendrian.

\begin{Definition}
We say that a form $\omega\in\Omega^d(\mathbb P_X)$ has horizontal degree at least $k$, written $\deg_H\omega\geq k$, if $\omega(v_1,\dots, v_d)$ vanishes whenever $d+1-k$ of the vectors $v_j$ are vertical, that is tangent to the fiber of $\pi:\mathbb P_X\to X$. 
\end{Definition}The following is a simple reformulation.

\begin{Lemma}\label{lem:horizontal_degree_wedge}
$\{\omega \in \Omega(\mathbb P_X): \deg_H\omega\geq k\}$ is the ideal in $\Omega(\mathbb P_X)$ generated by $\pi^*\Omega^k(X)$.
\end{Lemma}
\proof

Assume $\omega=\sum \pi^*\omega_j \wedge \eta_j$, where $\omega_j\in \Omega^k(X)$, $\eta_j\in\Omega^{d-k}(\mathbb P_X)$. Clearly if $v_1,\dots,v_{d+1-k}$ are vertical vectors then $\pi^*\omega_j\wedge \eta_j(v_1,\dots, v_{d+1-k},\dots)=0$.

For the opposite direction, let us choose a Riemannian structure on $X$. Then $\mathbb P_X$ is the sphere bundle on $X$. Fix coordinates $dx_j$ on $T_xX$ and $d\xi_j$ on $T_{\xi}S_xX$. 
Then $\pi^*\Omega^k(X)$ is spanned over $C^\infty(\mathbb P _X)$ by $\{\wedge_{i\in I} \pi^* dx_i: |I|=k\}$. Assume $\deg_H\omega\geq k$ and decompose $\omega=\sum_{I, J}  f_{IJ}\pi^*dx_I\wedge d\xi_J$. 

Assume a multi-index $I$ appears in the sum  with $|I|<k$, say $I=(i_1,\dots, i_l)$, $l<k$, with corresponding $J=(j_1,\dots,j_{d-l})$. Let $e_1 ^H,\dots, e_n^H, e_1^V, \dots, e_{n-1}^V$ be dual to $\pi^*dx_1, \dots, \pi^*dx_n, d\xi_1,\dots, d\xi_{n-1}$. Then  $d-l\geq d+1-k$ so by assumption
$ 0=\omega(e_{i_1}^H, \dots, e_{i_l}^H, e_{j_1}^V, \dots, e_{j_{d-l}}^V )=f_{IJ}(x,\xi)$, so $f_{IJ}=0$. 

\endproof

Let $M$ be a smooth manifold, which we assume oriented for simplicity of exposition, and refer the reader to \cite{alesker2, alesker3} for the general case. Denote by $\mathcal P(M)$ the compact differentiable polyhedra of $M$. We remark that manifolds with corners are an example of differentiable polyhedra, and refer to \cite{alesker2} for the definition of differentiable polyhedra.
The smooth valuations $\mathcal V^\infty(M)$ consist of functionals $\phi:\mathcal P(M)\to\R$ which can be presented in the form $\phi(X)=\int_X\mu+\int_{N^*X}\omega$ for some forms $\mu\in\Omega^n(M)$ and $\omega\in\Omega^{n-1}(\mathbb P_M)$. Here $N^*X$ is the conormal cycle of $X$. It consists of codirections $\xi\in \mathbb P_+(T^*_xM)$ which are non-positive on velocity vectors $\dot{\gamma}(0)\in T_xM$ of all curves $\gamma(t)\in X$ with $x=\gamma(0)$. The Euler characteristic $\chi$ is an important example of a smooth valuation. The smooth valuations over open subsets of $X$ constitute a soft sheaf over $X$, see \cite{alesker2}. We denote the compactly supported valuations by $\mathcal V_c^\infty(M)$, and $\mathcal W_{i,c}^\infty(M):=\mathcal V_c^\infty(M)\cap \mathcal W_i^\infty(M)$. 
We remark that whenever a valuation is evaluated on a subset, the subset is assumed compact.
There is a natural integration functional $\int_M:\mathcal V^\infty_c(M)\to \R$ which is essentially evaluation on $M$, but more precisely - on a sufficiently large compact set. Both $\mathcal V^\infty(M)$ and $\mathcal V^{\infty}_c(M)$ inherit natural topologies from the corresponding spaces of pairs of forms. 

There is a natural filtration $\mathcal W^\infty_n(M)\subset \dots\subset \mathcal W^\infty_0(M)=\mathcal V^\infty(M)$, introduced by Alesker \cite{alesker2}. We will use an equivalent description, which is the content of \cite[Corollary 3.1.10]{alesker2}. 

\begin{Definition}\label{def:filtration_form}
	$\mathcal W^\infty_k(M)$ consists of those valuations that can be represented by a pair $(\omega, \mu)$ with $\deg_H \omega \geq k$.
\end{Definition}
In particular, $\mathcal W_n^\infty(M)$ are just the smooth measures on $M$, denoted $\mathcal M^\infty(M)$.

Alesker defined a product structure $\mathcal W_i^\infty(M)\otimes \mathcal W_j^\infty(M)\to \mathcal W_{i+j}^\infty(M)$ which turns $\mathcal V^\infty(M)$ into a filtered algebra, whose unit is the Euler characteristic. It induces the following Alesker-Poincar\'e duality.
\begin{Theorem}[Alesker \cite{alesker4}]
	The pairing $\mathcal W^\infty_i(M)\otimes \mathcal W_{n-i,c}^\infty(M)\to \R$ given by
	$(\phi,\psi)\mapsto (\phi\cdot \psi)(M)$ is non-degenerate.
\end{Theorem}

The presentation of $\phi\in\mathcal V^\infty(M)$ by a pair of forms is not unique. There is an alternative faithful description due to Bernig and Br\"ocker \cite{bernig_brocker}. In the following, $a:\mathbb P_M\to\mathbb P_M$ is the antipodal map in every fiber, and $D$ is the Rumin differential introduced in \cite{Rumin}. We recall that $D\omega$ is the unique Legendrian form $d(\omega+\eta)$, where $\eta\in\Omega^{n-1}(\mathbb P_M)$ ranges over all Legendrian forms.
\begin{Theorem}\label{thm:current_injective}
	If $\phi$ is represented by the pair $(\omega, \mu)\in\Omega^{n-1}(\mathbb P_M)\times \Omega^n(M)$, then \begin{equation}\label{eq:currents_formula}(T,C):=(a^*(D\omega+\pi^*\mu), \pi_*\omega)\in \Omega^n(\mathbb P_M)\times C^\infty(M)\end{equation} is determined by $\phi$. They satisfy the relations $dT=0$, $\pi_*T=(-1)^ndC$, and $T$ is Legendrian. Moreover, any $(T,C)$ with those properties corresponds to a valuation.
\end{Theorem}

We refer to $(T,C)$ as the defining currents of $\phi$ (and often we refer just to $T$ as the defining current). The reason for this terminology will become evident once we introduce generalized valuations. 

The Alesker-Poincar\'e duality is easy to describe using the defining currents.
\begin{Theorem}[Bernig \cite{bernig_product}]
	Let $(\omega, \mu)$ represent $\phi_1$, and let $\phi_2$ have defining current $(T_2, C_2)$. Then
	\begin{equation}\label{eq:bernig_product}\int_M \phi_1\cdot \phi_2=\int_{\mathbb P_M}\omega_1\wedge T_2+\int_M C_2\mu_1.\end{equation}
\end{Theorem}

Let us describe the filtration on $\mathcal V^\infty(M)$ through the defining currents. 
\begin{Lemma}
	Let $(T, C)\in \Omega^n(\mathbb P_M)\times C^\infty(X)$ be the defining current of $\phi\in\mathcal V^\infty(M^n)$. 
	For $1\leq k\leq n$, $\phi\in \mathcal W_k^\infty(M^n)$ if and only if $C=0$ and $\deg_H T\geq k$.
\end{Lemma}
\proof
Consider first $k=n$. If $\phi\in \mathcal W_n^\infty(M)$, we have $C=0$ and $T=\pi^*\mu$ for some $\mu\in\Omega^n(M)$, which is clearly horizontal of degree at least $n$.
In the other direction, if $\deg_H T\geq n$, it follows that $T=f(\xi)\pi^*\mu$ for some $\mu\in \Omega^n(M)$ and $f\in C^\infty (\mathbb P_M)$. Now $0=dT=df\wedge \pi^*\mu$, that is $df$ vanishes when restricted to the vertical fiber, so $f=\pi^*f_1$ for some $f_1\in C^\infty(M)$. Hence $T=\pi^*(f_1\mu)$, so $\phi$ is just the measure $f_1\mu$.

Assume now $1\leq k\leq n-1$. Recall that by Alesker-Poincar\'e duality, $\phi\in\mathcal W_k^\infty(M)$ if and only if for all $\psi\in\mathcal W_{n+1-k,c}^\infty(M)$, $\int_M\phi\cdot \psi=0$. Combining Definition \ref{def:filtration_form} with eq. \eqref{eq:bernig_product}, we conclude that $\phi\in\mathcal W_k^\infty(M)$ if and only if for all $\omega\in\Omega^{n-1}(\mathbb P_M)$ with $\deg_H\omega\geq n+1-k$, and all $\mu\in\Omega^n(M)$, it holds that $\int_M(\pi_*(\omega\wedge T)+C\mu)=0$. 

Assume that $C=0$ and $\deg_H T\geq k$. Then by Lemma \ref{lem:horizontal_degree_wedge}, $\omega\wedge T=0$ for all $\omega$ as above, and hence $\phi\in\mathcal W_k^\infty(M)$.

In the other direction, assume $\phi\in\mathcal W_k^\infty(M)$. Taking $\omega=0$ and $\mu$ arbitrary we deduce $C=0$. It then follows that $\pi_*(\omega\wedge T)=0$, and since $\omega$ can have an arbitrarily small support, also that $\omega\wedge T=0$ whenever $\deg_H\omega\geq n-k+1$, so that $\deg_H T\geq k$.

\endproof

The valuations that appear naturally in contact manifolds are not smooth. To formally study them we will need the larger family of generalized valuations.

\begin{Definition}
The generalized valuations  $\mathcal V^{-\infty}(M)$ are the continuous dual of $\mathcal V_c^\infty(M)$, equipped with the weak topology.
\end{Definition}  
A generalized valuation is uniquely determined by its defining current $(T, C)\in\ \mathcal D_{n-1}(\mathbb P_M)\times \mathcal D_n(M)$, see \cite{alesker_bernig_product}. It can be an arbitrary pair of currents satisfying the three properties: $T$ is Legendrian, $\pi_*T=\partial C$, and $\partial T=0$. If $\phi_2\in\mathcal V^{-\infty}(M)$ has defining current $(T_2, C_2)$, it acts on smooth valuations $\phi_1\in\mathcal V_c^\infty(M)$ represented by forms $(\omega_1, \mu_1)$ through eq. \eqref{eq:bernig_product}. We will write $\phi_2=[T_2, C_2]$, $T_2=T(\phi_2)$, $C_2=C(\phi_2)$.

\begin{Example}
	1. Given $X\in\mathcal P(M)$, the evaluation at $X$ functional $\chi_X:\mathcal V^\infty(M)\to \R$ is a generalized valuation with defining currents $C=[[X]]$ and $T=[[N^*X]]$. 
	\\2. A generalized valuation $\phi_1$ can be represented by a pair of generalized forms $\omega_1\in\Omega_{-\infty}^{n-1}(\mathbb P_M)$ and $\mu_1\in \Omega_{-\infty}^{n}(M)$. It acts on  $\phi_2=[T_2, C_2]\in\mathcal V^\infty_c(M)$ through eq. \eqref{eq:bernig_product}. It has defining currents given by $T_1=a^*(D\omega_1+\pi^*\mu_1)$, $C_1=\pi_*\omega_1$.
\end{Example}

The filtration $\mathcal W_j^\infty(M)$ on $\mathcal V^{\infty}(M)$ extends to a filtration  $\mathcal W_j^{-\infty}(M)$ on $\mathcal V^{-\infty}(M)$ by taking $\mathcal W_j^{-\infty}(M)$ to be the annihilator of $\mathcal W_{n+1-j, c}^{-\infty}(M)$.

We refer to \cite{hormander, guillemin_sternberg} for the notion of the wavefront set of a distribution. Let us only record that for an oriented submanifold $X\subset M$, $\WF([[X]])=N^*X$.
\begin{Definition}
	The wavefront $\WF(\phi)$ of $\phi\in\mathcal V^{-\infty}(M)$ is the pair of wavefront sets $(\WF(T(\phi)), \WF(C(\phi)))$. When $\WF(C(\phi))=\emptyset$, we also write $\WF(\phi)=\WF(T(\phi))$.
\end{Definition}
Given a closed cone $\Gamma\subset T^*(\mathbb P_M)\setminus 0$, we let $\mathcal V^{-\infty}_\Gamma(M):=\{\phi\in  \mathcal V^{-\infty}(M):\WF(T(\phi))\subset\Gamma, \WF(C(\phi))=\emptyset \}$. It inherits a topology from H\"ormander's topology on the corresponding space of currents $\mathcal D_{n-1, \Gamma}(\mathbb P_M)$.

The Alesker-Poincar\'e duality extends to a pairing of generalized valuations, as long as the wavefront sets are in good relative position. We will only need the following weak version of Theorem 8.3 in \cite{alesker_bernig_product}.

For a subset $T\in T^*\mathbb P_M$, we write $-T$ for its image under the antipodal map in every cotangent space $T^*_{x,\xi} T^*(\mathbb P_M)$. Define also $T^s:=T\cup(-T)\cup aT\cup (-aT)$.

\begin{Theorem}
	Let $\Gamma_1,\Gamma_2\subset T^*(\mathbb P_M)\setminus 0$ be closed cones, $\Gamma_j=\Gamma_j^s$. Assume that $\Gamma_1\cap \Gamma_2=\emptyset$. Then there is a jointly sequentially continuous pairing
	\[ \mathcal V^{-\infty}_{\Gamma_1}(M)\otimes \mathcal V^{-\infty}_{\Gamma_2}(M)\to \R  \] 
	extending the Alesker-Poincar\'e duality on the dense subspaces of smooth valuations.
\end{Theorem}

It follows that if $\phi=[T,C]$ with $C$ smooth, and $X\in\mathcal P(M)$ s.t. $\WF[[N^*X]]\cap \WF(T^s)=\emptyset$, then one can naturally evaluate $\phi(X):=\langle \phi, \chi_X\rangle$. In particular, if $X, Y\in\mathcal P(M)$ have disjoint conormal cycles, one can evaluate $\int_M\chi_X\cdot \chi_Y$. In all reasonable settings, the result should equal $\chi(X\cap Y)$. This was shown to be the case when $X,Y$ are transversal submanifolds with corners in \cite[Theorem 6]{alesker_bernig_product}.

The Euler-Verdier involution $\sigma:\mathcal V^{-\infty}(M)\to \mathcal V^{-\infty}(M)$ was introduced by Alesker in \cite{alesker2, alesker4}. It can be described through its action on the defining currents: $\sigma [T, C]=[(-1)^na^*T, C]$.

We will need the following simple lemma, which we later use in the proof of Proposition \ref{prop:contact_linear_dependencies}.
\begin{Lemma}\label{lem:hypersurface_suffices}
	Take a generalized valuation $\phi\in\mathcal V^{-\infty}(M^{n})$ satisfying $\sigma\phi=(-1)^{n+1}\phi$ and $C(\phi)$ smooth. 
	Assume that $\phi(F)=0$ for all closed hypersurfaces $F$ for which $N^*(N^*F)$ is disjoint from $\WF(\phi)$. 
	Assume moreover that $M$ can be covered by open charts $U_\alpha\cong \R^n$, such that if $K\subset U_\alpha$ is a smooth, strictly convex body, then $N^*(N^*\partial K)\cap\WF(\phi)=\emptyset$. Then $\phi=0$.
\end{Lemma}
\proof
As generalized valuations form a sheaf over $M$ \cite[Proposition 7.2.2]{alesker4}, we may assume $M=U_\alpha=\R^{n}$.

By \cite{alesker1}, valuations of the form $\mu(\bullet - K)$, where $K$ is a convex body with smooth support function and $\mu\in\mathcal M_c^\infty(\R^n)$, span a dense subspace in $\mathcal V_c^\infty(\R^n)$. It remains dense if we only take strictly convex $K$. Fix such $\mu$ and $K$.

By \cite[Lemma 4.1.1]{alesker4} and approximating $\phi$ by smooth valuations in the H\"ormander topology with the same Euler-Verdier eigenvalue, we have 
$(-1)^{n+1}\phi(K)=\sigma\phi(K)=(-1)^n(\phi(K)-\phi(\partial K))$. By assumption, $\phi(\partial K)=0$, and hence $\phi(K)=0$. It follows that
\[ \langle \phi, \mu(\bullet - K)\rangle = \langle \phi, \int_{\R^n}\chi(\bullet \cap (K+x))d\mu(x)\rangle= \int_{\R^n} \phi(K+x)d\mu(x)=0 \]
concluding the proof.
\endproof

\subsection{Translation invariant valuations}
We will need the following standard facts when we study valuations in the dual Heisenberg algebra in section \ref{sec:flat_DH}.

The space of smooth translation invariant valuations $\Val^\infty(\R^n)$ can be defined as $\mathcal V^{\infty}(\R^n)^{tr}$. It is a Fr\'echet space. We will also consider generalized translation-invariant valuations, $\Val^{-\infty}(\R^n):=\mathcal V^{-\infty}(\R^n)^{tr}$.

The even/odd valuations are $\Val^{\pm, \pm\infty}$ which have eigenvalue $\pm1$ under the antipodal map. The $k$-homogeneous valuations are those $\phi\in  \Val^{\pm\infty}$ that have eigenvalue $\lambda^k$ under all rescalings of $\R^n$ by $\lambda>0$.

It is well-known that $\Val_0^{-\infty}(\R^n)=\Span\{\chi\}$, and by Hadwiger's theorem \cite{hadwiger}, $\Val_n^{-\infty}(\R^n)=\Span\{\vol_n\}$. It follows from McMullen's work \cite{mcmullen} that $\Val^{\pm\infty}(\R^n)=\oplus_{k=0}^n\Val_{k}^{\pm\infty}(\R^n)$.

The line bundle over $\Gr_k(\R^n)$ whose fiber over $E$ is $\Dens(E)$ is called the Klain bundle.
The Klain map $\Kl: \Val^{+,\infty}_k\to \Gamma^{\infty}(\Gr_k(V), \Dens(E))$ is given by $\Kl(\phi)(E)=\phi|_E$. It is well-defined by Hadwiger's theorem, and injective by a theorem of Klain \cite{klain}. By \cite{alesker_faifman}, it admits an injective extension $\Kl: \Val^{+,-\infty}_k\to \Gamma^{-\infty}(\Gr_k(V), \Dens(E))$.  

Dual to the Klain section is the Crofton map, given by $\Cr:\mathcal M^\infty(\AGr_{n-k}(V))^{tr}\to \Val^{+,\infty}_k(V)$, given by $\Cr(\mu)(K)=\int_{\AGr_{n-k}(V)}\chi(E\cap K)d\mu(E)$. It follows from Alesker's irreducibility theorem \cite{alesker_mcmullen} that $\Cr$ is surjective. It was shown in \cite{alesker_faifman} that it admits a surjective extension $\Cr:\mathcal M^{-\infty}(\AGr_{n-k}(V))^{tr}\to \Val^{+,-\infty}_k(V)$

The relation to valuations on manifolds is as follows:
\begin{Theorem}[\cite{alesker2}]\label{thm:filtration_quotients}
	$\mathcal W_k^{\infty}(M)/\mathcal W_{k+1}^{\infty}(M)$ is naturally isomorphic to $\Gamma^\infty(M, \Val_k^{\infty}(TM))$, the space of smooth sections of the bundle with fiber $\Val_k^{ \infty}(T_xM)$ over $x\in M$. 
	\end{Theorem}
Consequently by the Alesker-Poincar\'e duality we have \[\mathcal W_k^{-\infty}(M)/\mathcal W_{k+1}^{-\infty}(M)\simeq (\mathcal W_{n-k}^{\infty}(M)/\mathcal W_{n-k+1}^{\infty}(M))^* \simeq \Gamma^\infty(M, \Val_{n-k}^\infty(TM))^*.\]

\subsection{Crofton formulas on manifolds}  \label{sec:crofton}

Consider a double fibration
\[\xymatrix{& W\ar[dl]_{\tau} \ar[dr]^{\pi}&\\ Z & & X }\] where $\pi, \tau$ are proper, and $\tau\times \pi: W\to Z\times X$ is a closed embedding. We view $Z$ as a space parameterizing a family of submanifolds $\hat z:=\pi\tau^{-1}z\subset X$, of dimension $\dim \hat z=\dim W-\dim Z$.
Assume furthermore that the map $Z\to \mathcal D(\mathbb P_X)$, $z\mapsto [[N^*\hat z]]$ is smooth. 

Given a distribution $\mu\in\mathcal M^{-\infty}(Z)$, define the generalized valuation $\Cr(\mu)\in \mathcal V^{-\infty}(X)$ by $\Cr(\mu)=\pi_*\tau^*\mu$, which is the Radon transform with respect to the Euler characteristic of $\mu$ in the terminology of Alesker \cite{alesker_integral}. Explicitly, for $\psi\in\mathcal V^\infty_c(X)$,
\[ \langle \Cr(\mu),\psi\rangle=\int_{Z}\psi(\hat z) d\mu(z)  \]  
\begin{Lemma}\label{lem:crofton_filtration_level}
	Let $X$ be a manifold, and $\phi=\Cr(\mu)\in  \mathcal V^{-\infty}(X)$ for $\mu\in\mathcal M^{-\infty}(Z)$. Then $\phi\in \mathcal W_{\dim X-\dim W+\dim Z}^{-\infty}(X)$, and $\sigma \phi=(-1)^{\dim X-\dim W+\dim Z}\phi$, where $\sigma$ is the Euler-Verdier involution.
\end{Lemma}

\proof
	We use the Alesker-Poincar\'e duality: for every smooth, compactly supported $\psi\in \mathcal W_{\dim W-\dim Z+1}^{\infty}(X)$, it holds that $\langle \psi, \Cr(\mu)\rangle = \int_Z \psi(\hat z)d\mu(z)=0$. Now Proposition 7.3.2 in \cite{alesker4} implies the first statement.
	For the second statement, simply note that the defining current of $\phi$ is given by $\int_Z [[N^*\hat z]]d\mu(z)$. 
\endproof

\section{Dual Heisenberg manifolds}\label{sec:DH_manifolds}
\subsection{Geometric preliminaries}
Recall that the Heisenberg algebra is $\R^{2n+1}$ whose Lie bracket is given on the standard basis by $[e_i, e_{i+n}]=e_{2n+1}$ for $1\leq i\leq n$, and $[e_i, e_j]=0$ for all $i<j\neq i+n$. The center of the algebra is thus $Z:=\Span(e_{2n+1})$.

The dual Heisenberg algebra $U^{2n+1}$ is $\R^{2n+1}$, equipped with a distinguished hyperplane $H\subset U$, and $\omega\in  \wedge^2 H^*\otimes U/H$ is a non-degenerate (twisted) form. Let us justify this terminology.

\begin{Lemma}
	The space $U^*$ is naturally the Heisenberg Lie algebra.
\end{Lemma}
\proof
Denote $Z=(U/H)^*=H^\perp\subset U^*$. There is a $Z$-valued symplectic form on the quotient $U^*/Z=H^*$: we may consider $\omega$ as the isomorphism $H\simeq  H^*\otimes (U/H)$, or equivalently $ H^*\simeq H\otimes (U/H)^*$, which is a non-degenerate bilinear form $\omega^*\in  (H\otimes H)\otimes (U/H)^*$, which is clearly antisymmetric, that is $\omega^*\in \wedge ^2 (U^*/Z)^* \otimes Z$. The Lie bracket $[x, y]:=\omega^*(x+Z, y+Z)$ then defines on $U^*$ a Heisenberg Lie algebra structure with center $Z$.
\endproof

\begin{Definition}
	A manifold $X$ equipped with a hyperplane distribution $H$ and a smooth field of nowhere-degenerate forms $\omega\in \Gamma^\infty(X, \wedge ^2 H_x^*\otimes (T_xX/H_x))$ will be called a dual Heisenberg (DH) manifold. We call $H$ the horizontal distribution and $\omega$ the DH form.
\end{Definition}
\begin{Remark}
	Ovsienko \cite{ovsienko} defines the closely related notion of a local Heisenberg structure. Similar ideas lie in the foundation of Stein's Heisenberg calculus, see e.g. Ponge \cite{ponge1, ponge2}. 
\end{Remark}

\begin{Example}\label{exm:contact_is_DH}
	\emph{A contact manifold is naturally equipped with a DH structure. Let $(X, H)$ be contact, and let $\alpha$ be a contact form defined locally. Then  $\alpha:T_xX/H_x\xrightarrow{\sim} \R$ is an isomorphism, and we may define the 2-form $\omega(u,v)=\alpha^{-1}(d\alpha(u, v))$ for $u, v\in H_x$. It is independent of the choice of $\alpha$.}
\end{Example}

\begin{Definition}
Let $(X, H_X, \omega_X)$ be a DH manifold. A DH-submanifold is a submanifold $Y^{2k+1}\subset X^{2n+1}$ which is nowhere tangent to $H_X$. It inherits the structure of a DH manifold with hyperplane distribution $H_Y=H_X\cap TY$, and DH sform $\omega_Y=\omega_X|_{H_Y}$, where we naturally identify $TX/H_X=TY/H_Y$. 
\end{Definition}

\begin{Definition}
	Let $U^{2n+1}$ be the dual Heisenberg Lie algebra, with horizontal hyperplane $H\subset U$ and DH form $\omega\in \wedge^2 H^*\otimes U/H$. A Euclidean structure $P$ on $U$ is said to be compatible with the DH structure if one can find an orthonormal basis $X_j, Y_j, Z$ of $U$ such that $X_j, Y_j$ form a symplectic basis of $H$ w.r.t. $Z$: $\omega(X_i, Y_j)=0$, $\omega(X_j, Y_j)=Z+H\in U/H$ for all $1\leq i,j\leq n$ and $i\neq j$. A compatible form $\alpha\in U^*$ is any form such that $\Ker\alpha=H$. $P$ is compatible with $\alpha$ if $\alpha=P(Z, \bullet)$.
\end{Definition}

\begin{Definition}\label{def:compatible_riemannian}
	Let $M^{2n+1}$ be a DH manifold with horizontal distribution $H\subset TM$ and DH forms $\omega\in \wedge ^2 H^*\otimes TM/H$. A Riemannian metric $g$ on $M$ is compatible if $g_x$ is compatible for $(T_xM, H_x, \omega_x)$, for every $x\in M$. 
\end{Definition}

\begin{Remark}
	1. A similar notion of adapted metric was used by Chern and Hamilton in \cite{chern_hamilton} for contact 3-manifolds and in later works by other authors.
	\\
	2. In fact, for the purposes of this work we only need the metric we choose to be compatible with the contact structure at isolated points of interest.
\end{Remark}

\begin{Lemma}
	For any DH manifold $M$, a globally defined compatible Riemannian metric exists. Moreover, if a compatible form $\alpha\in\Omega^1$ is given, one may choose the metric to be compatible with $\alpha$.
\end{Lemma}
\proof
Cover $M$ by contractible open subsets $U_i$, and choose 1-forms $\alpha_i$ on $U_i$ such that $\alpha_i|_{U_i\cap U_j}=\pm \alpha_j|_{U_i\cap U_j}$ - this can be accomplished by fixing an arbitrary Riemannian metric and letting $\alpha_i(x)\in T_x^*M$ have unit norm for all $x\in U_i$. Set $\omega_i=\alpha_i\circ \omega\in\Gamma(U_i, \wedge^2 T^*M)$. Let us also fix a complementing line bundle $L$ to $H$: $L\oplus H=TM$.

For each $U_i$ we may choose a compatible complex structure $J_i$ for $\omega_i$ on $H|_{U_i}$ such that $J_i|_{U_i\cap U_j}=\pm J_j|_{U_i\cap U_j}$. We thus obtain the corresponding Euclidean forms $h_i$ on $H|_{U_i}$ which will satisfy $h_i=h_j$, and so can be patched to give a globally defined positive definite quadratic form $h$ on $H$. Now for every $U_i$, choose $Z_i\in\Gamma(U_i, L)$ satisfying $\alpha_i(Z_i)=1$, and define the Riemannian metric $g_i$ over $U_i$ by $g_i|_H=h$, $g_i(L, H)=0$, $g_i(Z_i, Z_i)=1$. Then clearly $g_i=g_j$ on $U_i\cap U_j $ and thus $(g_i)$ define a Riemmanian metric $g$ on $M$. Over $U_i$ we may locally choose an orthonormal symplectic basis $X^{i}_j, Y^i_j$ for $\omega_i$, and then $\omega(X^i_j, Y^i_j)=\alpha_i^{-1}(1)=Z_i$ for all $j$. This shows $g$ is a compatible metric.

Finally, if $\alpha$ is given, we just take $\alpha_i=\alpha|_{U_i}$ in the construction above.
\endproof

Let $M^{2n+1}$ be a DH manifold with horizontal distribution $H$ and DH form $\omega$. Set $M_H:=\{(x, \xi): x\in M, \xi\perp H_x\}\subset\mathbb P_M(=\mathbb P_+(T^*M))$.

\begin{Definition}\label{def:normally transversal}
	
	A $C^2$ boundaryless submanifold $F\subset M$ is said to be normally transversal (to $H$) if $N^*F$ intersects $M_H$ transversally.
		
	More generally, a $C^2$ submanifold with boundary  $F\subset M$ is normally transversal if $T_xF\not\subset H_x$ for all $x\in\partial F$, and both $N^*(\textrm{int} F)$ and $N^*\partial F$ intersect $M_H$ transversally.

\end{Definition}
Note that if $F$ is normally transversal then so is $\partial F$.

\begin{Lemma}\label{lem:wavefront_normally_transversal}
	If a submanifold with boundary  $F\subset M$ intersects $H$ normally transversally, then $\WF([[N^*F]])\cap N^*M_H=\emptyset$.
\end{Lemma}
\proof
For a boundarlyess submanifold this is immediate: $N^*F$ is a boundaryless submanifold, and $\WF([[N^*F]])=N^*(N^*F)$.

Now assume $F$ has boundary.  Consider the sets $N_{\partial F}^*F=\cup_{x\in\partial F} N^*_xF$ and $N_{\textrm{int}}^*F:=\cup_{x\in  F}\{ \xi\in \mathbb P_+(T_x^*M): \xi\in (T_xF)^\perp\}$. Both are manifolds with the same boundary $B=\{(x,\xi): x\in\partial F, \xi\in (T_xF)^\perp\}$, and $N^*F=N_{\partial F}^*F\cup N_{\textrm{int}}^*F$. Also, $N_{\partial F}^*F\subset N^*\partial F$ is a submanifold with the same boundary $B$.

By assumption, $M_H$ intersects the interior points of $N_{\partial F}^*F$ and $N_{\textrm{int}}^*F$ transversally, and does not intersect $B$. It also holds for $(x,\xi)\notin B$ that $\WF_{x,\xi}([[N^*F]])$ coincides with $ N_{x,\xi}^*(N_{\partial F}^*F)$ if $x\in \partial F$, and with $N_{x,\xi}^*( N_{\textrm{int}}^*F)$ if $x\in\textrm{int}F$. This concludes the proof.
\endproof

Normal transversality is of course generic:
\begin{Lemma}\label{lem:M_H_perturbation}
	Any hypersurface with boundary $F\subset M$ can be perturbed by an arbitrarily $C^\infty$-small  amount to become normally transversal.
\end{Lemma}
\proof
First, we may perturb $F$ so that it is tangent to the contact distribution at isolated points.
Now near interior contact points, we only need to perturb $F$ locally near those points to get normal transversality.
Next, we may perturb $F$ near the boundary to have no contact points of $F$ near the boundary, and isolated contact points of $\partial F$.
We then perturb $F$ near those contact points s.t. $\partial F$ has normal transversality at all its contact points.
\endproof

Normal transversality is universal to embeddings, as the following lemma shows.
\begin{Lemma}\label{lem:normal_universality}
	Let $(Y, H^Y)\subset (X, H^X)$ be DH manifolds, and $F\subset Y$ a normally transversal submanifold with boundary. Then $F\subset X$ is also normally transversal.
\end{Lemma}
\proof
We will distinguish the conormal bundles of $F$ in the different ambient manifolds by writing e.g. $N^*_YF\subset \mathbb P_Y$. Recall that $H^Y_x=T_xX\cap H^X_x$. Thus for $x\in\partial F$, $T_xF\not\subset H^X_x$, since $T_xF\subset H^Y_x$. 

Next take $(x,\xi)\in N_X^*F$. For $x\in \textrm{int} F$ we should check $T_{x,\xi}N^*_XF\cap T_{x,\xi}X_{H^X}=\{0\}$. 
Consider the natural surjective map $\beta:\mathbb P_X|_Y\setminus N^*X\to \mathbb P_Y$. It holds that $\beta(N^*_XF)=N^*_YF$.
Denote $X_{H^X}|_Y:=X_{H^X}\cap  \pi_X^{-1}(Y)$ where $\pi_X:\mathbb P_X\to X$. Since $X\subset Y$ is a DH submanifold, it holds that $\beta:X_{H^X}|_Y\to Y_{H^Y}$ is a diffeomorphism. If $v\in T_{x,\xi}N^*_XF\cap T_{x,\xi}X_{H^X}|_Y$ then $d\beta(v)\in T_{x,\xi}Y_{H^Y}\cap T_{x,\xi}N^*_YF$ is non-zero, a contradiction to the assumption of normal transversality of $F\subset M$.
The case of $x\in\partial F$ is virtually identical.

\endproof

\subsection{Constructing the valuations}
Consider a point $q=(p,\xi)\in M_H$. The vertical subspace is $\Ker (d\pi)=\xi^*\otimes T^*_pM/\xi=(H_p^\perp)^*\otimes H_p^*$. Recall that $\omega_p\in\wedge^2 H_p^* \otimes T_pM/H_p$ is non-degenerate, so that $H_p^*\simeq H_p\otimes H_p^\perp$. Thus $\Ker(d\pi)\simeq H_p$.
Note also that $d\pi: T_qM_H\to T_pM$ is an isomorphism, so that we have a natural decomposition

\begin{equation}\label{eqn:contact_decomposition}T_q \mathbb P_M\simeq T_pM\oplus H_p\end{equation}
Denoting the contact hyperplane of $\mathbb P_M$ by $\widehat H_q$, we get 
\begin{equation}\label{eqn:contact_hyperplane}
\widehat H_q\simeq H_p\oplus H_p.
\end{equation}

We will define certain  $\eta_{k}\in\Gamma^{-\infty}_{M_H}(\mathbb P_M, \Omega^{2n+1}\otimes \pi^*o_M)$, $1\leq k\leq 2n$ (twisted by the pull-back of the orientation bundle of $M$) supported on $M_H$.

Take $\tilde \psi\in \Omega_c^{2n+1}(\mathbb P_M)$ a Legendrian form. We have \[\tilde \psi|_q\in \wedge^{2n+1}T_q^*(\mathbb P_M)\simeq\bigoplus_{k=0}^{2n} \wedge ^{2n+1-k}T^*_pM\otimes \wedge^{k} H_p^*\]
Moreover, since $\tilde \psi$ is Legendrian, its first factor in the $k$-th component belongs to the kernel of the restriction map $\wedge^{2n+1-k}T_p^*M\to\wedge^{2n+1-k}H_p^*$, which is naturally isomorphic to $H_p^\perp\otimes\wedge^{2n-k} H_p^* $. Thus the $k$-th component of $\tilde\psi$, denoted  $\tilde\psi^{2n+1-k,k}$, lies in 
$ H_p^\perp \otimes\wedge^{2n-k}H_p^*\otimes \wedge^k H_p^*$. 

Wedging the last two factors, we then get an element $(\tilde\psi)_k\in \wedge^{2n}H_p^*\otimes( T_pM/H_p)^* \simeq \wedge ^{2n+1}T^*_pM\simeq \wedge^{2n+1}T^*_q M_H$. Now set
\[ \langle \tilde\eta_{k},\tilde\psi \otimes \epsilon_M \rangle :=\int_{M_H} (\tilde\psi)_k \otimes \epsilon_M \]
where $\epsilon_M$ is a section of $\pi^* o_M$, so that $\tilde{\psi}_k \otimes \epsilon_M\in \Dens(T_q M_H)$. Thus $\tilde \eta_k$ is a linear functional on compactly supported, Legendrian, $\pi^*o_M$-twisted $(2n+1)$-forms, and it is supported on $M_H$. For $k=0$, $\tilde{\eta}_0=[[M_H]]$ is in fact a well-defined (twisted) generalized $2n$-form.

Finally for $1\leq k\leq 2n$, define $\eta_k\in \Omega_{-\infty}^{2n+1}(\mathbb P_M)\otimes \pi^*o_M$ by setting for $\psi\in \Omega_c^{2n}(\mathbb P_M)$ \[\langle \eta_k, \psi\otimes\epsilon_M\rangle:= \langle \tilde \eta_k, D\psi\otimes \epsilon_M\rangle\] where $D$ denotes the Rumin differential. It is clear that $\eta_k$ is a (twisted) Legendrian cycle, and $\pi_*\eta_k=0$.

\begin{Definition}
	Define $\phi_k\in\mathcal V^{-\infty}(M)$ for $0\leq k\leq 2n$ as follows: for $1\leq k\leq 2n$ by $\phi_k=[\frac12 \eta_k, 0]$; for $k=0$, $\phi_0$ is represented by $(\frac12\tilde \eta_0, 0)$.
\end{Definition}
\begin{Remark}
	It is in fact possible to define natural curvature measures globalizing to $\phi_k$: extend $\tilde \eta_k$ arbitrarily as a functional on all forms, that is $\tilde \eta_k\in\Omega_{-\infty}^{2n}(\mathbb P_M)$, and define the generalized curvature measure $\Phi_k$ represented by the pair of forms $(\frac12 \tilde \eta_k, 0)$. The resulting curvature measure is independent of the extension, as $\tilde \eta_k$ is only applied to subsets of conormal cycles.
\end{Remark}
Note that $\phi_0(X)=\int_{N^*X}\frac{1}{2}\tilde \eta_0=\frac12 \int_{N^*X}[[M_H]]$ is one-half the intersection index of $M_H$ and $N^*X$, both oriented locally by a fixed local orientation on $M$. 
In the following proofs, we often assume for simplicity $M$ is oriented. They are easily adjusted for the general case.
\begin{Lemma}
	$\phi_0$ is the Euler characteristic.
\end{Lemma}
\proof
Note that $[[M_H]]$ is a closed current, and denoting $\pi:\mathbb P_M\to M$ clearly $\pi_*[[M_H]]=2$. It follows that $\phi_0=\chi$ by an obvious extension  to $\mathcal V^{-\infty}$ of \cite[Corollary 1.5 ]{bernig_brocker}.
\endproof

Recall that $a:\mathbb P_M\to\mathbb P_M$ is the fiberwise antipodal map, and $\sigma:\mathcal V^{-\infty}(M)\to \mathcal V^{-\infty}(M)$ denotes the Euler-Verdier involution.
\begin{Proposition}\label{prop:euler_verdier_one}
	For all $0\leq k\leq 2n$ it holds that $\sigma \phi_k=\phi_k$.
\end{Proposition}
\proof
 For $\psi\in\Omega^{2n+1}(\mathbb P_M)$ it holds that $\int_{M_H}\psi=\int_{M_H}a^*\psi$, since $a:M_H\to M_H$ is clearly orientation preserving. Note also that $a:\mathbb P_M\to \mathbb P_M$ is orientation reversing as the fibers of $\mathbb P_M\to M$ are even-dimensional spheres. Hence for a test form $\psi\in\Omega_c^{2n}(\mathbb P_M)$,
\[\langle a^*\eta_k, \psi\rangle=- \langle\eta_k, a^*\psi\rangle=-\int_{M_H}[Da^*\psi]_k=-\int_{M_H}a^*[D\psi]_k=-\int_{M_H}[D\psi]_k=-\langle \eta_k, \psi\rangle\]
That is, $a^*\eta_k=-\eta_k$ for all $k$, and hence 
$\sigma \phi_k=(-1)^{2n+1}(-\phi_k)=\phi_k$ as asserted.
\endproof

\begin{Proposition}\label{prop:phi_k_filtration_level}
	It holds that $\phi_{k}\in \mathcal W_{k}^{-\infty}(M)$.
\end{Proposition}
\proof
Take a closed Legendrian form $\psi\in \Omega^{2n+1}(\mathbb P_M)$ defining a smooth valuation $\Psi\in \mathcal W_{2n+2-k}^\infty(M)$. We ought to show that $\phi_{k}\cdot \Psi=0$, equivalently $\langle \tilde \eta_{k}, \psi\rangle=0$.
But $\Psi\in \mathcal W_{2n+2-k}^\infty(M)$ implies that $\psi$ has a horizontal degree at least $2n+2-k$, and the claim follows from the definition of $\tilde\eta_k$.
\endproof

\begin{Lemma}\label{lem:small_wavefront}
	The wavefront set of $\phi_k$ for $1\leq k\leq 2n$ is $N^*M_H$.
\end{Lemma}
\proof
This is immediate from definition: the wavefront set of $[[M_H]]$ is $N^*M_H$, and restriction to $M_H$ is the only source of singularities of $\eta_k$.
\endproof
It follows by Lemma \ref{lem:wavefront_normally_transversal} that we may evaluate $\phi_k$ on any normally transversal submanifold with boundary.

Let $F\subset M$ be a smooth hypersurface in a DH-manifold $M^{2n+1}$, assume that $F$ is normally transversal to the horizontal distribution, and tangent to it at $p\in M$.

Working in a small open ball $U$ near $p$ with no other contact points, let us fix a 1-form $\alpha\in\Omega^1(U)$ defining the horizontal distribution, which also trivializes the DH form: $\omega\in \Gamma(U,\wedge^2 H^*)$. Let $g$ be a Riemannian metric on $U$ which is compatible with $\alpha$ and $\omega$ on $T_pM$ (but not necessarily elsewhere), with Levi-Civita connection $\nabla$. Let $R:U\to SM$ be the vector field given by $R(x)\perp H_x$,  $g(R, R)=1$, $\alpha(R)>0$. Let $(X_j)_{j=0}^{2n}$ be an orthonormal frame in $U$ with $X_0|_F=\nu:F\cap U\to SM$ the unit normal oriented by $\alpha(\nu)>0$, and $(X_j(p))_{j=1}^{2n}$ a symplectic basis of $H_p$. Let $\theta_j\in\Omega^1(U)$ be the dual coframe to $X_j$. Let $S=(s_{ij})_{i,j=1}^{2n}$, $s_{ij}=\theta_i(\nabla_{X_j}\nu)$ be the second fundamental form given by $\nabla_{X_j}{\nu} = \sum_{i=0}^{2n}s_{ij} X_i$. 
\begin{Definition}\label{def:second_fund_form}
	Let $h_{ij}$, $1\leq i\leq 2n$, $0\leq j\leq 2n$ be given by $h_{ij}= \theta_i(\nabla_{X_j}R)$.    
	The matrix $h=(h_{ij})_{i,j=1}^{2n}$ is the second fundamental form of the contact structure.
\end{Definition}
\begin{Remark}
This is a slightly different definition than the one in \cite{reinhart}, where the second fundamental form is symmetrized.
\end{Remark}
Define the matrix $A_p=(h_{ij}-s_{ij})_{i,j=1}^{2n}$.

\begin{Proposition}\label{prop:DH_general_formula}
	
	Let $F$ be a normally transversal closed hypersurface. It then holds that
	
	\begin{equation}\label{eq:DH_general_formula}\phi_{k}(F)= {2n \choose k}\sum_{p: T_pF=H_p} |\det A_p|^{-1}D(A_p[2n-k], J[k]) 	\end{equation}
	where $D$ is the mixed discriminant.
	
\end{Proposition}
\begin{Remark}
	We thus see that geometrically, $\phi_k$ is reminiscent of the $k$-th elementary symmetric polynomial in the principal radii of a hypersurface in a Riemannian manifold. In particular, $\phi_{2n}^{-1}$ plays the role of the absolute value of the gaussian curvature. In fact, for $M=U$ - the dual Heisenberg algebra, $\phi_{2n}(F)$ is precisely the inverse absolute value of the gaussian curvature, summed over all contact points.
	\end{Remark}

\proof
Write $\omega_i$ for the Ehresmann connection on $SU$, namely $\omega_i=\pi_V^*\theta_i$, where $\pi_V:T_{p,\xi}SU\to T_p M$ is the projection to the vertical tangent space.
It then holds that $R^*\omega_i=\sum_{j=0}^{2n} h_{ij}\theta_j$, and  $\nu^*\omega_i=\sum_{j=1}^{2n}s_{ij}\theta_j$. Recall that $(s_{ij})$ is symmetric. For notational simplicity, we write $\theta_j$ also for $\pi^* \theta_j\in \Omega^1(SU)$.

First we describe the normal cycle $NF=\{(x, \pm \nu(x)): x\in F\}\subset SM$ near $(p,\xi)=X_0(p)\in M_H$, explicitly as a generalized form $\omega_{NF}:=[[NF]]\in \Omega_{-\infty}(SM)$. 
Let $W\subset \pi^{-1}U$ be a small neighborhood of $(p, \xi)$.
Write $\sigma_{SM}:=\wedge_{i=0}^{2n}\theta_i \wedge\wedge_{i=1}^{2n} \omega_i$, and $\theta_F=\wedge_{i=1}^{2n}\theta_i\in \Omega^{2n}(F\cap U)$.

 \textit{Claim}. $ \omega_{NF}^W:=\omega_{NF}|_{W}\in \Omega_{-\infty}(W)$ is given by
\[ \omega_{NF}= \theta_0\wedge \wedge_{i=1}^{2n}(\omega_i-\sum_{j=1}^{2n}s_{ij}\theta_j)\delta_{NF}^W \]
where

\[\langle \delta^W_{NF}, \mu(x,\zeta)\sigma_{SM}\rangle=\int_{F\cap \pi(W)} \mu(x, \nu(x)) \theta_F \]

\textit{Proof of Claim.}
Take $\psi\in \Omega^{2n}_c(W)$. We should check that 

\[ \int_{NF}\psi= \int_F \nu^*\left(\frac{\psi\wedge \theta_0\wedge  \wedge_{i=1}^{2n}(\omega_i-\sum_{j=1}^{2n}s_{ij}\theta_j)}{\sigma_{SM}}\right) \theta_F \]
which reduces to the pointwise verification
\[ \frac{\nu^*\psi}{\theta_F}=\nu^*\left(\frac{\theta_0\wedge \psi\wedge  \wedge_{i=1}^{2n}(\omega_i-\sum_{j=1}^{2n}s_{ij}\theta_j)}{\sigma_{SM}}\right)  \]
We will check this equality for a basis of the $2n$ forms, which is given by 2n-wedges of $(\theta_i)_{i=0}^{2n}$, $(\omega_j)_{j=1}^{2n}$. 
If $\psi$ contains a $\theta_0$ factor, clearly both sides vanish.

Assume $\psi = \theta_1\wedge\dots\wedge \theta_k\wedge \omega_{k+1}\wedge\dots\wedge\omega_{2n}$.
The left hand side is easily seen to be equal $\det(s_{ij})_{i,j=k+1}^{2n}$. The right hand side is 
\[\nu^*\left(\frac{(-1)^{2n-k}\det(s_{ij})_{i,j=k+1}^{2n}\theta_0\wedge\dots\wedge \theta_k\wedge \omega_{k+1}\wedge\dots\wedge\omega_{2n}\wedge \omega_1\wedge\dots\wedge \omega_k
	\wedge \theta_{k+1}\wedge\dots\wedge\theta_{2n}}{\sigma_{SM}}\right) \]
and after reordering we get $(-1)^{2n-k+k(2n-k)}\det(s_{ij})_{i,j=k+1}^{2n}=\det(s_{ij})_{i,j=k+1}^{2n}$, concluding the proof of the claim.\qed
\\\\
We wish to apply $\tilde \eta_k$ to $ \omega_{NF}^W$. Since $NF$ intersects $M_H$ transversally at the isolated point $p$, we need only look at the value of $ \omega_{NF}^W$ at $(p,\pm X_0(p))$.

Recall the decomposition $T_{p,\xi}\mathbb P_M=T_{p,\xi}M_H\oplus T_\xi(\mathbb P_+(T^*_pM))$.
We denote by $W_{p,\xi}=T_{p,\xi}M_H=T_pM$ the contact-horizontal subspace, and by $V_{p,\xi}=T_\xi(\mathbb P_+(T^*_pM))=(H_p\otimes H_p^\perp)^*$ the vertical subspace.
Using the Riemannian structure, we get a decomposition $T_{p,\xi}(SM)=W_{p,\xi}\oplus V_{p,\xi}$.
The DH form $\omega\in \wedge^2 H_p^*\otimes T_pM/H_p$ can be recast as the isomorphism $\Omega:H_p\otimes H_p^\perp\to H_p^*$, that is $\Omega:V_{p,\xi}^*\to H_p^*$.

We first need to identify the $(2n+1-k, k)$ component (with respect to decomposition \eqref{eqn:contact_decomposition}) of $ \omega_{NF}^W$ over $(p, \xi)\in M_H$. 
The forms $\theta_i$ vanish on the vertical subspace.
The contact-horizontal component of $\omega_i$ is given by $R^*\omega_i=\sum_{j=0}^{2n} h_{ij}\theta_j$. We will write the resulting decomposition as $\omega_i=R^*\omega_i\oplus\omega_i^V$
where $R^*\omega_i$ is the contact-horizontal component, while $\omega_i^V$, the restriction of $\omega_i$ to  $V_{p,\xi}$, is the vertical component. 

Define $(\gamma_{ij})_{i,j=1}^{2n}$ by $\Omega(\omega_i^V)=\sum_{j=1}^{2n}\gamma_{ij}\theta_j\in H_p^*$. Note that under the Euclidean identification, $V_{p,\xi}=H_p$ and $\omega_i^V=\theta_i$ for $i=1,\dots,2n$.  Compute $\gamma_{ij}=\Omega(\theta_i)(X_j)=\omega(X_i, X_j)$. By assumption, $X_j$ is a symplectic basis of $H_p$, so that $\gamma_{ij}$ is the standard $2n\times 2n$ matrix $J$ representing $\sqrt{-1}$.

Write $a_{ij}=h_{ij}-s_{ij}$, $1\leq i,j\leq 2n$. Then 

\[  \omega_{NF}^W= \delta_{NF}^W\cdot\theta_0\wedge \wedge_{i=1}^{2n} (\omega_i^V+\sum_{j=1}^{2n}a_{ij}\theta_j)\] so that
\[ (\omega_{NF}^W)^{2n+1-k,k}=\delta_{NF}^W\cdot\theta_0\wedge  \sum_{|I|=k}\epsilon_I \wedge_{i\notin I} \left(\sum_{j=1}^{2n}a_{ij}\theta_j\right)\wedge\wedge_{i\in I}\omega_i^V \]
where $\epsilon_I=(-1)^{i_2-i_1+\dots + i_{2j}-i_{2j-1}-j}$ is the sign of the permutation $\sigma_I={1\dots 2n \choose I^c, I}$, where $I=\{i_1<\dots <i_k\}$, $I^c=\{1,\dots,2n\}\setminus I$ is ordered increasingly, and $j=\lfloor \frac k 2 \rfloor$.

Applying $\Omega$ to the last $k$ factors and subsequently wedging all the factors to get a top form on $T_pM$, we get

\begin{align*}\langle\tilde \eta_k,  \omega_{NF}^W\rangle&=\int_{M_H}\theta_0\wedge  \sum_{|I|=k}\epsilon_I \wedge_{i\notin I} \left(\sum_{j=1}^{2n}a_{ij}\theta_j\right)\wedge\wedge_{i\in I}\left(\sum_{j=1}^{2n}\gamma_{ij}\theta_j\right)\cdot \delta_{NF} ^W \\&
= \sum_{|I|=k}\sum_{\sigma\in S_{2n}}\sign\sigma \prod_{i\notin I}a_{i\sigma(i)}\prod_{i\in I}\gamma_{i\sigma(i)}\int_{M_H} \theta_0\wedge\dots\wedge \theta_{2n}\cdot \delta_{NF}^W   \end{align*}

Recall $A_p=(a_{ij})=h-S$ where $S=(s_{ij})$, $h=(h_{ij})$, $1\leq i,j\leq 2n$, and let $B^\alpha$ be the $\alpha$-row of a matrix $B$. For an ordered subset $I\subset \{1,\dots,2n\}$, let $(A_p^{I^c}, J^I)$  denote the matrix with the corresponding columns. Then
\[\langle\tilde \eta_k,  \omega_{NF}^W\rangle= \sum_{|I|=k} \det (A_p^{I^c}, J^{I})\int_{M_H} \theta_0\wedge\dots\wedge \theta_{2n}\cdot \delta_{NF}^W  \]

Define \[B(\alpha)=\left\{\begin{array}{cc} A_p,& \alpha\leq 2n-k\\ J,& \alpha>2n-k\end{array} \right.\] The mixed discriminant is given by

\[ D(A_p[2n-k], J[k])=\frac{1}{(2n)!}\sum_{\tau\in S_{2n}}\det( B(\tau(i))_{i,j})=\frac{1}{(2n!)}\sum_{|I|=k} k!(2n-k)!\det (A_p^{I^c}, J^{I})  \] 
that is
\[ \langle\tilde\eta_k, \omega_{NF}^W\rangle={2n\choose k} D(A_p[2n-k], J[k]) \int_{M_H} \theta_0\wedge\dots\wedge \theta_{2n}\cdot \delta_{NF}^W \]

Recall that we should fix a section of the orientation bundle of $M$ over $SM$ to get numerical values for the integral. 
Let us verify that for a choice of $\epsilon_\theta\in \pi^*o_M$ corresponding to $\theta_0\wedge\dots\wedge \theta_{2n}$ we get the identity 
\[ \int_{M_H} \theta_0\wedge\dots\wedge \theta_{2n}\cdot \delta_{NF}^W\otimes \epsilon_\theta=|\det A_p|^{-1}. \]
One can compute it directly, but we can do something simpler: Observe that $\langle \tilde \eta_0, \omega_{NF}^W\rangle$ is just the intersection index $I_{p,\xi}$ of $NF$ and $M_H$ at $(p,\xi)$. Note that the order of intersection is not important as $\dim NF$ is even, while the orientations of $M_H$ and $NF$ are determined by $\epsilon_\theta$.
\\
By what we have proved, we see that \[I_{p,\xi}(NF,M_H)=\langle\tilde\eta_0,  \omega_{NF}^W\rangle =\det A_p \int_{M_H}\theta_0\wedge\dots\wedge\theta_{2n}\delta_{NF}^W.\]
It thus remains to verify that $I_{p,\xi}(NF, M_H)=\sign \det A_p$.

The positive orientation on $T_{p, \xi}SM$ is given by the dual basis $\theta_0,\dots, \theta_{2n}, \omega_1, \dots,\omega_{2n}$. To see that, consider a homotopy of the Riemannian metric between our metric and a flat one, and some corresponding homotopy of the orthonormal frame $X_j$. The dual basis above remains a basis throughout the homotopy, and clearly defines the positive orientation in the flat case.  

Considering $\nu:F\to SM$, $R:F\to SM$ as maps, we get a positive basis of $T_{p,\xi}M_H$ given by $D_pR(X_0),\dots, D_pR(X_{2n})$, and a positive basis of $T_{p,\xi}NF$ given by $D_p\nu(X_1),\dots, D_p\nu(X_{2n})$.

Form the $(4n+1)\times (4n+1)$ matrix 

\[B=\left(\begin{array}{cc} \theta_i(D_pR(X_j))_{i=0\dots 2n}^{j=0\dots 2n} & \theta_i(D_p\nu(X_j))_{i=0\dots 2n}^{j=1\dots 2n}\\ \omega_i(D_pR(X_j))_{i=1\dots 2n}^{j=0\dots 2n} &  \omega_i(D_p\nu(X_j))_{i=1\dots 2n}^{j=1\dots 2n} \end{array}\right) =  \left(\begin{array}{cc} \delta_{ij} & \delta_{ij}\\ h_{ij} &  s_{ij}\end{array}\right)  \]
By definition, $I_{p,\xi}(NF,M_H)=\sign\det B$.
Now for $1\leq i\leq 2n$ substract column $2n+1+i$ from column $i+1$. Finally, interchange each culumn of index $2\dots 2n+1$ with the respective column from the last $2n$ columns, resulting in $4n^2$ swaps. The resulting matrix is block-triangular, and has determinant $\det B=\det (h_{ij}-s_{ij})_{i=1\dots 2n}^{j=1\dots 2n}=\det A_p$, verifying our assertion.

Accounting also for the point $\xi=H^\perp$ with the opposite orientation, we conclude that

\[\langle \tilde \eta_k, \omega_{NF}\rangle =2 {2n \choose k}\sum_{p: T_pF=H_p} |\det A_p|^{-1}D(A_p[2n-k], J[k]) \]
and the statement follows.

\endproof
\begin{Definition}\label{def:local_areas} For a normally transversal contact point $p$ of $F$, the local contact areas are \[\phi_k(F, p):= {2n \choose k} |\det A_p|^{-1}D(A_p[2n-k], J[k]).\]
\end{Definition}

\begin{Corollary}\label{cor:DH_general_formula}
	When $k=2n$, equation \eqref{eq:DH_general_formula} remains valid also for a normally transversal hypersurface with boundary $F$.
\end{Corollary}
\proof
By Lemma \ref{lem:wavefront_normally_transversal}, $\phi_{2n}(F)$ is well-defined. Using the notation of the proof of Lemma \ref{lem:wavefront_normally_transversal}, we may write 
\[\phi_{2n}(F)=\frac12 \langle \tilde \eta_{2n}, [[N^*_{\textrm {int}}F]] \rangle+\frac12 \langle\tilde \eta_{2n}, [[N^*_{\partial F}F]]\rangle.    \]
 The first summand is computed in Proposition \ref{prop:DH_general_formula}. The second summand trivially vanishes.
\endproof

\begin{Corollary}\label{cor:positive_area}
For a normally transversal hypersurface with boundary $F$, $\phi_{2n}(F)\geq 0$,  with equality if and only if it is nowhere tangent to $H$.
\end{Corollary}
\proof
Immediate from Corollary \ref{cor:DH_general_formula} and eq. \eqref{eq:DH_general_formula} .
\endproof
 It will be convenient to extend $\phi_{2n}$ to general hypersurfaces by setting $\phi_{2n}(F)=\infty$ when $F$ is not normally transversal.
 
Next we establish the universality of $\phi_{k}$ with respect to embeddings.

Let $(X^{2n+1}, H^X, \omega)$ be a DH-manifold, and $Y^{2m+1}\subset X^{2n+1}$ a DH-submanifold. Let $i:Y\to X$ denote the embedding. We will write $\eta_{2k}^X$, $\phi_{2k}^X$ for the generalized forms and valuations defined over $X$, and similiarly for $Y$. 

\begin{Theorem}\label{thm:weyl}
	It holds that $i^*\phi_{k}^X=\phi_{k}^Y$.
\end{Theorem}
\proof
The following constructions appeared in \cite{alesker_integral}, and we refer therein for more details. Consider the natural submersion $\beta:\mathbb P_X|_Y\setminus N^*Y\twoheadrightarrow \mathbb P_Y$ and the inclusion $\alpha:\mathbb P_X|_Y\hookrightarrow \mathbb P_X$. Let $\tilde \pi:W\to\mathbb P_X|_Y$ be the oriented blow-up of $\mathbb P_X|_Y$ along the conormal bundle $N^*Y\subset\mathbb P_X$, and $\tilde \alpha: W\to \mathbb P_X$ the corresponding lift. Let $\tilde{\beta}:W\to\mathbb P_Y$ be induced by the restriction map $X\times _Y T^*Y\to T^*X$. For a valuation $\Psi\in \mathcal V^\infty(Y)$ defined by the closed Legendrian form $\psi\in \Omega^{2m+1}(\mathbb P_Y)$, $i^*\Psi$ is given by the current $\tilde{\alpha}_*\tilde \beta^*\psi$.  
Now $\langle \phi_{k}^Y, \Psi\rangle= \langle \tilde{\eta}_{k}^Y, \psi\rangle$ and similarly for $X$. We should thus verify that 
\[ \langle \tilde{\eta}_{k}^Y, \psi\rangle=\langle \tilde \eta_k^X,\tilde{\alpha}_*\tilde \beta^*\psi\rangle  \]
Note that $\tilde \eta_{k}^X$ is supported on $X_{H^X}\subset \mathbb P_X$, which by assumption is disjoint from $N^*Y$. Hence the right hand side can be replaced by $\langle \tilde \eta_k^X,\alpha_* \beta^*\psi\rangle$.

Next take $p\in Y$, $\xi=(H^Y_p)^\perp\in\mathbb P_Y$, and $(p,\xi')\in\mathbb P_X|_Y$ s.t. $\beta(p,\xi')=(p,\xi)$. We may decompose $T_{p,\xi'}\mathbb P_X|_Y=T_pY\oplus (H_p^X)^*\otimes (H_p^X)^{\perp*}$. Consider $d\beta: T_{p,\xi'}\mathbb P_X|_Y\to T_{p,\xi}\mathbb P_Y $, so that
\[ d\beta^*: T^*_pY\oplus H_p^Y\otimes (T_pY/H_p^Y)^*\to T^*_pY\oplus  H_p^X\otimes (T_pX/H_p^X)^*  \]
acts as the identity on the first summand. On the second summand, recalling we have the canonical identificaion $T_pY/H_p^Y\simeq T_pX/H_p^X$ induced by the inclusion $T_pY\subset T_pX$, $d\beta$ simply acts as the inclusion $H_p^Y\hookrightarrow H_p^X$. 
Let us denote by $\Omega_k^Y: \wedge^k H^Y_p\otimes (H^Y_p)^\perp\to \wedge ^k (H^Y_p)^*$ the isomorphism induced by the DH form $\omega$, and similarly for $X$. It follows that the following diagram commutes:

\begin{displaymath}
\xymatrix{
	H_p^Y\otimes (H_p^Y)^\perp \ar[r]^--{d\beta^*}\ar[d]^{\Omega_Y} &  H_p^X\otimes (H_p^X)^\perp \ar[d]^{\Omega_X}\\
	(H_p^Y)^*& ( H_p^X)^*\ar[l]^{\textrm{res}} }
\end{displaymath}
We are left with verifying the identity

\[ \int_{Y_H}\psi_0\otimes \epsilon_Y= \int_{X_H}\alpha_*\beta^*\psi_0\otimes \epsilon_X  \] for an arbitrary Legendrian form $\psi_0\in \Omega^{2m+1}(\mathbb P_Y)$.
But this is now equivalent to the statement $i^*\phi_0^X=\phi_0^Y$, which holds as both sides are just the Euler characteristic.

\endproof

Weyl's principle for DH manifolds, which we just established, is readily applicable in conjunction with the following technical lemma.

\begin{Lemma}\label{lem:intermediate_DH}
	Consider a compact submanifold with boundary $F^{2k}\subset M^{2n+1}$ of a DH manifold $(M, H, \omega)$. 
	\begin{enumerate}
		\item 	If $2k\leq n$, $F$ lies in fact inside a DH submanifold $N^{2k+1}\subset M$. 
		\item For arbitrary $k<n$, one may find a pair of DH-manifolds $N^{2k+1}\subset X^{4n+1}$ such that we get a commuting diagram of DH manifolds
	
		\begin{displaymath}
		\xymatrix{& N^{2k+1}\ar@{^{(}->}[rd] & \\  F^{2k} \ar@{^{(}->}[ru] \ar@{^{(}->}[r]& M^{2n+1}\ar@{^{(}->}[r] & X^{4n+1} &
			 }
		\end{displaymath}
		where all inclusions are DH-embeddings. 
		
	\end{enumerate}

\end{Lemma}
\proof
Assume first $2k\leq n$. Choose a Riemannian metric $g$ on $M$, and let $L:F\to \mathbb P(TM)|_F$ be given by $L(x)=H_x^\perp$, the orthogonal complement with respect to $g$. Note that $L(F)$ does not intersect $\mathbb P(H)|_F\subset\mathbb P(TM)|_F$. By the transversality theorem, and since \[\dim L(F)+\dim \mathbb P(TF)<\dim \mathbb P(TM)|_F\iff 2k+4k-1<2n+2k\iff 2k\leq n,\] we may perturb $g$ so that $L(F)$ avoids $\mathbb P(TF)\subset \mathbb P(TM)$. We then take $N$ to be the image under the exponential map of a small neighborhood of the zero section in $TF\oplus L(F)\subset TM|_F$. It is clearly a DH-submanifold containing $F$.

Now in the general case, consider $\widehat M_H\subset\mathbb P(T^*M)$, which is the quotient of $M_H\subset \mathbb P_M$ under the two-covering map $\mathbb P_M\to \mathbb P(T^*M)$. Define a DH structure on a neighborhood of $\widehat M_H$ as follows: the horizontal structure $\widehat H$ will be the canonical contact structure of $\mathbb P(T^*M)$. For $(x,\xi)\in\widehat M_H$, by eqns. \eqref{eqn:contact_decomposition} and \eqref{eqn:contact_hyperplane}, $\widehat H_{x,\xi}$ is canonically identified with $H_x\oplus H_x\subset T_xM\oplus H_x$. Noting that $T_{x,\xi}\mathbb P(T^*M)/\widehat H_{x,\xi}\simeq T_xM/H_x$, define $\omega_X\in \wedge^2\widehat H_{x,\xi}^*\otimes T_{x,\xi}\mathbb P(T^*M)/\widehat H_{x,\xi}$ by \[\omega_X((u_1,v_1), (u_2, v_2))=\omega(u_1, v_1)+\omega(u_2, v_2)\quad \forall (u_1, v_1), (u_2, v_2)\in \widehat H_{x,\xi}.\] Extend $\omega_X$ arbitrarily to a global section of $\wedge ^2\widehat H^*\otimes T\mathbb P(T^*M)/\widehat H$. Take $X$ to be a small neighborhood of $\widehat M_H$ in which $\omega_X$ is non-degenerate. Then $(X, \widehat H, \omega_X)$ is a DH manifold, and the DH-submanifold $\widehat M_H\subset X$ is clearly isomorphic to $M$. It remains to note that $\dim X=2\cdot2n+1$ and $2k< 2n$, and so we may choose the desired $N^{2k+1}\subset X$ by the first case.
\endproof

\begin{Corollary}\label{cor:general_positivity}
	Let $F^{2k}\subset M$ be a submanifold with boundary, normally transversal to the horizontal distribution $H$. Then $\phi_{2k}(F)\geq0$, with equality if and only if it is nowhere tangent to $H$.
\end{Corollary}
\proof
This follows at once from Lemmas \ref{lem:intermediate_DH} and \ref{lem:normal_universality}, Corollary \ref{cor:positive_area} and Theorem \ref{thm:weyl}. 
\endproof

Finally, we consider the relations between the valuations we constructed.

\begin{Proposition}\label{prop:valuations_on_DH_manifold}
	(i) On any DH manifold,  $\phi_{2b}\in\mathcal W^{-\infty}_{2b}\setminus \mathcal W^{-\infty}_{2b+1}$ as $0\leq b\leq n$. In particular, $\phi_{2b}$, $0\leq b\leq n$ are linearly independent valuations.
	\\(ii) On a generic DH manifold (in a sense to be made precise in the proof), $(\phi_b)_{b=0}^{2n}$ are all linearly independent.
\end{Proposition}
\proof

By Proposition \ref{prop:phi_k_filtration_level}, $\phi_k\in\mathcal W_k^{-\infty}(M)$. For a normally transversal hypersurface $F$ write \begin{align*}\phi_k(F,p)&=|\det A_p|^{-1}D(A_p[2n-k], J[k])\\
&=(-1)^k|\det A_p|^{-1} \sum_{j=0}^{2n-k}{2n-k \choose j} D(S_p[2n-k-j], h_p[j], J[k])\end{align*} 

For even $k=2b$, the term of highest order in $S_p$ (corresponding to $j=0$) is non-zero for generic $S_p$. In particular, $\phi_{2n}(F)\neq 0$ for generic hypersurfaces $F$. 
Now choose any DH submanifold $N\subset M$ of dimension $2b+1$. Using Theorem \ref{thm:weyl}, Lemma \ref{lem:normal_universality} and the last observation, we may find a $2b$-dimensional submanifold $F^{2b}\subset N$ s.t. $\phi^M_{2b}(F)=\phi^N_{2b}(F)\neq 0$. Thus $\phi_{2b}\notin \mathcal W_{2b+1}^{-\infty}(M)$, proving (i).

For (ii), it now suffices to check that $\phi_{2b-1}, \phi_{2b}$ are linearly independent in $\mathcal W^{-\infty}_{2b}(M)/\mathcal W^{-\infty}_{2b+1}(M)$, for a generic $M$.  In turn this is implied by the following statement: as a function on normally transversal closed submanifolds of dimension $2b$ lying inside a fixed DH submanifold $N^{2b+1}\subset M$ , $\phi_{2b-1}, \phi_{2b}$ are linearly independent. By Theorem \ref{thm:weyl}, we may assume $b=n$. Now note that for $\phi_{2n-1}$, the summand of degree $1$ in $S_p$ vanishes, as $S_p$ is symmetric while $J$ is antisymmetric. The statement now follows by examining the summand whose numerator contains no $S_p$ for both valuations, which has different coefficients (depending on $p$ for a general horizontal distribution).

\endproof

\subsection{Extending $\phi_{2k}$ to arbitrary $2k$-submanifolds}
Define the contact area $\CA_{2k}(F)\in[0,\infty]$ for any $C^2$, $2k$-dimensional submanifold with boundary $F\subset M$ by 

\[\CA_{2k}(F) : = \liminf_{F_\epsilon\to F} \phi_{2k}(F_\epsilon)\]
where $F_\epsilon$ is a $C^2$, normally transversal submanifold with boundary that $C^2$-converges to $F$. 
Note that Proposition \ref{prop:DH_general_formula} and Theorem \ref{thm:weyl}, and the proof of Lemma \ref{lem:intermediate_DH}, imply that $\phi_{2k}(F)$ is $C^2$-continuous on normally transversal submanifolds $F$. Hence on such submanifolds, $\CA_{2k}(F)=\phi_{2k}(F)$.

We now show that the vanishing of $\CA_{2k}(F^{2k})$ is a necessary condition for the existence of an arbitrarily small perturbation with no contact points. For closed $F$ this supplements the topological necessary condition $\chi(F)=0$.

\begin{Proposition}
	For a $2k$-dimensional submanifold with boundary $F\subset M$, $\CA_{2k}(F^{2k})=0$ if and only if there is an arbitrarily small $C^2$-perturbation of $F$ which is nowhere tangent to the horizontal structure.
\end{Proposition}
\proof
The if direction follows from Corollary \ref{cor:general_positivity}. For the other direction, let us assume $\CA_{2k}(F)=0$. Again by Corollary \ref{cor:general_positivity}, $F$ cannot be normally transversal to the horizontal distribution unless it is nowhere tangent to it. Assuming the contrary to the assertion, we conclude $F$ is not normally transversal. Let $F_\epsilon\to F$, $\epsilon\to 0$ be a normally transversal family of smooth perturbations of $F$ s.t. $\phi_{2k}(F_\epsilon)\to 0$. If a normally transversal contact point exists for $F$ (which is then necessarily isolated), it persists to $F_\epsilon$ for small $\epsilon$, and examining eq. \eqref{eq:DH_general_formula}, one can find $c>0$ such that $\phi_{2k}(F_\epsilon)\geq c$ for all small $\epsilon$, a contradiction. Thus $F$ is not normally transversal to the horizontal distribution at any of its contact points. It follows that a sequence $x_j$ of contact points of $F_{\epsilon_j}$ must approach a necessarily degenerate tangency point, and hence by eq. \eqref{eq:DH_general_formula}, $\phi_{2k}(F_{\epsilon_j})\to \infty$. This again is a contradiction, implying that for small $\epsilon$, $F_\epsilon$ has no horizontal tangent spaces, as claimed.
\endproof
\section{The dual Heisenberg Lie algebra}
\label{sec:flat_DH}

\subsection{Linear algebra}

Recall that $U^{2n+1}$ denotes the dual Heisenberg algebra, with $H\subset U$ a fixed linear hyperplane, and $\omega\in \wedge^2 H^*\otimes U/H$ is a non-degenerate (twisted) form. It is the simplest DH manifold.

Let the group of automorphisms of the Heisenberg algebra be denoted by $\Sp_H(U)$ or $\Sp_H(2n+1)$. 
There is an $\Sp_H(U)$-equivariant isomorphism $\Dens(H)\simeq\Dens(U/H)^{n}$.

Define the subgroups $\Sp_H^+(U)=\{ g\in \Sp_H(U): g|_{U/H}>0 \}$, $\Sp^1_H(U)=\{ g\in \Sp_H(U): g|_{U/H}=\Id\}$ and $\textrm{Scal}_H(U)=\{g\in \Sp_H(U): g|_H\in \R^*\}\subset \Sp_H^+(U)$. Note that for $g\in\Sp_H^1(U)$ clearly $\det g =1$.
For an $\Sp_H(U)$-module $X$ and $x\in X$, we will write $\Stab(x)$, resp. $\Stab^+(x)$, $\Stab^1(x)$ for its stabilizer in the corresponding subgroup. 
Let us record the following trivial fact. 
\begin{Lemma}\label{lem:Sp_generators}
	$\Sp_H^+(U)$ is generated by  $\Sp^1_H(U)$ and $\emph{Scal}_H(U)$
\end{Lemma}
\begin{Corollary}
	Let $\delta_\lambda\in \emph{Scal}_H(U)$ act by $\lambda\in\R$ on $H$ and by $\lambda^2$ on a fixed direction $L$ complementing $H$. Then $\Sp_H^1(U)$ and $(\delta_{\lambda})_{\lambda\neq 0}$ generate $\Sp_H^+(U)$
\end{Corollary}
\proof

Take any $g\in \textrm{Scal}_H(U)$, set $\lambda:=g|_H$. Then $\delta_{\lambda}^{-1}\circ g\in \Sp_H^1$. The claim now follows from Lemma \ref{lem:Sp_generators}.
\endproof

\begin{Lemma}\label{lem:unrelated_S_T}
	For any subspace $E\subset H$, set $E_0=E\cap E^\omega$. Then one can find:
	\begin{enumerate}
		\item $S\in \Stab^1(E)$ s.t. $S(E_0)=E_0$ and $S|_{E_0}=2$.
		
		\item $T\in \Stab^+(E)$ s.t. $T(E_0)=E_0$, $T|_{E_0}=1$ and $\det T< 1$.
	\end{enumerate}
	
\end{Lemma}

\proof
i)
Decompose  $E\cap H=E_0\oplus F$, where $F\subset H$ is a non-degenerate subspace. Fix a vector $z\in E\setminus H$. There is then an induced symplectic form $\omega_z$ on $H$ given by $\omega_z(u,v)(z+H)=\omega(u,v)$. Define $S\in \GL(E)$  by setting $S|_{E_0}=2$, $S|_F=1$, and note that $S$ leaves $\omega_z|_E$ invariant. By Witt's extension theorem, we can find an extension  $S\in \Sp(H, \omega_z)$, and finally setting $S(z)=z$ yields $S\in\Stab^1(E)$ as required.
\\\\
ii)  Using $S$ from i), define $T\in \Sp_H^+(U)$ by $T|_H=\frac12 S|_H$, $Tz=\frac14 z$.

\endproof

The orbits $ Y^k_{\epsilon, r}$ of $\Sp_H^+(U)$ in $\Gr_k(U)$ are classified by the pairs $(\epsilon, r)$ as follows: \[E\in Y^k_{\epsilon, r}\iff \epsilon=k-\dim H\cap E\in \{0,1\},2r=\kappa-\dim\ker (\omega|_{E\cap H})\] where $\epsilon\in \{0,1\}$, and  $r\in \{ 0, \dots, \lfloor\kappa/2\rfloor\}$, $\kappa =\min (k-\epsilon, 2n-(k-\epsilon))$. . The unique open orbit has $\epsilon=1$ and $r=\lfloor\kappa/2\rfloor$, the unique closed orbit has $\epsilon=0$ and $r=0$.

\subsection{Translation invariant valuations on the dual Heisenberg algebra}

We start by classifying the $\Sp_H^+(U)$-invariant Klain sections.

\begin{Proposition}\label{prop:klain_upper_bounds} Fix $1\leq k\leq 2n$.
	\\	i) For even $k$, there is at most a one-dimensional space of $\Sp_H^+(U)$-invariant generalized Klain sections over $\Gr_k(U)$.
	\\ii) For odd k, there are no $\Sp_H^+(U)$-invariant generalized Klain sections over $\Gr_k(U)$.
\end{Proposition}

\proof We will make repeated use of Lemma \ref{lem:LocallyClosedOrbits} without explicit mention, wherein also the bundle $F_Y^\alpha$ of principal symbols transversal to $Y$ is defined.

Take $E\in Y^k_{\epsilon,r}$, and denote  $E_0=(E\cap H)\cap (E\cap H)^\omega$, $\dim E_0=\kappa-2r$.

\textit{Step 1.} $\epsilon=1$. We consider first the open orbit: $r=r_{\max}$. If $k$ is odd,  $E\cap H$ is non-degenerate, and $\omega$ gives an isomorphism $\Dens(E\cap H)=\Dens(U/H)^{\frac{k-1}{2}}$. Thus
\begin{align*}\Dens(E)&=\Dens(E/E\cap H)\otimes \Dens(E\cap H) = \Dens(U/H)\otimes \Dens(U/H)^{\frac{k-1}{2}}\\&=\Dens(U/H)^{\frac{k+1}{2}}=\Dens(U)^{\frac{k+1}{2(n+1)}}. \end{align*}
It follows that for odd $k$ there are no invariant sections on the open orbit.

For even $k$, let $E_0\subset E\cap H$ be the kernel of $\omega|_{E\cap H}$, which is a line. Then

\begin{align*}  \Dens(E)&=\Dens(E/E\cap H)\otimes \Dens(E\cap H) \\&=\Dens(U/H)\otimes \Dens(E_0)\otimes \Dens((E\cap H)/E_0)\\& = \Dens(U/H)^{\frac k 2} \otimes \Dens(E_0)=\Dens(U)^{\frac{k}{2(n+1)}}\otimes \Dens(E_0) \end{align*}

By Lemma \ref{lem:unrelated_S_T}, we may find $S\in\Stab^1(E)$ acting trivially on the first factor and rescaling the second factor. It follows that there are no invariant sections on the open orbit $Y^k_{1, r_{\max}}$.

Now let us consider the orbit $Y=Y^k_{1,r}$ with $r<r_{\max}$. We have $N_E Y=N_{E\cap H} X^{k-1}_r(H)=\wedge^2E_0^*$. Consider the bundle over $Y$ with fiber \[F_Y^\alpha|_E=\Dens(E)\otimes \Dens^*(N_E Y)\otimes \Sym^\alpha (N_E Y)
=\Dens(E)\otimes \Dens(\wedge^2 E_0)\otimes \Sym^\alpha(\wedge^* E_0)\]
Since $E/E\cap H$ is $\Stab(E)$-isomorphic to $U/H$, 

\[\Dens(E)=\Dens(E\cap H)\otimes \Dens(U/H)=\Dens((E\cap H )/ E_0)\otimes \Dens(E_0)\otimes \Dens(U/H)  \]
Now $\dim (E\cap H )/ E_0=2r$, and $\omega$ readily yields a non-degenerate form \[\tilde \omega\in \wedge^2((E\cap H )/ E_0)^*\otimes U/H\] so that there is a $\Stab(E)$-isomorphism $\Dens((E\cap H)/E_0)=\Dens(U/H)^{r}$. Thus

\[F^\alpha_Y|_E=\Dens(U/H)^{r+1}\otimes  \Dens(E_0) \otimes \Dens(\wedge^2 E_0)\otimes \Sym^\alpha(\wedge ^2 E_0^*),\]
Again by Lemma \ref{lem:unrelated_S_T}, there are no $\Stab^+(E)$-invariants in $F_Y^\alpha|_E$ when $r<r_{\max}$ (and so $E_0\neq \{0\}$). 

We conclude that for no $k$ are there invariant generalized sections whose support intersects $Y^k_{1,r}$, for any $r$.
We assume from now on that $\epsilon=0$, so $E\subset H$.

\textit{Step 2.} Consider $Y=Y^k_{0,r}$ with $r=r_{\max}$. Then $N_E Y=E^*\otimes U/H$, and \[F^0_Y|_E=\Dens(E)\otimes\Dens^*(N_EY) = 
\Dens(E)^2\otimes\Dens^*(U/H)^k
\]
Taking $g\in \textrm{Scal}_H(U)$ with $g|_H=\lambda$, we see that it acts on $F^0_E\otimes \Sym^\alpha(N_EY)=\Dens(E)^2\otimes\Dens^*(U/H)^k\otimes (U/H)^\alpha\otimes\Sym^\alpha E^*  $ by $\lambda^{-2k+2k+2\alpha-\alpha}=\lambda^\alpha$.  Thus $\alpha=0$ is the only possible transversal order of an invariant section. We now consider separately the different parities of $k$.

If $k$ is odd, set $E_0=E\cap E^\omega$, $\dim E_0=1$. Then
 
\[F_E^0= \Dens(E_0)^2\otimes \Dens(E/E_0)^2\otimes \Dens^*(U)^{\frac{k}{n+1 }}=\Dens(E_0)^2\otimes \Dens^*(U)^{\frac{k+1}{2(n+1)}  }  \]
since $\Dens(E/E_0)^2\simeq\Dens(U/H)^{k-1}= \Dens(U)^{\frac{k-1}{n+1}}$. Thus by Lemma \ref{lem:unrelated_S_T}, the action of $\Stab^+(E)$ on $F^0_E$ is clearly non-trivial, and so there are no invariant sections whose support intersects $Y$.

If $k$ is even, the restriction of $\omega$ to $E$ gives an isomorphism $\Dens(E)=\Dens(U/H)^{\frac k 2}$, so that $\Stab(E)$ acts trivially on $F_E^0$. We know by now that all invariant sections are supported inside $\overline Y$. Thus we conclude that the space of restrictions of the space of invariant sections to $\Gr_k(U)\setminus (\overline Y\setminus Y)$ is at most one-dimensional.

\textit{Step 3.} It remains to show there are no invariant sections supported on the closure of either of the orbits $Y_{0, r}^k(U)$ with $r<r_{\max}$, for any $k$. In particular, we have $1<k<2n$ and $\kappa-2r=\dim E_0\geq 2$, with $\kappa=\min(k, 2n-k)$.

One has the chain of inclusions $T_EY^k_{0,r}(U)=T_EX^k_r(H)\subset T_E\Gr_k(H)\subset T_E \Gr_k(U)$.  Hence $N_E Y^k_{0,r}(U)$ fits into the exact sequence 
\[0\to T_E\Gr_k(H)/T_EX^k_r(H)\to  N_E Y^k_{0,r}(U) \to T_E \Gr_k(U)/T_E\Gr_k(H)\to 0 \]
which is $\Stab(E)$-isomorphic to
\[0\to\wedge ^2 E_0^*\to  N_E Y^k_{0,r}(U) \to E^*\otimes U/H\to 0 \]
Write $Y=Y^k_{0,r}$. For $\alpha=0$, 
\begin{align*} F^0_Y|_E&=\Dens(E)\otimes \Dens^*(N_EY^k_{0,r}(U))\\&=\Dens(E_0)^2\otimes \Dens(E/E_0)^2 \otimes \Dens^*(U/H)^{k}\otimes \Dens^*(\wedge ^2 E_0^*)\\&=\Dens(E_0)^2\otimes \Dens^*(U/H)^{k-r}\otimes \Dens^*(\wedge ^2 E_0^*) . \end{align*}
Thus $ (F^0_Y|_E)^{\Stab^+(E)}=\{0\}$ since one can find $g\in\Stab^1(E)$ with $g|_{E_0}=\lambda\neq1$.

When $\alpha>0$, a $\Stab^+(E)$-invariant element of $F^\alpha_Y|_E$ would imply the existence of an invariant element in 
\[F^0_Y|_E\otimes (\wedge ^2 E_0^*)^a\otimes (E^*)^b\otimes (U/H)^b\]
which in turn implies the existence of an invariant element in

\[\Dens(E_0)^2\otimes \Dens(\wedge ^2 E_0)\otimes (\wedge ^2 E_0^*)^a\otimes (E_0^*)^{b'}\otimes  ((E/E_0)^*)^{b''}\otimes \Dens(U)^{\lambda}\]
for some non-negative integers $a, b, b', b''$ and $\lambda\in\R$.

By the proof of Lemma \ref{lem:unrelated_S_T}, we can find $S\in \Stab^1(E)$ such that $S|_{E_0}=2$ and $S:E/E_0\to E/E_0$ is the identity. Thus there are no invariants in this space.

\endproof

For a vector bundle $E$ over a manifold $B$, we let $\Gamma_m(B, E):=\Gamma_c(B, E^*\otimes |\omega_B|)^*$ denote the space of generalized sections that are given locally by a regular Borel measure.  

\begin{Proposition}\label{prop:klain_construction}
	There are $\Sp_H^+(U)$-invariant generalized Klain sections $\kappa_k$ realizing the upper bounds obtained in Proposition \ref{prop:klain_upper_bounds}. Moreover,  $\kappa_k\in\Gamma_m(\Gr_k(U), \Dens(E))$, they are supported on $Y^k_H:=\{E\subset H\}$ and are $\Sp_H(U)$-invariant.
\end{Proposition}
\proof
Take $f(E)=|\omega^{\wedge\frac k 2}|^{\otimes 2}\in \Dens(E)^2\otimes \Dens^*(U/H)^k$, which is a continuous section over $\Gr_k(H)$. We may rewrite $f$ as an absolutely continuous measure on $\Gr_k(H)$ with values in the bundle $\Dens^2(E)\otimes \Dens^*(U/H)^k\otimes\Dens^*(T_E\Gr_k(H))=\Dens^*(U)^k\otimes \Dens(E)^{2n+2}$. Writing $i:\Gr_k(H)\to \Gr_k(U)$ for the natural embedding, we get \begin{align*}i_*f &\in \mathcal M_{\Gr_k(H)}(\Gr_k(U), \Dens^*(U)^k\otimes \Dens(E)^{2n+2})\\&\simeq\Gamma_m(\Gr_k(U), \Dens^*(U)^k\otimes \Dens(E)^{2n+2}\otimes \Dens(T_E{\Gr_k(U)}))\end{align*}
The latter bundle is just the Klain bundle:
\[ \Dens^*(U)^k\otimes \Dens(E)^{2n+2}\otimes \Dens(E^*\otimes U/E)=\Dens(E) \]
and it remains to note that we only used $\Sp_H(U)$-equivariant identifications.

\endproof

Let us fix an involution $R\in \Sp_H(U)$ acting by $-1$ on $U/H$. Since $R^{-1}\Sp_H^+(U)R=\Sp_H^+(U)$, it follows that 
$R$ acts on the space of $\Sp_H^+(U)$-invariants in any $\Sp_H(U)$-module $M$. We will call the $\pm1$ eigenspaces of $R$ in $M^{\Sp_H^+(U)}$ $R$-even, resp. $R$-odd.

\begin{Proposition}\label{prop:DH_upper_bound}
	For $0\leq k\leq n$ it holds that $\dim \Val^{-\infty}_{2k}(U)^{\Sp_H^+(U)}\leq 2$ while $\dim \Val^{-\infty}_{2k+1}(U)^{\Sp_H^+(U)}=0$. For every $k$, the spaces of $R$-even and $R$-odd invariants are each at most one-dimensional. 
	Moreover, any $R$-even invariant valuation is even, and any $R$-odd invariant valuation is odd. 
\end{Proposition}
	In particular, any $\Sp_H(U)$-invariant valuation must be even.
\proof
\textit{Step 0.} In the following, all forms are translation-invariant. For $\xi \in \mathbb P_+(U^*)$, let $\xi_\perp\subset U$ be its annihilator. There is a natural identification of $\Omega_{-\infty}^{j,2n+1-j}(U\times \mathbb P_+(U^*))^{tr}$ with the generalized sections of the bundle over $\mathbb P_+(U^*)$ with fiber $\wedge^j U^*\otimes \wedge^{2n+1-j}\xi_\perp\otimes \xi^{2n+1-j}$ over $\xi$.
A Legendrian form corresponds to a section of the subbundle $\xi\otimes \wedge ^{j-1}(U^*/\xi)\otimes \wedge^{2n+1-j}\xi_\perp\otimes \xi^{2n+1-j}$, which we then call a Legendrian section. We will find the $\Sp^+_H(U)$-invariant generalized Legendrian sections. Let $\psi(\xi)$ be such a section. 

There are three orbits under $\Sp_H^+(U)$: the open orbit $X_o=\{\xi\neq H^\perp\}$ and two closed orbits $X_c^\pm=\{\pm H^\perp \}$.

\textit{Step 1.} Let us first show that $\psi$ vanishes when restricted to the open orbit. Note that $\psi|_{X_o}$ is smooth by $\Sp_H^+(U)$-invariance, and fix $\xi\in X_o$.

Consider the stabilizer $G_\xi:=\Stab^+(\xi)\subset \Sp_H^+(U)$, and $G_\xi^1=\{g\in G_\xi: g|_\xi=1\}$. Then $\psi(\xi)$ is a $G_\xi$-invariant element in $ \wedge ^{j-1}(U^*/\xi)\otimes \wedge^{2n+1-j}\xi_\perp\otimes \xi^{2n+2-j}$.
By considering the various invariant subquotients, we deuduce the existence of an invariant element in one of the following spaces:
\begin{align*}V_1&=(H^\perp\oplus \xi)/\xi \otimes \wedge^{j-2} (U^*/(\xi\oplus H^\perp))\otimes  \wedge^{2n+1-j} (\xi_\perp \cap H)\otimes \xi^{2n+2-j}
\\V_2&= (H^\perp\oplus \xi)/\xi \otimes \wedge^{j-2} (U^*/(\xi\oplus H^\perp))\otimes  \wedge^{2n-j}(\xi_\perp \cap H)\otimes (\xi_\perp/\xi_\perp\cap H)\otimes\xi^{2n+2-j}\\V_3&=  \wedge^{j-1} (U^*/(\xi\oplus H^\perp))\otimes  \wedge^{2n+1-j} (\xi_\perp \cap H)\otimes\xi^{2n+2-j}\\V_4&= \wedge^{j-1} (U^*/(\xi\oplus H^\perp))\otimes  \wedge^{2n-j}(\xi_\perp \cap H)\otimes (\xi_\perp/\xi_\perp\cap H)\otimes\xi^{2n+2-j}\end{align*}

Now take $g_\lambda\in G_\xi$ such that $g_\lambda|_H=\lambda$. Since $g_\lambda$ has a $2n$-dimensional eigenspace of eigenvalue $\lambda$ on $U$, the action of $g_\lambda$ on $U^*$, which is by $(g_\lambda^{-1})^*$, has a $2n$-dimensional invariant subspace of eigenvalue $\lambda^{-1}$. In the following we will simply write $g_\lambda$ for this action. We may choose $g_\lambda$ such that  $g_\lambda|_\xi=\lambda^{-1}$, $g_\lambda|_{H^\perp}=\lambda^{-2}$. Since $\xi_\perp/(\xi_\perp\cap H)\simeq U/H$, $g_\lambda$ acts on it by $\lambda^2$.

The action of $g_\lambda$ on $V_i$ is as follows: $g_\lambda|_{V_1}=\lambda^{-2-(j-2)+2n+1-j-(2n+2-j)}=\lambda^{-j-1}$, 
$g_\lambda|_{V_2}=\lambda^{-j}$, $g_\lambda|_{V_3}=\lambda^{-(j-1)+2n+1-j-(2n+2-j)}=\lambda^{-j}$, $g_\lambda|_{V_4}=\lambda^{-j+1}$.

Since $j\geq1$, we conclude that a $G_\xi$-invariant can only exist in $V_4$, and only if $j=1$. However, $\wedge^{2n}\xi_\perp \otimes \xi^{2n+1}$ has no  $G_\xi$-invariants: one may choose an element $g\in G_\xi$ with $g|_\xi=1$, while $\det g\neq 1$. Then $\det (g:\xi_\perp\to \xi_\perp) \neq 1$, and thus no such invariant exists. 

\textit{Step 2.} We conclude that $\psi$ is supported on $X_c:=\{\pm H^\perp\}$. Since $-\Id$ commutes with $\Sp_H^+(U)$, it acts on any space of $\Sp_H^+(U)$-invariants, which thus decomposes into a sum of even and odd invariants.
Now $\det(-\Id)=(-1)^{2n+1}=-1$, so that odd closed Legendrian forms correspond to even valuations and vice versa. 

Let us show that for odd $j$ there are no $\Sp_H^+(U)$-invariant closed Legendrian forms. Otherwise by what we proved, there is an $\Sp_H^+(U)$-invariant such form $T$ supported on $\xi=H^\perp$ with a fixed orientation. Then $T_s:=T-(-\Id)^*T$ is odd, closed, Legendrian and $\Sp_H^+(U)$-invariant. Thus $T_s$ defines an even $\Sp_H^+(U)$-invariant generalized valuation which is $j$-homogeneous. This contradicts Proposition \ref{prop:klain_upper_bounds}.

\textit{Step 3.}
The principal symbol of $\psi$, denoted $F^\alpha_{X_c}|_\xi$  (see Appendix \ref{app:invariants}) over each point $\xi\in X_c$  is an element of \begin{align*}W&:= H^\perp\otimes \wedge^{j-1}H^* \otimes \wedge ^{2n+1-j}H\otimes (H^\perp)^{2n+1-j}\otimes \Dens^*(T_\xi \mathbb P_+(U^*))\otimes \Sym^\alpha(T_\xi \mathbb P_+(U^*)) \\&=\wedge^{j-1}H^*\otimes \wedge ^{2n+1-j}H\otimes (H^\perp)^{2n+2-j-\alpha}\otimes \Dens(H)\otimes \Dens(H^\perp)^{2n}\otimes\Sym^\alpha(H^*) . \end{align*}

Take $\delta_\lambda\in\textrm{Scal}_H(U)$ acting by $\lambda \in \R$ on $H$ and by $\lambda^2$ on some fixed vector $w\in U\setminus H$. Then $\delta_\lambda$ acts on $H^\perp$ by $\lambda^{-2}$, and

\[\delta_\lambda|_{W}= \lambda^{-(j-1)}\lambda^{2n+1-j}\lambda^{-2(2n+2-j-\alpha)}\lambda^{-2n}\lambda^{4n}\lambda^{-\alpha}=\lambda^{\alpha-2}\] hence $\alpha=2$ is necessary for an invariant to exist. 

We may identify $\wedge^{2n+1-j}H\simeq \wedge^{j-1}H^*\otimes \wedge^{2n}H$. Thus 
\[W= \wedge^{j-1}H^*\otimes \wedge^{j-1}H^*\otimes \Sym^2(H^*) \otimes\wedge^{2n}H\otimes (H^\perp)^{2n-j}\otimes \Dens(H)\otimes \Dens(H^\perp)^{2n} \]
Let us find all $w\in  W$ that are invariant under $\Sp^1_H(U)$. We may fix an $\Sp^1_H(U)$-invariant symplectic form $\omega_H$ on $H$, so that $\Sp^1_H(U)$-equivariantly, 
\[W\simeq \wedge ^{j-1}H^*\otimes \wedge ^{j-1}H^*\otimes \Sym^2 (H^*).\]

By the fundamental theorem of invariant theory, an invariant element of $W$ is given by fixing pairings of all the factors using $\omega$, and then symmetrizing/antisymmetrizing accoridngly (some pairings could give zero). There would be $j$ pairings. 

Note that $R^*\omega_H=-\omega_H$. Hence $\det R|_H=(-1)^n$ and $\det R=(-1)^{n+1}$. We see that $Rw=(-1)^j(-1)^n(-1)^j=(-1)^nw$ for any $w\in W^{\Sp^1_H(U)}$.

\textit{Step 4.}
Now take $\phi\in\Val^{-\infty}_j(U)^{\Sp_H^+(U)}$ which is either even or odd, and denote $T=T(\phi)$.
Assume first that $\phi$ is $R$-even, that is $\Sp_H(U)$-invariant.
Then $R^*T=\det \cdot T$. Thus its principal symbol is $\sigma(T)(+H^\perp, -H^\perp)=(w_0, -w_0)$ for some $w_0\in W^{\Sp^1_H(U)}$. 
Noting that $-\Id$ acts by $(-1)^j$ on $W$, it follows that $\phi$ is an even valuation if and only if $j$ is even, which we may assume by step 2. By Proposition \ref{prop:klain_upper_bounds}, $\dim \Val_j^{-\infty}(U)^{\Sp_H(U)}\leq 1$.

Similarly if $\phi$ is $R$-odd, its principal symbol is $\sigma(T)=(w_0, w_0)$ for some $w_0\in W^{\Sp^1_H(U)}$, and since $j$ is even, $\phi$ must be odd. It remains to prove the dimension of the space of $R$-odd, $j$-homogeneous valuations is at most $1$. It suffices to show by Theorem \ref{thm:current_injective} that the space of $\Sp_H^+(U)$-invariant, $R$-odd, $j$-homogeneous, closed Legendrian generalized forms supported on $X_c$ is at most $1$-dimensional. Any such form has the form $T_+-R^*T_+$ with $T_+$ supported on $+H^\perp$ and $\Sp_H^+(U)$-invariant. Clearly $T_+-R^*T_+\mapsto T_++R^*T_+$ is a bijection onto the corresponding $R$-even forms, which in turn bijectively correspond to $\Sp_H(U)$-invariant valuations. The $R$-even case now concludes the proof.

\endproof

The following proposition completes the proof of Theorems \ref{mainthm:DH_Hadwiger} and \ref{mainthm:DH+_Hadwiger}. 
 
\begin{Proposition}\label{prop:DH_space_Hadwiger}
For $0\leq k\leq n$, it holds that $\dim\Val^{-\infty}_{2k}(U)^{\Sp_H(U)}=1$ while $\dim\Val^{-\infty}_{2k}(U)^{\Sp_H^+(U)}=2$.
\end{Proposition}
\proof
The first assertion follows from the general construction of valuations on DH manifolds, as follows. By Proposition \ref{prop:valuations_on_DH_manifold}, we have at least $n+1$ linearly independent valuations $\phi_{2j}\in \mathcal W _{2j}^{-\infty}(U)$ which are invariant under the symmetries of $U$ as a DH-manifold, in particular under translations and $\Sp_H(U)$. Now Proposition \ref{prop:DH_upper_bound} concludes the proof.

The second assertion follows from the first one: by the last paragraph of the proof of Proposition \ref{prop:DH_upper_bound}, there is a bijection between the $\Sp_H(U)$-invariant valuations and the $R$-odd $\Sp_H^+(U)$-invariant valuations.
\endproof

\begin{Remark}
	1. Alternatively, one may use the Corfton valuations $\psi_{2j}$ on the contact sphere constructed in section \ref{sec:contact_sphere} to prove the statement, by considering the translation-invariant valuations they define on the tangent spaces by Theorem \ref{thm:filtration_quotients}.
	\\	
	2. Examining the construction of $\phi_{2k}$ for general DH manifolds, we immediately see that $\phi^U_{2k}$ is in fact $2k$-homogeneous.
\end{Remark}

\begin{Example}
	Let us describe explicitly the valuation $\phi_{2n}$. Consider $U=\R^{2n+1}$ with the standard Euclidean structure, and take $H=\{x_{2n+1}=0\}$. For $u,v\in H$ we let $\omega(u,v)=\sum_{j=1}^n (u_{2j-1}v_{2j}-u_{2j}v_{2j-1})$ be the standard symplectic form.
	Identifying $U/H=H^\perp=\R e_{2n+1}=\R$ we recover the DH structure on $U$.
	
	Examining the proof of Proposition \ref{prop:klain_construction}, we see that for convex $K\in\mathcal K(V)$, $\phi_{2n}(K)=\frac12(\sigma_K(H^\perp)+\sigma_K(-H^\perp))$, where $\sigma_Kd\theta\in\mathcal M(S^{2n})$ stands for the surface area measure of $K$, which is essentially just the push-forward of the $2n$-dimensional Hausdorff measure on $\partial K$ to $S^{2n}$, see \cite{schneider} for exact definition.

\end{Example}

\subsection{Homogeneity $2n$}
The valuation of homogeneity $2n$ in $U$, which is just the surface area measure at the points of tangency to $H\subset U$, has also a different natural setting, valid in dimension of any parity.

Let $H^m\subset V^{m+1}$ be a hyperplane, and let $\vol_H$, $\vol_V$ be fixed volume forms on $H$, resp. $V$. Consider $\SL_H(V):=\{g\in \SL(V): g(H)=H,g|_H\in\SL(H)\}$. Fix a Euclidean product $P$ on $V$ inducing the given Lebesgue measures on $V$ and $H$. Assume for simplicity $m\geq 2$.
	\begin{Proposition}
		$\phi_{m}(K):=\frac12(\sigma_K(H^\perp)+\sigma_K(-H^\perp))$ is the unique $m$-homogeneous, even $\SL_H(V)$-invariant generalized valuation.
	\end{Proposition}	
\proof The Klain section of $\phi_{m}$ is a delta-measure on $\Gr_m(V)$ supported on $H$. It is given by an element of $\Dens(H)\otimes \Dens^*(T_H\Gr_{m}(V))=\Dens(H)^{m+2}\otimes\Dens^*(V)^{m}$. Clearly $\SL_H(V)$ acts trivially on this space.

The proof of uniqueness is similar to Proposition \ref{prop:DH_upper_bound}. We first note that there are no invariant sections of the Klain bundle over the open orbit $\{E: E\neq H\}\subset\Gr_m(V)$, since for such $E$ there is a natural isomorphism $E/E\cap H=V/H$, so a $\Stab(E)$-invariant density on $E$ would produce a $\Stab(E)$-invariant density on $E\cap H$, which clearly does not exist since $\Stab(E)$ can rescale this space.
Next we use Lemma \ref{lem:LocallyClosedOrbits} to study Klain sections supported on the closed orbit $\{H\}$: for $\alpha\geq 0$, \[F^\alpha_{\{H\}}=\Dens(H)\otimes \Dens^*(H^*\otimes V/H)\otimes \Sym^\alpha(H^*\otimes V/H) \] which is $\SL_H(V)$-equivariantly isomorphic to $\Sym^\alpha(H^*)$. There are no invariant polynomials on $H$, so we must have $\alpha=0$, which readily yields a one-dimensional space of invariants. 
\endproof
		 
	Next, we write an invariant Crofton formula for this valuation.
	\begin{Proposition}\label{prop:crofton_gaussian}
	There is an invariant Crofton measure $\mu_{\SL}$ over $\Gr_1(V)$ which defines $\phi_{m}$.  
	\end{Proposition}

\proof

Let $\theta(E)\in [0,\pi]$ be the Euclidean angle between $E\in \Gr^+_1(V)$ and $H^\perp$ with fixed orientation. Define the meromorphic family of functions $f_s(E)=(\cos \theta)^*|t|^s$ where $|t|^s\in C^{-\infty}(\R)$ is the standard meromorphic family of even homogeneous generalized functions. Note that $\cos\theta$ is submersive whenever $t=\cos\theta=0$, so the pull-back is well-defined. Since $f_s(E)$ is invariant to orientation reversal, we get a generalized family on $\Gr_1(V)$, still denoted $f_s(E)$.

Let us identify translation-invariant measures (distributions) on $\AGr_1(V)$ with (generalized) sections of the Crofton bundle $\Cr_{m}$ over $\Gr_1(V)$, whose fiber over $E\in \Gr_1(V)$ is $\Dens(V/E)\otimes \Dens(T_E\Gr_1(V))$.

Define a generalized section of $\Cr_{m}$ over $\Gr_1(V)\setminus \Gr_1(H)$ by $\mu_{\SL}:=f_{-m-2}(E)$ when $m$ is even, and $\mu_{\SL}:=\Res_{s=-m-2}f_s$ when $m$ is odd.

Let us check $\mu_{\SL}$ is $\SL_H(V)$-invariant. 
For $g\in \GL(V)$ consider the Jacobian $\psi_g(E)=\textrm{Jac}(g:E\to gE)^{-2}$, where $E, gE$ are endowed with the volume induced by the Euclidean product $P$. It follows that for $g\in \SL_H(V)$ and $E\not\subset H$, \[\psi_g(E)=\frac{\cos^2 \angle (g E, H^\perp)}{\cos^2 \angle (E, H^\perp)}.\]  Now $f(E)$ represents (with respect to the Euclidean trivialization) an $\SL_H(V)$-invariant Crofton measure precisely if   
$g^*f(E)=\psi_g(E)^{-\frac{m+2}{2}}f(E)$. Since $g^* f_s=\psi_g(E)^{\frac{s}{2}} f_s(E)$, we conclude $\mu_{\SL}$ is indeed $\SL_H(V)$-invariant.

It remains to verify these Crofton measures define non-zero valuations. For this, we evaluate $\phi_{m}(B)$ for the unit ball $B^{m+1}$.
By definition, $\phi_{m}(B)=1$. On the other hand, writing $\omega_m$ for the volume of the Euclidean ball $B^{m}$, we get 

\begin{align*} &\Cr(f_s(E)dE)(B)= \omega_m \int_{\Gr_1(V)}  |\cos\theta(E)|^s dE=\omega_m  \frac{ \int_0^\pi |\cos\theta|^s\sin^{m-1}\theta d\theta}{\int_0^\pi \sin^{m-1}\theta d\theta}\\&=\omega_m   \frac{\int_0^1 t^s (1-t^2)^{\frac{m-2}{2}}dt }{\int_0^1(1-t^2)^{\frac{m-2}{2}}dt}=\omega_m \frac{\int_0^1 u^{\frac{s-1}{2}} (1-u)^{\frac{m-2}{2}}du }{\int_0^1 u^{-\frac12} (1-u)^{\frac{m-2}{2}}du} =\omega_m\frac{B(\frac{s+1}{2}, \frac m 2)}{B(\frac12, \frac m2)} \end{align*}
Thus \[C_m=\frac{\pi^{\frac m 2}}{\Gamma(1+\frac m 2)}  \frac{\Gamma(\frac{s+1}{2})}{\Gamma(\frac{s+1}{2}+\frac m 2)}\frac{\Gamma(\frac m2+\frac12) }{\Gamma(\frac12)}=\frac{\pi^{\frac m 2-\frac12}}{\Gamma(1+\frac m 2)}  \frac{\Gamma(\frac{s+1}{2})\Gamma(\frac m2+\frac12)}{\Gamma(\frac{s+1}{2}+\frac m 2)}\]
When $m=2n$ is even, $s=-2n-2$ and  \[C_{2n}=-\frac12 \frac{\pi^{n-1}}{n!}\Gamma(n+\frac12)\Gamma(-n-\frac12)\]
When $m=2n-1$ is odd, $s=-2n-1$, $\Gamma(\frac{1+s}{2})$ has a simple pole with residue $2\Res_{z=-n}\Gamma(z)=2(-1)^n\frac{1}{n!}$, and 
\[ C_{2n-1}=2(-1)^n\frac{1}{n!}\frac{\pi^{n-1}}{\Gamma(n+\frac 12)} \frac{\Gamma(n)}{\Gamma(-\frac12)}=(-1)^{n-1}\frac{\pi^{n-\frac32}}{n\Gamma(n+\frac 12)} \] 
\endproof

We can now recover effortlessly the inverse form of the Koldobsky-Ryabogin-Zvavitch formula  \cite {KRZ}.
\begin{Corollary}
	Let $K\subset \R^{m+1}$ be a smooth, symmetric convex body, and let $\kappa(x)$ be its Gaussian curvature and $\nu(x)$ the unit normal at $x\in\partial K$. Then there is a universal explicit constant $C_m$ (computed in the proof of Proposition \ref{prop:crofton_gaussian}) s.t. for even $m=2n$, \[\kappa(x)^{-1}=C_m ^{-1}\int_{S^{m}}\vol_{m}({\Pr}_{\theta^\perp}(K))\frac{d\sigma_1(\theta)}{\langle \theta,\nu(x)\rangle^{m+2}} \] 
	while for odd $m$, 
	\begin{align*}\kappa(x)^{-1}&=C_m^{-1} \left\langle \delta^{(m+1)}(\langle \nu(x), \theta\rangle),\vol_{m}({\Pr}_{\theta^\perp}(K))\right\rangle \\&=\frac{2}{(m+1)!} C_m^{-1}\left.\frac{d^{m+1}}{dt^{m+1}}\right|_{t=0}\left((1-t^2)^{\frac{m-2}{2}}\int_{\theta:\langle \theta, \nu(x)\rangle=\cos t}\vol_m(\textrm{Pr}_{\theta^\perp}(K))d\sigma_1(\theta)\right)  \end{align*}
	
\end{Corollary}
\begin{Remark}
Since $\vol_m(\Pr_{\eta^\perp}(K))$ is the cosine transform of $\sigma_K(\theta)d\sigma_1(\theta)$, one can think of this formula as simply the inversion of the cosine transform. The integral can be viewed as the $(-2n-2)$-cosine transform of the support function of the projection body.

\end{Remark}

\section{Contact manifolds}\label{sec:contact}
\subsection{Specializing from general DH manifolds.}
Recall \ref{exm:contact_is_DH} that a contact manifold has a canonical structure of a DH manifold. Let us recall a few basic facts.

The standard contact structure on $\R^{2n+1}=\{(x_1,\dots, x_n, y_1, \dots, y_n,z)\}$ is given by $dz=\sum_{j=1}^n x_jdy_j-y_jdx_j$.
Its symmetry group contains translations along the $z$-axis and rotations in the $(x_j, y_j)$ planes.

\begin{Theorem}[Darboux-Pfaff] Take a contact manifold $M^{2n+1}$. Then any $x\in M$ has a neighborhood contactomorphic to an open subset of the standard contact $\R^{2n+1}$.
\end{Theorem}

We will make use of the following simple observation.
\begin{Lemma} \label{lem:standard_contact_convex}
	In the standard contact space $\R^{2n+1}$, it holds that a closed, smooth, strictly convex hypersurface $F\subset \R^{2n+1}$ is normally transversal.
\end{Lemma}
\proof
We will use the standard Euclidean structure to identify $\mathbb P_{\R^{2n+1}}$ with $\R^{2n+1}\times S^{2n}$.
Let $p$ be a contact point of $F$. We may assume $p=(a_1,0,\dots, a_n,0,0)$. As $T_pF=H_p=\{dz=\sum a_j dy_j\}$, $(x_j, y_j)_{j=1}^n$ form a system of coordinates for $F$ near $p$, and so locally $F$ is the graph $z=f(x_1,\dots, y_n)$ of a strictly convex or strictly concave function $f$. We will denote $w=w(x)=(x_1,\dots, y_n)$, $x=(w,f(w))$.
Let $\nu_w=\nu_x$ be the normal to $F$ at $x$, and $h_w=h_x$ the normal to $H_{x}$, both normalized to have $z$ coordinate $1$ (the other pair of normals has $z=-1$ and is treated identically). Thus $\nu_w=(-f'_{x_1},\dots, -f'_{y_n}, 1)$, and $h_w=(y_1, -x_1,\dots, y_n, -x_n, 1)$. 

\textit{Claim.} There exists $c=c(f,p)>0$ such that $\|\nu_w-h_w\|\geq c\|w-w(p)\|$ in a neighborhood of $p$.

\textit{Proof.} Write 
\[E=\|\nu_w-h_w\|= \sum( x_j-f'_{y_j})^2+ (y_j+f'_{x_j})^2\]
Replacing $f$ with $f-\sum a_j y_j$ and $x_j$ with $x_j-a_j$, $E$ remains unchanged, and the convexity of $f$ is retained. We thus may assume $p=0$, and $\nabla f(p)=0$. Hence
\[E=\|w\|^2+\|\nabla f(w)\|^2-2\sum_j (x_jf'_{y_j}(w)-y_jf'_{x_j}(w))=\|w\|^2+\|\nabla f(w)\|^2-2\omega(w,\nabla f(w))  \]
where $\omega(w,w')=\sum (x_jy'_j-x'_jy_j)$ is the standard symplectic form on $\R^{2n}$. We clearly may replace $f$ with its quadratic approximation at $p$, namely $f(w)=\frac12\langle Aw, w\rangle $ for $A=\textrm{Hess}f(p)$. Then $\nabla f(w)=Aw$, and by Writinger's inequality $|\omega(w,\nabla f(w))|\leq \|w\|\|A w\|\sin\angle (w, Aw)$. We conclude that 
\begin{align*}E&\geq\|w\|^2+\|Aw\|^2-2\|w\|\|Aw\|\sin\angle(w,Aw)\\&=(\|w\|-\|Aw\|)^2+2\|w\|\|Aw\|(1-\sin\angle(w,Aw)) \end{align*}
Since $A$ is sign-definite, it holds that $\lambda_{\min}\|w\|\leq\|Aw\|\leq \lambda_{\max}\|w\|$ and $|\langle Aw,w\rangle|\geq \lambda_{\min}\|w\|^2$,  
where $\lambda_{\min\textrm{/}\max}$ is the minimal/maximal eigenvalue of $|A|$. It follows that $|\cos\angle (w, Aw)|=\frac{|\langle w, Aw\rangle|}{\|w\|\|Aw\|}\geq\frac{\lambda_{\min}}{\lambda_{\max}}$, so $1-\sin\angle(w,Aw)>1-\sqrt{1-\frac{\lambda_{\min}^2}{\lambda_{\max}^2}}=:c_0$, and hence
$E \geq 2c_0\lambda_{\min} \|w\|^2$ as claimed.
\endproof
Recall $M_H=\{(x, h_x):x\in \R^{2n+1}\}$ is $2n+1$ dimensional, and $NF=\{(x,\nu_x): x\in F\}$ has dimension $2n$. Thus we ought to show $T_{(p, \nu_p)}M_H\cap T_{(p, \nu_p)}NF=\{0\}$. If not, there is a curve $x(t)\subset F$, $x(0)=p$ with non-zero velocity vector, such that $h_{x(t)}$ and $\nu_{x(t)}$ have equal velocity vectors at $t=0$. But then we should have $\|h_{x(t)}-\nu_{x(t)}\|=O(t^2)$ as $t\to 0$, contradicting the claim. \qed

\begin{Proposition} \label{prop:contact_linear_dependencies}
	Assume $M^{2n+1}$ is a contact manifold. Then $\phi_{2j}$, $0\leq j\leq n$ are all linearly independent, while $\phi_{2j-1}$ is a linear combination of $\phi_{2i}$, $j\leq i\leq n$.
\end{Proposition}
\proof The first statement is just part i) of Proposition \ref{prop:valuations_on_DH_manifold}. For the second statement, we will use the same notation as in the proof of Proposition \ref{prop:DH_general_formula}.

Write $h=(h_{ij})_{i,j=1}^{2n}=H_s+H_a$, the symmetric and antisymmetric part.

First, we note that in a contact manifold, the relation between the symplectic and horizontal structures is given by $H_a=-J$. Indeed, $\alpha(v)=g(R, v)$ so we have

\begin{align*}d\alpha(u,v)&=L_u\alpha(v)-L_v\alpha (u)-\alpha([u,v])=L_u g(R,v)-L_vg(R, u)-\alpha([u,v])\\&=g(\nabla _u R, v)+g(R, \nabla_u v)-g(\nabla_v R, u)-g(R, \nabla_v u)-g(R, [u,v])\\&=g(\nabla_u R, v)-g(\nabla_v R, u). \end{align*}

Hence $h_{ji}-h_{ij}=\theta_j(\nabla_{X_i}R)-\theta_i(\nabla_{X_j}R)=d\alpha(X_i, X_j)=J_{ij}$.

Take $F$ a closed hypersurface which is normally transversal. Introduce the notation $\Psi_{k}(F, p):={2n \choose k}^{-1}|\det A_p| \phi_{2n-k}(F, p)$. Then

\begin{align*}\Psi_k(F, p)&=D(H-S[k], J[2n-k])=(-1)^kD(S-H_s+J[k], J[2n-k]\\&= (-1)^k\sum_{i=0}^{k}{k\choose i}D(S-H_s[i], J[2n-i])  \end{align*}
Observe that for $2n\times 2n$ matrices $X, Y$ for which $X^T=X$, $Y^T=-Y$ it holds that $D(X[2n-i], Y[i])=0$ for odd $i$. We conclude  
\[\Psi_k(F, p)= (-1)^k\sum_{i=0}^{\lfloor\frac k2\rfloor}{k\choose 2i}D(S-H_s[2i], J[2n-2i]) \]
that is,
\begin{equation}\label{eq:phi_relations} \phi_k(F,p)=(-1)^k{2n \choose k}|\det A_p|^{-1}\sum_{i=0}^{\lfloor\frac {2n-k}2\rfloor}{2n-k\choose 2i}D(S-H_s[2i], J[2n-2i]) \end{equation}
Thus for normally transversal closed hypersurface $F$, $\phi_{2j-1}$ is a fixed linear combination of the $\phi_{2i}$ with $j\leq i\leq n$.
By Darboux-Pfaff, we may cover $M$ by charts $U_\alpha$ that are  contactomorphic to open subsets  $V_\alpha$ of the standard contact space $\R^{2n+1}$. By Lemma \ref{lem:standard_contact_convex}, a closed, smooth, strictly convex hypersurface in $V_\alpha$ is normally transversal to the contact structure. The statement now follows from Lemma \ref{lem:small_wavefront}, Proposition \ref{prop:euler_verdier_one}, and Lemma \ref{lem:hypersurface_suffices}.

\endproof

We now prove a Hadwiger theorem for contact manifolds.

\begin{Theorem}\label{thm:contact_valuations_classification}
	Let $M^{2n+1}$ be a contact manifold. Then $\mathcal V^{-\infty}(M)^{Cont(M)}$ is spanned by $\phi_{2k}$, $0\leq k\leq n$.
\end{Theorem}
\proof

For an element $\phi\in\mathcal W_k^{-\infty}$, we will write $[\phi]_k$ for its image in $\mathcal W_k^{-\infty}/\mathcal W_{k+1}^{-\infty}$.
We will show by induction on $2n+1-k$ that $\mathcal W_k^{-\infty}(M)^{Cont(M)}$ is spanned by $(\phi_{2j})_{2j\geq k}$. 
Take $\phi\in \mathcal W_{2n+1}^{-\infty}(M)^{Cont(M)}=\mathcal M^{-\infty}(M)^{Cont(M)}=\mathcal M^{\infty}(M)^{Cont(M)}$,  where the latter equality holds since $Cont(M)$ acts transitively on $M$, so that an invariant distribution is automatically smooth.
There is no such invariant measure, hence $\phi=0$. 

Assume now that  $\phi\in \mathcal W_{2n+1-k}^{-\infty}(M)^{Cont(M)}$. 
Recall \[\mathcal W_{2n+1-k}^{-\infty}(M)/\mathcal W_{2n+2-k}^{-\infty}(M)=(\mathcal W_{k,c}^{\infty}(M)/\mathcal W_{k+1,c}^{\infty}(M))^*=\Gamma_c (M, \Val_{k}^\infty(TM))^*\]
Hence by Proposition \ref{prop:infinite_dimensional_generalized}, \[[\phi]_{2n+1-k}\in(\Gamma_c (M, \Val_{k}^\infty(TM))^*)^{Cont(M)}=\Gamma^\infty(M, \Val_{2n+1-k}^{-\infty}(TM))^{Cont(M)}=\] where we used the Alesker-Poincar\'e isomorphism 
$\Val^\infty_{k}(T_xM)^*=\Val_{2n+1-k}^{-\infty}(T_xM)\otimes \Dens(T_xM)$.
For a point $x\in M$, its stabilizer in $Cont(M)$ is $\Stab(x)=\Sp_H(T_x M)$. Hence 
$\Gamma^\infty(M, \Val_{2n+1-k}^{-\infty}(TM))^{Cont(M)}=\Val_{2n+1-k}^{-\infty}(T_xM)^{\Sp_H{T_xM}}$.
By Proposition \ref{prop:DH_space_Hadwiger}, the latter space of invariants is trivial if $k$ is even, and one dimensional if $k$ is odd. 

Thus if $k$ is even, $\phi\in \mathcal W_{2n+2-k}^{-\infty}(M)$. In case $k$ is odd, we use part (i) of Proposition \ref{prop:valuations_on_DH_manifold} to find a multiple of $\phi^M_{2n+1-k}$ such that $[\phi]_{2n+1-k}=c[\phi^M_{2n+1-k}]_{2n+1-k}\Rightarrow\phi-c\phi^M_{2n+1-k}\in \mathcal W_{2n+2-k}^{-\infty}(M)^{Cont(M)}$. The induction assumption now completes the proof.
\endproof

\begin{Remark}
	Unlike the contact case, we do not have a uniqueness result in the DH category, where the symmetry group is in general trivial. In the Riemannian setting, uniqueness of isometry-invariant valuation assignment can be deduced from the classical Hadwiger theorem in conjunction with the Nash embedding theorem. This last piece is missing in the DH setting. A different type of uniqueness in terms of the Cartan frame apparatus, which is again tailored to the Riemannian setting, was established by Fu-Wannerer \cite{fu_wannerer}.
	
\end{Remark}

\subsection{A dynamical point of view}

Let $F\subset M^{2n+1}$ be a hypersurface, and let $p\in F$ be a normally transversal contact point: $T_pF=H_p$. Denote by $F_H$ the singular hyperplane field on $F$ given by $F_H|_p=H_p\cap T_pF$. When $\dim M=3$, this field integrates to the characteristic foliation.

One can describe $\phi_{2k}(F, p)$ explicitly through the singular bundle $H_x$ near $x$.

Let $\beta\in\Omega^1(F)$ be a form defined near $p$ such that $\Ker \beta=F_H$. Since $M$ is contact, we may assume $d\beta\neq 0$ near $p$ (e.g. by taking $\beta=\alpha|_F$ for some contact form $\alpha$ on $M$), and there is a unique vector field $B\in \mathfrak X(F)$ near $p$ such that $i_Bd\beta=\beta$. 
In particular, $B(p)=0$ and $B(x)$ is tangent to the characteristic foliation.
If $\beta'=f\beta$ is a different form with $d_p\beta'\neq 0\iff f(p)\neq 0$, the corresponding vector field is 
\[ B' = \left(1+\frac{df(B)}{f}\right)^{-1}B  \]
Since $df(B)(p)=df(B(p))=0$, the differential $d_pB\in \mathfrak{gl}(T_pB)$ only depends on $F_H$. 

\begin{Remark}
	Note  that $\sign \det d_pB$ determines whether $p$ is an elliptic or hyperbolic singular point of the characteristic foliation.
\end{Remark}

\begin{Proposition}\label{prop:contact_curvature_formula}
	$F$ is normally transversal at $p$ if and only if $d_pB$ is non-singular. In that event	

	\[\phi_{k}(F, p) = |\det d_pB|^{-1} \tr(\wedge ^{2n-k} d_p B)\]
\end{Proposition}
\proof
Again we use notation from the proof of Proposition \ref{prop:DH_general_formula}.

Since all contact manifolds are locally isomorphic, we work in $\R^{2n+1}$ with coordinates $(x_1,\dots, x_n, y_1, \dots, y_n, z)$, and contact form $\alpha=-dz+\sum x_jdy_j$. We may assume further that $p=0$ and $T_pF=\{z=0\}$. Then $(x_j, y_j)$ are local coordinates on $F$, and $\beta=\sum -\frac{\partial f}{\partial x_j} dx_j + (x_j-\frac{\partial f}{\partial y_j}) dy_j$, $d\beta|_0=\sum dx_j\wedge dy_j$. It follows that $B(x,y)^T=(x,0)^T+J\nabla f$, where $\nabla f= (\frac{\partial f}{\partial x_j}, \frac{\partial f}{\partial y_j})^T$. Thus \[d_0B= \left(\begin{array}{cc}I_n & 0\\ 0 & 0\end{array}\right) +J H^2f\]

On the other hand, one immediately computes that \[h= \left(\begin{array}{cc}0 & 0\\ I_n & 0\end{array}\right)\] and $S=H^2f$, hence $h-S=Jd_0B$. This readily shows that $F$ is normally transversal at $p$, which is equivalent to the non-singularity of $h-S$, if and only if $\det d_p B\neq 0$; and \[D(h-S[2n-k], J[k])=\det J\cdot D(d_0B[2n-k], I_n[k])={2n \choose k}^{-1}\tr(\wedge^{2n-k}(d_0 B)).\]
By Definition \ref{def:local_areas}, we are done.
\endproof

Computatiton of $\phi_{2k}$ is straightforward with this approach. Here is a simple proof of a well-known fact.
\begin{Corollary}
In the standard contact space $\R^{3}$ with contact structure given by $dz=xdy-ydx$, spheres of different radii are not equivalent through a contactomorphism of the ambient space.
\end{Corollary}
\proof
One computes that $\phi_2(S_R)=8(1+\frac14 R^{-2})^{-1}$.
\endproof

\begin{Example}[The contact sphere]{\emph 
Let us compute $\phi_{2k}(S^{2m})$ in $S^{2n+1}$. By Theorem \ref{thm:weyl}, we may assume $n=m$. Consider $S^{2m+1}\subset \C^{m+1}$, with coordinates $x_1,y_1, \dots, x_{m+1}, y_{m+1}$. The contact form is given by \[\alpha_p(v)=\langle v, \sqrt{-1}p\rangle=\sum_{j=1}^{m+1}(-y_jdx_j+x_jdy_j).\] so that $d\alpha=2\sum_{j=1}^{m+1} dx_j\wedge dy_j$. 
Fix $S^{2m}=\{y_{m+1}=0\}$. Then the two unique contact points of $S^{2m}$ are given by $x_{m+1}=\pm1$, and we use the coordinates $(x_j, y_j)_{j=1}^m$ near those points.
In those coordinates, $\beta=\alpha|_{S^{2m}}=\sum_{j=1}^{m}(-y_jdx_j+x_jdy_j)$, and $d\beta=2\sum_{j=1}^{m} dx_j\wedge dy_j$.
Then 
\[i_Bd\beta=\beta\Rightarrow 2\sum (dx_j(B)dy_j-dy_j(B)dx_j)=\sum (x_jdy_j-y_j dx_j) \]
so that $B(p)=\frac{1}{2}p$. Thus $d_0B=\frac12 I_{2m}$, and
\begin{equation}\label{eqn:contact_sphere_valuations} \phi_{2k }(S^{2m})= 2\tr(\wedge^{2m-2k}I_{2m})=2{2m \choose 2k}\end{equation}}
\end{Example}

\section{Symplectic-invariant distributions}\label{sec:grassmannians}

\subsection{Linear algebra}

The real anti-symmetric matrices of size $2N\times 2N$ will be denoted $\Alt_{2N}$. 
Define $\textrm{SDiag}(\lambda_1,\dots,\lambda_N)\in\Alt_{2N}(\R)$ to be the block-diagonal matrix consisting of the $2\times 2$ blocks 
\[\textrm{SDiag}(\lambda_j)=\left(\begin{array}{cc} 0 & \lambda_j\\
-\lambda_j & 0\end{array}\right) \]
The following is a standard fact from linear algebra.
\begin{Lemma}
	For $A\in \Alt_{2N}$ there is a matrix $B\in\OO_{2N}(\R)$ and $D=\emph{SDiag}(\lambda_1,\dots,\lambda_N)$ such that $A=B^TDB$. The vector $(\lambda_1,\dots,\lambda_N)$ is uniquely defined up to permutations and signs of the $\lambda_j$.
\end{Lemma} 
Denote $\Delta_N=\{\lambda_1\geq \dots\geq \lambda_N\geq 0\}$. For $A\in\Alt_{2N}$, let $\Lambda(A)\in\Delta_N$ be the unique vector such that $A=B^T\cdot\textrm {SDiag}(\Lambda(A))\cdot B$ for some $B\in\OO_{2N}(\R)$.

Recall the multi-K\"ahler angles $0\leq \theta_1\leq \dots\leq \theta_\kappa\leq \frac{\pi}{2}$ of $E\in\Gr^{\R}_{2k}(\C^n)$ introduced by Tasaki \cite{tasaki}, where $\kappa=\min(k, n-k)$. They are defined as follows: choose an orthonormal basis $(e_i)_{i=1}^{2k}$ of $E$, and define the symplectic Gramm matrix $A=\omega(e_i, e_j)$. Then $\Lambda(A)=(\cos\theta_1,\dots,\cos\theta_\kappa)$.
\begin{Proposition}\label{prop:uniform_distribution}
	Let $\theta_i$, $i=1,\dots,\kappa$ be the multi-K\"ahler angles of a subspace $E\in\Gr^{\R}_{2k}(\C^n)$ chosen at random (with respect to the $\SO(2n)$-invariant probability measure). Then the probability distribution of  $(\cos\theta_i)_{i=1}^\kappa$ is uniform in $\Delta_\kappa$.
\end{Proposition}
\proof
We may assume $2k\leq n$ so that $\kappa=k$.
We first observe that the distribution is independent of $n$. Indeed we may condition on the event $E\subset F$ where $F$ is any fixed complex $k$-dimensional subspace in $\C^n$, but the distribution of the multi-K\"ahler angles is clearly independent of $F$. Thus we assume $n=2k$.

Next notice that for a symplectic subspace $E$, there is a unique, up to order, decomposition $E=E_1\oplus\dots\oplus E_k$, where $\dim E_j=2$ and all $F_j:=\mathbb C E_j$ are pairwise orthogonal.
This decomposition can be found as follows: for an orthonormal basis $e=(e_j)_{j=1}^{2k}$ of $E$, consider the matrix $M(E, e)=(\omega(e_i, e_j))_{i,j=1}^{2k}\in \Alt_{2k}$. 
If $(e')$ is a different orthonormal basis, there is an equality of row vectors $(e'_i)=(e_j)B$ for some $B\in\OO(2k)$, and one checks that $M(E, e')=B^TM(E, e)B$. 
Thus for generic $E$ there is a unique, up to order, orthonormal basis $e$ such that $M(E, e)=\textrm{SDiag}(\Lambda(M(E,e)))$. We then set $E_j=\Span(e_{2j-1}, e_{2j})$ for $1\leq j\leq k$.

For a given decomposition $\C^{2k}=F_1\oplus\dots\oplus F_k$ into orthogonal copies of $\C^2$, the multi-K\"ahler angles of $E$ are the collection of K\"ahler angles of $E\cap F_j\in \Gr_2^{\R}(F_j)$. Thus conditioning on the decomposition, we conclude the multi-K\"ahler angles of $E$ are independent and identically distributed, and it remains to find the distribution of the K\"ahler angle of a real 2-plane $E\subset\C^2$. 

We will work with the oriented Grassmannian.
Consider $\R^4=\C^2$ with \[\omega((x_1+ix_3, x_2+ix_4), (x'_1+ix'_3, x_2+ix'_4))=x_1x_3'-x_1'x_3+x_2x_4'-x_2'x_4.\] 
Let $(u,v)$ be an orthonormal basis for $E\in\Gr^+_2(\R^4)$.
We identify $\Gr_2^+(\R^4)=S^2\times S^2$, $E\mapsto (z, w)$ using the standard Euclidean structure and the Pl\"ucker embedding $\Gr^+_2(\R^4)\subset S(\wedge^2 \R^4)$, $E\mapsto u\wedge v=(x_{12},\dots, x_{34})$ where $x_{kl}=u_kv_l-u_lv_k$, followed by the change of coordinates 

\begin{align*}
	x_{12} & :=\frac{w_1+z_1}{2} \quad x_{34}:=\frac{w_1-z_1}{2} \quad x_{13}:=-\frac{w_2+z_2}{2}\\
	x_{24} & :=\frac{w_2-z_2}{2} \quad x_{14}:=\frac{w_3+z_3}{2} \quad x_{23}:=\frac{w_3-z_3}{2}. 
\end{align*}

The corresponding measure on $S^2\times S^2$ is the standard one. Then for $E=(z,w)\in S^2\times S^2$, and $u, v$ an oriented orthonormal basis of $E$, $\omega(u, v)=x_{13}+x_{24}=-z_2$. Denoting by $\theta_E$ the K\"ahler angle, we conclude that $\cos\theta_E=|\Lambda(E)|=|\omega(u,v)|=|z_2|$ is distributed uniformly in $[0,1]$ by the theorem of Archimedes.

\endproof

Let $\Delta^1_\kappa=\{1\geq \lambda_1\geq\dots\geq\lambda_\kappa\geq 0\}$.
Denote by $\Lambda:\Gr_{2k}^\R(\C^n)\to \Delta^1_\kappa$ the vector $(\cos\theta_i)_{i=1}^\kappa$. 
We conclude that $\Lambda_*(d\sigma_1)=\kappa!\prod_{i=1}^\kappa d\lambda_i$.

\subsection{Powers of the Pfaffian}

Consider the meromorphic families of generalized functions $|\Pfaff|_\pm^s\in C^{-\infty}(\Alt_{2N}(\R))$, $s\in\mathbb C$.
They can be constructed by first considering $\Re s>-1$, whence $|\Pfaff|_\pm^s$ is an integrable function, and then using the Cayley-type identity (see e.g. \cite{pfaffian})
\[\Pfaff(\partial) \Pfaff(X)^{s+1}= (s+1)(s+3)\dots(s+2N-1)\Pfaff(X)^s \]
for a meromorphic extension to $s\in\C$. Thus the poles of both families are at $s=-1,-2,\dots$. 

\begin{Theorem}[Muro \cite{muro}]\label{thm:muro}
	The linear combination $|\Pfaff(X)|^s:=|\Pfaff(X)|_+^s+|\Pfaff(X)|_-^s$ is analytic at even $s\in\mathbb Z$ and has a simple pole at odd $s\leq -1$. The linear combination $\sign\Pfaff(X)|\Pfaff(X)|^s:=|\Pfaff(X)|_+^s-|\Pfaff(X)|_-^s$ is analytic at odd $s\in\mathbb Z$ and has a simple pole at even $s\leq -2$.
\end{Theorem}

Let $(V,\omega)$ be a $2n$-dimensional symplectic space, $P$ a compatible Euclidean structure with corresponding complex structure $J$, so that $\omega(u,v)=P(Ju,v)$. Then $P$ induces a Lebesgue measure $\vol_P(E)$ on all subspaces $E\subset V$. Define $\sigma_{\omega,P}:\Gr^+_{2k}(V)\to[-1,1]$ by $\sigma_{\omega, P}(E)=\frac{\omega^k|_E}{\vol_P(E)}$. We will often omit $P$ from the index when no confusion can arise. 

We will now define meromorphic families of distributions on $\Gr^+_{2k}(\R^{2n})$ given by $\mu_\pm(s)=|\sigma_\omega(E)|_\pm^sd\sigma_1(E)$ for large $\Re (s)$. The construction is virtually identical to the one carried out in \cite{faifman_opq}, with $\Sp(2n)$ replacing $\OO(p,q)$. We present it here for the reader's convenience.

We will assume for now that $2k\leq n$. 

Let $U\subset \Gr_{2k}(V)$ be an open set. Let $B_E=(u_1(E),\dots,u_{2k}(E)):U\to V^k$ be a smooth field of $P$-orthonormal bases of $E\in U$. Define the function $M_P:U\to \Alt_{2k}(\mathbb R)$ given by $M_P(E)=\omega(u_i(E),u_j(E))_{i,j=1}^{2k}$. Note that $\sigma_{\omega}(E)=\Pfaff M_P(E)$.

Denote by $U_P\subset\Gr_k(V)$ the open, dense subset of subspaces $E\in \Gr_{2k}(V)$ for which $E\cap JE=\{0\}$. Clearly $E\in U_P$ if and only if $1\not\in \Lambda(M_P(E))$.

\begin{Lemma}\label{lem:matrix_is_submersion}
	$M_P$ is a proper submersion at every $E\in U\cap U_P$.
\end{Lemma}

\proof
Consider a curve $\gamma_1$ through $E$ given by \[\gamma_1(t)=\Span(u_1(t),u_2,\dots,u_{2k})\] with $u_2,\dots, u_{2k}$ fixed, and $\xi=\dot u_1(0)\in E^P$ arbitrary. It follows that 
\[D_E M_P (\dot{\gamma_1}(0))=
\begin{pmatrix}
0 & \omega(\xi, u_2) & \cdots & \omega(\xi,u_{2k}) \\
-\omega(\xi,u_2) & 0 & \cdots & 0 \\
\vdots  & \vdots  & \ddots & \vdots  \\
-\omega(\xi,u_{2k}) & 0 & \cdots & 0
\end{pmatrix}
\]
Since $E\in U_P$, $\omega(u_j, \bullet)_{j=2}^{2k}$ are linearly independent functionals in $\xi\in E^P$, and so the first row of $D_\Lambda M_P(\dot{\gamma_1})$ is arbitrary, while the other entries in the upper triangle vanish. Replacing  $\gamma_1$ with $\gamma_j$ in the obvious way, we conclude $D_E M_P(\sum \alpha_j \dot\gamma_j(0))$ can be arbitrary, thus concluding the proof.
\endproof

\begin{Lemma}\label{lem:compatible_cover}
	One can choose finitely many $\omega$-compatible Euclidean structures $P_i$ s.t. $\{U_{P_i}\}$ cover $\Gr_{2k}(V)$.
\end{Lemma}
\proof
Given $E\in \Gr_{2k}(V)$, a generic choice of an $\omega$-compatible $(P, J)$ would have $E\cap JE=\{0\}$ by a trivial dimension count. Fixing one such $J=J(E)$ with corresponding $P(E)$ for every $E$, we get an open cover $(U_{P(E)})_{E\in\Gr_{2k}(V)}$ of $\Gr_{2k}(V)$. The claim now  follows by the compactness of $\Gr_{2k}(V)$.
\endproof

We now explain how to pull-back $|\Pfaff(X)|^s_\pm$ to $\Gr_{2k}^+(V)$, using the locally-defined submersion $M_P$.
\begin{Definition}
	For $s\in\mathbb C$, let $\mathcal D^s$ be the line bundle of $s$-densities over $\Gr_{2k}^+(V)$, which has fiber $\Dens^s(E)$ over $E\in \Gr_{2k}^+(V)$. 
	We say that a choice of generalized section $f(s)\in\Gamma^{-\infty}(U, \mathcal D^s)$ over $U\subset \Gr_{2k}^+(V)$ for $s\in\Omega\subset\C$ is meromorphic in $s$ if, having fixed a Euclidean metric $P$ and using it to identify all bundles $\mathcal D^s$, one obtains a map $f_P:\Omega\to  C^{-\infty}(U)$ which is meromorphic in $s$. 
	
	We denote by $\mathfrak{M}^{-\infty}(\mathcal D^s)$ the sheaf for which $\Gamma(U, \mathfrak{M}^{-\infty}(\mathcal D^s))$ is the space of meromorphic in $s$ maps $\mathbb C\to \Gamma^{-\infty}(U, \mathcal D^s)$. 
	
\end{Definition}

Recall the orbits $X^{2k}_\pm$,  $(X^{2k,+}_r)_{r=0}^{\kappa-1}$  of $\Gr^+_{2k}(V)$ under $\Sp(V)$ defined in the paragraph following eq. \eqref{eq:orbits}, where $\kappa=\min(k,n-k)$. In terms of multi-K\"ahler angles, $E\in X^{2k,+}_r$ precisely when exactly $r$ of the angles are distinct from $\frac \pi 2$.

We are now ready to construct the meromorphic families.
\begin{Proposition}\label{prop:muro_pullback}
	There are global sections $f_\pm(s)=M_P^*|\Pfaff|^s_\pm$ of $\mathfrak{M}^{-\infty}(\mathcal D^s)$ supported on $\overline {X^{2k}_\pm}$, s.t. whenever $s$ is not a pole of $f_\pm$,  $f_\pm(s)$ is $\Sp(V)$-invariant.
\end{Proposition}
\proof
Assume first $2k\leq n$. Let $P_i$ be a finite collection of $\omega$-compatible Euclidean structures as in Lemma \ref{lem:compatible_cover}, and let $U_i=U_{P_i}\subset\Gr^+_{2k}(V)$ be the corresponding open sets of generic subspaces. For each $i$, cover $U_i$ by open sets $U_{ij}\subset U_i$ so that $M_{ij}=M_{P_i}:U_{ij}\to \Alt_{2k}$ can be defined by some smooth field of orthonormal bases of $E$ over $ U_{ij}$. Now since $M_{ij}$ is a proper submersion, one obtains a meromorphic in $s$ family of functions $\tilde f^\pm_{ij}(E;  s)\in C^{-\infty} (U_{ij})$ given by $\tilde f_{ij}^\pm(\bullet; s)=M_{ij}^*|\Pfaff|_\pm^s$.

It then obviously holds that on $U_{ij}\cap U_{ij'}$, $\tilde f^\pm_{ij}(\Lambda; s)$ and $\tilde f_{ij'}^\pm(\Lambda; s)$ coincide as continuous functions for $\Re(s)>0$. Therefore, they coincide on $U_{ij}\cap U_{ij'}$ as meromorphic functions, and we may merge all $\tilde f_{ij}^\pm$ into one meromorphic family $\tilde f_i^\pm(\bullet; s)\in C^{-\infty} (U_i)$. The corresponding (through $P_i$) section $f_i^\pm\in \Gamma(U_i, \mathfrak{M}^{-\infty}(D^s))$ is obviously $\Sp(V)\cap \OO(P_i)$-invariant. Moreover, it is $\mathfrak{sp}(V)$-invariant.

Next, we claim that $f_i^\pm$ and $f_{i'}^\pm$ coincide on $U_i\cap U_{i'}$. Since both are meromorphic, we may assume in the following that $\Re(s)>0$.

It is easy to see, using Lemma \ref{lem:LocallyClosedOrbits} as in the proof of Proposition \ref{prop:klain_upper_bounds}, that for $\Re(s)>0$, no $\Sp(V)$-invariant generalized sections of $\mathcal D^s$ can be supported on a set of positive codimension: 
using Lemma \ref{lem:LocallyClosedOrbits}, we consider for $\alpha\geq0$ the bundle $F^\alpha_{X^{2k, +}_{r}}$ over $X^{2k,+}_{r}$ where $r>0$, that has fiber $\Dens(E)^s\otimes \Dens^*(N_EX^{2k,+}_{r})\otimes \Sym^\alpha (N_EX^{2k,+}_{r})$ over $E$.
Denoting $E_0=E\cap E^\omega$, by Corollary \ref{cor:normal_space}
\[F^\alpha_{X^{2k,+}_{r}}|_E=\Dens(E)^s\otimes \Dens(\wedge^2E_0)\otimes \Sym^\alpha (\wedge^2E_0^*)\]
which clearly contains no $\Sp(V)$-invariants.  It follows that the space of $\Sp(V)$-invariants in $\Gamma^{-\infty}(\Gr^+_{2k}(V), \mathcal D^s)$ supported on ${\overline{X^{2k}_{\epsilon}}}$ is at most 1-dimensional for each $\epsilon\in\{\pm\}$. 

Since $U_i\subset \Gr^+_{2k}(V)$ is dense, it follows by construction that for $\Re(s)>0$, $f^\pm_i(\bullet; s)$ extends by continuity to an $\Sp(V)$-invariant section of $\mathcal D^s$ over $\Gr^+_{2k}(V)$ supported on $\overline{X^{2k}_{\pm}}$, and by the previous paragraph we can find meromorphic functions $c_i(s)$, such that $c_1(s)=1$ and the sections $c_i(s) f_i^\pm(\Lambda;  s)$ coincide for all $i$. Denoting by $p_i(s)\in\Gamma^\infty (\Gr^+_{2k}(V), \mathcal D^s)$ the Euclidean section defined by $P_i$,  it holds that \[f^\pm_i(E;s)= |\Pfaff M_{P_i}(E)|^{s}p_i(s)\] for $E\in X^{2k}_{\pm}$, so that \[|\Pfaff M_{P_1}(E)|^{s}p_1(s)=c_i(s) |\Pfaff M_{P_i}(E)|^{s}p_i(s)\]
implying 
\[c_i(s)=\left(\frac{|\Pfaff M_{P_1}(E)|}{|\Pfaff M_{P_i}(E)|}\frac{p_1(1)}{p_i(1)}\right)^s\] 
for all $E\in X^{2k}_{\pm}$. Since $c_i(s)$ is independent of $E$, one has $c_i(s)=c_i^s$ for some $c_i>0$. 

Finally, for $s=1$, $\mathcal D^1$ is the Klain bundle of Lebesgue measures, and it is easy to see that all $f_i^\pm(1)$ represent the Liouville measure induced by the symplectic form on every symplectic subspace $E$.  It follows that $c_i=1\Rightarrow c_i(s)\equiv 1$. Thus we have shown that $f_i^\pm$, $f_{i'}^\pm$ coincide in $\Gamma(U_i\cap U_{i'}, \mathfrak M^{-\infty}(\mathcal D^s))$.

We conclude there is a globally defined section $f_\pm$ of $\mathfrak{M}^{-\infty}(\mathcal D^s)$ which is $\mathfrak{sp}(V)$-invariant and supported on $\overline {X^{2k}_\pm}$, respectively. For $\Sp(V)$-invariance, we observe it holds for $\Re(s)>0$ and then invoke uniqueness of meromorphic continuation. 

This concludes the proof when $2k\leq n$. For the case $2k>n$, we simply use the oriented skew-orthogonal complement map $\Gr^+_{2k}(V)\to \Gr^+_{2n-2k}(V)$ to pull-back $f_s$. This is a valid operation since we have the equivariant identification $E^\omega\simeq (V/E)^*\Rightarrow \Dens(E^\omega)\simeq\Dens(E)$.

\endproof

\begin{Definition}
	Set $|\sigma_\omega|_\pm^s\in C^{-\infty}(\Gr_{2k}^+(V))$ to be the value of $f_\pm(s)$ under the Euclidean trivialization.
\end{Definition}

\begin{Lemma}\label{lem:total_measure}Write $\kappa=\min(k,n-k)$.
One has $\int_{\Gr_{2k}(\R^{2n})}|\sigma_{\omega}|_\pm^s d\sigma_1(E)=\frac{1}{2\kappa!(s+1)^\kappa}$.
\end{Lemma}
\proof
Using the cosines of the multi-K\"ahler angles in decreasing order, denoted $\Lambda(E)=(\lambda_1,\dots, \lambda_\kappa)$, $\lambda_j=\cos\theta_j$, we have $|\sigma_\omega(E)|=\prod_{j=1}^\kappa \lambda_j$. Then for $\Re(s)>0$, using Proposition \ref{prop:uniform_distribution} we get
\[ \int_{\Gr_{2k}(\R^{2n})}|\sigma_{\omega}|^s d\sigma_1(E)= \frac{1}{\kappa!}\int_{[0,1]^\kappa}\prod_{j=1}^\kappa \lambda_j^sd\lambda_j  = \frac{1}{\kappa!(s+1)^\kappa} \]
and the result follows by uniqueness of meromorphic extension.
\endproof

The two values $s=-2n,-(2n+1)$ are of particular interest, as evidenced in the two theorems below. The first theorem concerns the linear Grassmannian:

\begin{Theorem}\label{thm:even_measure}
	The distribution $\mu^+_\omega\in\mathcal M^{-\infty}(\Gr^+_{2k}(V))$ given by $\mu^+_\omega:=|\sigma_\omega(E)|^{-2n}d\sigma_1(E)$ is $\Sp(V)$-invariant, has full support and $\int_{\Gr_{2k}(V)}\mu_\omega^+\neq 0$. It is even with respect to orientation reversal.
\end{Theorem}
In particular, we get a canonically normalized $\Sp(V)$-invariant distribution $\mu_\omega:=\pi_*\mu_\omega^+$ on $\Gr_{2k}(2n)$, where $\pi:\Gr^+_{2k}(V)\to\Gr_{2k}(V)$ is the double cover map.
\proof
A distribution over $\Gr^+_{2k}(V)$ is a generalized section of the bundle with fiber $\Dens(T_E\Gr_{2k}^+(V))$ over $E$, which is $\Sp(V)$-isomorphic to $\mathcal D^{-2n}$. All statements follow immediately from Proposition \ref{prop:muro_pullback}, Theorem \ref{thm:muro} and Lemma \ref{lem:total_measure}.
\endproof
We conjecture that $\mu_\omega$ is the unique $\Sp(V)$-invariant distribution on $\Gr_{2k}(V)$. This was shown by Gourevitch, Sahi and Sayag in \cite{GSS} for $k=n$ when $n$ is even.

Similarly, we have a statement for the affine Grassmannian. We define $|\sigma_\omega(E)|^s\in C^{-\infty}(\AGr^+_{2k}(V))^{tr}$ by pulling-back by the projection map $\AGr^+_{2k}(V)\to\Gr^+_{2k}(V)$. Let $dE$ be the measure on $\AGr_{2k}^+(V)$ which is built out of $d\sigma_1$ on $\Gr_{2k}^+(V)$, and the Euclidean measure on translations. Define the odd distribution \begin{equation}
\label{eq:affine_distribution}
\overline \mu_\omega:=\sign(\sigma_\omega)|\sigma_\omega(E)|^{-2n-1}dE\in\mathcal M^{-\infty}(\AGr^+_{2k}(V))\end{equation} and the even distribution  \begin{equation}
\label{eq:affine_distribution_singular}\overline\mu_0:=\Res_{s=-2n-1}|\sigma_\omega(E)|_+^{s}dE\in\mathcal M^{-\infty}(\AGr^+_{2k}(V)).\end{equation}

\begin{Theorem}
The distributions $\overline \mu_\omega$, $\overline \mu_0$ are $\overline{\Sp(V)}$-invariant. $\overline \mu_\omega$ has full support. 
\end{Theorem}
The proof is as in the linear case.
In particular, there is a an even, canonically normalized $\overline{\Sp(V)}$-invariant distribution $\overline \mu_0:=\pi_*\mu_0^+$ on $\AGr_{2k}(V)$, supported on the $\omega$-degenerate subspaces.

We will need the following a-priori information about the wavefront set of $\overline \mu_\omega$.

\begin{Proposition}
	The wavefront set of $\overline\mu_\omega$ belongs to $\cup_rN^*X^{2k}_r$.
\end{Proposition}
\proof
This is immediate from $\overline{\Sp(V)}$-invariance.
\endproof

\section{The contact sphere}\label{sec:contact_sphere}

\subsection{A Crofton basis for $S^{2n+1}$}

In this section, $V=\R^{2n+2}$, $M=S^{2n+1}=\mathbb P_+(V)$. Take $\xi \in \mathbb P_+(V)$ so that $T_\xi S^{2n+1}=\xi ^*\otimes V/\xi$. The contact hyperplane is $H_\xi=\xi^*\otimes \xi^\omega/\xi\subset T_\xi  \mathbb P_+(V)$. Then $T_\xi  S^{2n+1}/H_\xi=\xi^*\otimes V/\xi^\omega\simeq(\xi^*)^{\otimes 2}$ $\Stab(\xi)$-equivariantly. The symplectic form on $\xi^\omega/\xi$ defines a form 
$\omega_\xi\in \wedge^2 H_\xi^*\otimes (\xi^*)^{\otimes 2}=\wedge^2 H_\xi^*\otimes T_\xi  S^{2n+1}/H_\xi$.
We remark that the form $\omega_\xi$ is determined by the contact distribution alone, without reference to the form $\omega$ on $V$, see Example \ref{exm:contact_is_DH}.
The linear symplectic group $\Sp(V)$ acts on $S^{2n+1}$ by contactomorphisms. The stabilizer of $\xi\in S^{2n+1}$ is easily seen to act on $T_\xi S^{2n+1}$, which is a dual Heisenberg algebra, by its full group $\Sp_H(2n+1)$ of automorphism.

Now $\dim \Sp(2n+2)-\dim \Sp_H(2n+1)=(n+1)(2n+3)-(n+1)(2n+1)=2n+2=\dim S^{2n+1}+1$, so that $\dim \Stab(\xi)=\dim  \Sp_H(2n+1)+1$, and we denote $\Stab(\xi)=\Sp^*_H(2n+1)$. The kernel of the restriction homomorphism  $\Sp^*_H(2n+1)\to \Sp_H(2n+1)$ consists of the linear maps that fix all $v\in \xi^\omega$ and acts on some $w\in V\setminus \xi^\omega$ by $w\mapsto w+\lambda\xi$, for some $\lambda\in \R$.

Thus $S^{2n+1}=\Sp(2n+2)/\Sp_H^*(2n+1)$. Inspired by the analogy to the Riemannian symmetric space presentation $S^{2n+1}=\SO(2n+2)/\SO(2n+1)$, we look for Crofton-type formulas for the contact valuations on the sphere.

Consider the double fibration 
\[\xymatrix{& W\ar[dl]_{\tau} \ar[dr]^{\pi}&\\ \Gr_{2k}(V) & & S^{2n+1} }\] 
where $W$ is the partial flag manifold $\{(E, \theta)\in  \Gr_{2k}(V)\times S^{2n+1}: \theta\in E \}$

\begin{Definition}
	
	For $0\leq k\leq n$, define the generalized valuation $\psi_{2k}$ through the Crofton formula (in the sense of subsection \ref{sec:crofton}) \[\psi_{2k}:=\pi_*\tau^*(\mu_\omega)=\int_{\Gr_{2n+2-2k}(V)}\chi(\bullet\cap E)\sigma_\omega(E)^{-2n-2}d\sigma_1(E)\in \mathcal V^{-\infty}(S^{2n+1})^{\Sp(V)}\] 

\end{Definition}

It follows from Lemma \ref{lem:crofton_filtration_level} that $\psi_{2k}\in \mathcal W^{-\infty}_{2k}(S^{2n+1})$. Moreover, since $\psi_{2k}(S^{2k})\neq 0$ by Theorem \ref{thm:even_measure} (see the computation preceding eq. \eqref{eq:sphere_integral_values} below for precise value), we conclude $\psi_{2k}\in \mathcal W^{-\infty}_{2k}(S^{2n+1})\setminus \mathcal W^{-\infty}_{2k+1}(S^{2n+1})$, and in particular all $\psi_{2k}$ are linearly independent.

\begin{Proposition}
	$\mathcal V^{-\infty}(S^{2n+1})^{\Sp(V)}$ is spanned by $\chi,\psi_2,\dots,\psi_{2n}$.
\end{Proposition}
\proof
The proof is identical to that of Theorem \ref{thm:contact_valuations_classification}, with $\Sp(V)$ replacing the full group of contactomorpisms. The proof remains valid since $\Sp(V)$ acts transitively, with the same action of the stabilizer on the tangent space by $\Sp_H(2n+1)$.
\endproof
This completes the proof of Theorem \ref{mainthm:sphere_Hadwiger}.
In light of Theorem \ref{thm:contact_valuations_classification} we get
 \begin{Corollary}\label{cor:abstract_crofton_sphere}
 	For $0\leq k\leq n$, $\psi_{2k}$ are linear combinations of $\phi_{2j}$, $0\leq j\leq n$. In particular, $\psi_{2k}$ is invariant under all contactomorphisms of $\mathbb P_+(V)$. 
 \end{Corollary}
Thus we establish Theorem \ref{mainthm:crofton_sphere}, except for the explicit determination of the coefficients which is deferred to the next subsection.

It follows also that $\psi_{2k}(F)$ is well-defined for subsets $F\subset S^{2n+1}$ normally transversal to the contact distribution. We will make use of the following lemma.

\begin{Lemma}\label{lem:Crofton_wavefront}
	Assume that a closed submanifold $F\subset S^{2n+1}$ is normally transversal, and $\chi_F:=\chi(F\cap\bullet)\in C^{-\infty}(\Gr_{2n+2-k}(V))$ has wavefront disjoint from $\WF(\mu_\omega)$. Then
	\[ \psi_{2k}(F)=\int_{\Gr_{2n+2-2k}(V)} \chi(F\cap E)d\mu_\omega(E)  \]
\end{Lemma}
\proof
Choose an approximate identity $\rho_\epsilon$ on $\GL(V)$, and define $\phi_\epsilon= \int_{\GL(V)}g^*\chi_F \cdot\rho_\epsilon(g)dg\in\mathcal V^\infty(S^{2n+1})$. Then 
\[ \psi_{2k}(F)=\lim_{\epsilon\to 0}\langle \psi_{2k}, \phi_{\epsilon}\rangle= \lim_{\epsilon\to 0} \int_{\Gr_{2n+2-2k}(V)} \phi_\epsilon(E)  d\mu_\omega(E)\]
Now $\phi_{\epsilon}(E)=\int_{\GL(V)}\chi_F(gE)\rho_\epsilon(g)dg$ converges to $\chi_F\in C^{-\infty}(\Gr_{2n+2-2k}(V))$  in H\"ormander's topology on the space of generalized functions with wavefront set contained in $\WF(\chi_F)$. It follows that $\psi_{2k}(F)= \int_{\Gr_{2n+2-2k}(V)} \chi_F(E)d\mu_\omega(E)$.

\endproof

 \subsection{Integral geometry of the contact sphere}
 Here we determine the coefficients in Theorem \ref{mainthm:crofton_sphere}. We have two different bases of contact-invariant valuations on the contact sphere, indexed by $0\leq k\leq n$. Namely, we have $\phi_{2k}$ defined in terms of curvature at the contact points, and $\psi_{2k}$ given by Crofton integrals. Since  $\phi_{2k}, \psi_{2k}\in \mathcal W_{2k}^{-\infty}(S^{2n+1})$, they are related by a triangular matrix, that is \begin{equation}\label{eq:crofton_sphere}\psi_{2k}=\sum_{j=k}^n c^n_{kj}\phi_{2j}.\end{equation}
 We will compute $c_{kj}^n$ by evaluating both bases on all great spheres $S^{2m}$.
 
 Take $F=S^{2n}=S^{2n+1}\cap \Pi$ where $\Pi\subset S^{2n+1}$ is a fixed hyperplane. It is normally transversal. Note also that $\chi(F\cap E)=2$ for a generic $E\in \Gr_{2n+2-2k}(\R^{2n+2})$, that is $\chi_F$ is the constant $2$ on $\Gr_{2n+2-2k}(\R^{2n+2})$. Thus by Lemma \ref{lem:Crofton_wavefront} we may compute $\psi_{2k}(F)$ using the explicit Crofton formula. Denote $\kappa_k=\min (k, n+1-k)$. By Lemma \ref{lem:total_measure}, \begin{align*}\psi_{2k}(S^{2n})&=\int_{\Gr_{2n+2-2k}(\R^{2n+2})}\chi(S^{2n}\cap E)\sigma_\omega(E)^{-2n-2}d\sigma_1(E)\\&=2\int_{\Gr_{2n+2-2k}(\R^{2n+2})}\sigma_\omega(E)^{-2n-2}d\sigma_1(E)=(-1)^{\kappa_k}\frac{2}{\kappa_k!}\frac{1}{(2n+1)^{\kappa_k} }\end{align*}
 
 Considering spheres of lower dimension, we see that $\psi_{2k}(S^{2m})=\psi_{2k}(S^{2n})$ if $m\geq k$ and zero otherwise. Thus
 
 \begin{equation}\label{eq:sphere_integral_values}
 \psi_{2k}(S^{2m})=\left\{\begin{array}{cc}(-1)^{\kappa_k}\frac{2}{\kappa_k!}\frac{1}{(2n+1)^{\kappa_k}},& k\leq m\\ 0,& k>m\end{array}\right.
 \end{equation}
 On the other hand by eq. \eqref{eqn:contact_sphere_valuations}, 
 
 \[  \phi_{2j}(S^{2m})=2{2m \choose 2j} \]
 
 We now plug those values into \eqref{eq:crofton_sphere}.
 Define the lower-triangular matrix $A$ by $A(m, j)={2m\choose 2j}$, $0\leq m,j\leq n$. Its inverse is given by $A^{-1}(j, m)={2j\choose 2m}E_{2j-2m}$, where $E_i$ is the $i$-th Euler (secant) number.
 Set $b_k=(-1)^{\kappa_k}\frac{1}{\kappa_k!}\frac{1}{(2n+1)^{\kappa_k}}$. Then $c_{kj}^n=b_k\sum_{m=k}^jA^{-1}(j,m)$. In particular, $c^n_{kk}=b_k$, $c^n_{k, k+1}=b_k(1-{2k+2\choose 2})$.
 
\subsection{Contact curvature of convex sets}

Here we prove an upper bound on the contact valuations of a convex set.
\begin{Theorem}
	Let $K\subset  S^{2n+1}$ be a convex subset with $C^2$ boundary. Then for all $0\leq k\leq n$, $\phi_{2k}(\partial K)\leq \phi_{2k}(\partial K_0)$, where $K_0$ is the hemisphere.  
\end{Theorem}
\proof
We will use an auxiliary complex structure: $V=\C^{n+1}$, $S^{2n+1}$ is identified with the unit sphere therein.
We may assume that $\partial K$ is normally transversal. 

First, let us verify that $\partial K$ has exactly two contact tangent hyperplanes. Since the Euler characteristic $\chi(\partial K)=2$, this amounts to verifying that the intersection index of $N^*K$ and $M_H$ in $\mathbb P_{S^{2n+1}}$ at every contact point is $+1$. We will refer to this number as the contact index, denoted $I_H(\partial K, x)$.

Consider a sphere $\mathcal E=S^{2n}\subset S^{2n+1}$ given by a quadratic equation in $V$. It is easy to check by an explicit computation that $\mathcal E$ has exactly two contact tangent points. 
Since $\chi(\mathcal E)=2$, we conclude that the contact index of $\mathcal E$ at each of those points is $+1$.
Now let $x\in \mathcal K$ be a contact tangent point, and let $\mathcal E$ be the osculating sphere at $x$. Then $I_H(\partial K, x)=I_H(\mathcal E, x)=1$.

Next, let $p\in \partial K$ be a point where $T_p\partial K=H_p$. We now project the hemisphere $U$ centered at $p$ to $\R^{2n+1}=T_pS^{2n+1}$ by a central projection $\pi$ from the origin $\pi$, so that $p$ is mapped to the origin. Clearly $\pi(K\cap U)$ is convex near the origin: if $K=S^{2n+1}\cap C$ where $C$ is a convex cone, then $\pi(K\cap U)=C\cap T_pS^{2n+1}$. By a standard computation and assuming $p=\{y_{n+1}=1\}$, the resulting contact structure in $\R^{2n+1}$ is given by the contact form $\alpha=-dx_{n+1}+\sum_{j=1}^n(-y_jdx_j+x_jdy_j)$. We will write $z=-x_{n+1}$.

We thus consider a convex body $K$ with $C^2$ boundary in the contact space $(\R^{2n+1},\alpha)$. We assume $K$ is tangent to the contact distribution at the origin, which is $\R^{2n}=\{z=0\}$, and further assume without loss of generality that $K$ lies below it.
The normal to the contact distribution is $\nu_H=(-y_1,\dots, -y_n, x_1, \dots, x_n, 1)$, the normal to $\partial K$ is $\nu_K$. Then by Proposition \ref{prop:DH_general_formula}, \[\phi_{2n}(\partial K, 0)=|\det(d_0\nu_K-d_0\nu_H)|^{-1}D(d_0\nu_K-d_0\nu_H[2n-2k], J[2k])\] where $d_0\nu_K, d_0\nu_H: \R^{2n}\to T_{e_{2n+1}}S^{2n}$, the latter space identified with $\R^{2n}$. 
Using the coordinates $x_1,\dots, x_n, y_1,\dots, y_n$ on $\R^{2n}$ we get \[d_0\nu_H=\left( \begin{array}{cc}
0 & -I_n\\ I_n & 0
\end{array}  \right)=J\]  Write also $S=d_0\nu_K$. Thus $\phi_{2n}(\partial K, 0)=|\det (S-J)|^{-1}D(S-J[2n-2k], J[2k])=|\det(I+SJ)|^{-1}D(I+SJ[2n-2k], I[2k])=|\det(I+SJ)|^{-1}\tr \wedge^{2n-2k}(I+SJ)$. 

Note that $S\geq 0$ since $K$ is convex. Then $JS$ and $\sqrt SJ\sqrt S$ have the same characteristic polynomial. The latter matrix is antisymmetric, hence the roots of its characteristic polynomial appear in purely imaginary pairs $\pm i\lambda_j$, $j=1,\dots, n$. Let us write $\mu_1,\dots, \mu_{2n}$ for these eigenvalues in some order. Note that for $K_0$, all $\lambda_j=0$. Writing $m=2n-2k$, we ought to show that
\[ \tr \wedge^{m} (I+SJ)\leq {2n \choose m} \det(I+SJ)\]
\[\iff \sum_{|T|=m}\prod_{t\in T} (1+\mu_t) \leq  {2n \choose m}\prod_{j=1}^n(1+\lambda_j^2) \]
We will use induction on $n$. For $n=0,1$ the verification is trivial.
Partition the sum over subsets $T$ as $S_0+S_1+S_2$, where $S_j$ ($j=0,1,2$) is composed of those summands where $\pm i\lambda_1$ appears $j$ times inside $\{\mu_t:t\in T\}$.
By the induction assumption, $S_0\leq {2n-2\choose m}\prod_{j=2}^n(1+\lambda_j^2)$, $S_1\leq 2 {2n-2\choose m-1}\prod_{j=2}^n(1+\lambda_j^2)$, $S_2\leq {2n-2\choose m-2}\prod_{j=1}^n(1+\lambda_j^2)$. It remains to check that
\[ {2n-2\choose m}+2{2n-2\choose m-1}+{2n-2\choose m-2}(1+\lambda_1^2)\leq {2n\choose m}(1+\lambda_1^2) \] 
which clearly follows from the equality
\[ {2n-2\choose m}+2{2n-2\choose m-1}+{2n-2\choose m-2}={2n\choose m}\]
concluding the induction and the proof.
\endproof

\begin{Remark}
	The case of equality is far from unique: any convex subset which is flat to second order at its two contact points would have the same values of $\phi_{2k}$.
\end{Remark}

\section{Symplectic integral geometry}\label{sec:symplectic}

\subsection{Symplectic space}

Let us first show there is no interesting symplectic valuation theory.
\begin{Theorem}
	There is no $\Sp(2n)$-invariant, translation invariant generalized valuation except for linear combinations of $\chi$ and $\vol_{2n}$.
\end{Theorem}

\proof
We will find all invariant $k$-homogeneous valuations. Let $N=\lfloor\frac{\min(k,2n-k)}{2}\rfloor$ be the number of multi-K\"ahler angles for $E\in\Gr_k(\R^{2n})$. Since $\textrm U(n)\subset \Sp(2n)$, an $\Sp(2n)$-invariant valuations $\phi$ would also be $\textrm U(n)$-invariant. In particular, it is smooth by Alesker's theorem \cite{alesker_multiplicative}, as $\textrm U(n)$ acts transitively on the unit sphere.
We consider two separate cases. When $k=2l+1$ and $E\in\Gr_k(V)$ is maximally non-degenerate, then there is no $\Stab(E)$-invariant Lebesgue measure on $E$. Since such subspaces are dense in $\Gr_k(V)$ we conclude $\Kl(\phi)=0$, and hence $\phi=0$.

When $k=2l$, an $\Sp(2n)$-invariant section of the Klain bundle should be proportional to $|\omega^l |_E|$, that is after Euclidean trivialization it is proportional to $|\sigma(E)|=\prod_{i=1}^{N}\cos\theta_i$. But by \cite{bernig_fu_hermitian}, the Klain section of a $\textrm U(n)$-invariant valuations should be given by a symmetric polynomial of $\cos^2\theta_i$.
Thus there can be no non-trivial $\Sp(2n)$-invariant valuations unless $k=0, 2n$.

\endproof
\begin{Remark}
	Instead of using the description of Bernig-Fu, one can simply notice that $\prod_{i=1}^{N}\cos\theta_i\in C(\Gr_k(V))$ is not smooth, violating the smoothness of $\phi$.
\end{Remark}
Recall the distributions $\overline \mu_\omega, \overline \mu_0\in\mathcal M^{-\infty}(\AGr_{2n-2k}^+(V))^{\overline {\Sp(V)}}$ given by eqs. \eqref{eq:affine_distribution} and \eqref{eq:affine_distribution_singular}. 

\begin{Corollary}
	For $1\leq k\leq n-1$ and a smooth convex body $K\subset V$ it holds that \[\int_{\AGr^+_{2k}(V)}\chi(K\cap E) d\overline \mu_0(E)=0.\] 
\end{Corollary}
Put differently, $\overline\mu_0$ lies in the kernel of the cosine transform.
\\\\
Nevertheless, we can write symplectic Crofton formulas with the oriented valuation theory approach detailed in Appendix \ref{app:oriented_valuations}.
For a compact oriented $C^1$-submanifold with boundary $F\subset V$ of codimension $2k$, set $\ind_F(E)=I(E, F)$ for $E\in\AGr^+_{2k}(V)$, which is well-defined whenever $E$ and $F$ intersect transversally.

\begin{Lemma}
	For $F$ as above, $\ind_F\in C^{-\infty}(\AGr^+_{2k}(V))$, and $\WF(\ind_F)\cap N^*X^{2k, +}_{r}=\emptyset$ for all $r<\min(k,n-k)$.
\end{Lemma}

\proof

Let $Z=\{(x,E):x\in E\}\subset V\times\AGr^+_{2k}(V)$ be the incidence manifold, which has a natural orientation, and denote by $\tau:Z\to V$, $\pi:Z\to \AGr^+_{2k}(V)$ the obvious submersions. Consider $Z_F=\tau^{-1}(F)$, which is a $C^1$ oriented submanifold of $Z$ of codimension $2k$. Define $\delta_{Z_F}:=\tau^*[[F]]=[[Z_F]]$.  Note that $\pi_*\tau^*[[F]]\in C^{-\infty}(\AGr^+_{2k}(V))$, and $\pi_*\tau^*[[F]](E)=I(E, F)$ whenever $E\pitchfork F$. The first statement follows. Note that $\WF(\delta_{Z_F})\subset\Im(\tau^*)$, and therefore also $\pi^*\WF(\ind_F)\subset \Im(\tau^*)$.
As $\pi^*$ is injective, the statement of the lemma would follow from $ \Im(\tau^*)\cap \pi^*N^*X^{2k, +}_{r}=\{0\}\iff \Im(\tau^*)\cap N^*\pi^{-1}X^{2k, +}_{r}=\{0\}$. 

Take $(x,E)\in Z_F\cap \pi^{-1} X^{2k, +}_{r}$. Define $\tilde Z_x:=\tau^{-1}(x)\subset Z$. Then $\Im(\tau^*)\cap T^*_{x,E}=\tau^* (T_x^*V)=N^*\tilde Z_x$. Let us check that $\tilde Z_x$ and $\pi^{-1}X^{2k, +}_{r}$ intersect transversally at $(x, E)$.
For a tangent vector $(v,\Xi)\in T_{x,E}Z\subset T_x V\times T_E\AGr_{2k}^+(V)$, use Lemma \ref{lem:transversal_intersection} to decompose $\Xi=\Xi_x+\Xi_L$ with $\Xi_x\in T_EZ_x$ where $Z_x:=\pi\tilde Z_x$, and $\Xi_L\in T_EX^{2k, +}_{r}$.
In the product manifold $V\times\AGr^+_{2k}(V)$ we get the equality $(v, \Xi)=(0,\Xi_x)+(v, \Xi_L)$. 
Since $(0, \Xi_x)\in TZ$ while by assumption $(v, \Xi)\in TZ$, we conclude $(v, \Xi_L)\in TZ$. Thus $T_{x, E}Z=T_{x, E}\tilde Z_x+T_{x,E}\pi^{-1}X^{2k, +}_{r}\Rightarrow N^*Z_F\cap N^*\pi^{-1}X^{2k, +}_{r}=\{0\}$, concluding the proof.

\endproof

\begin{Theorem} Let $F\subset V$ be a $C^1$ compact, oriented submanifold with boundary. Then
	\[  \int_F \omega^k=\frac{(-1)^\kappa}{2}{n\choose k}{2n \choose 2k}^{-1}n^\kappa \int_{\AGr^+_{2n-2k}(V)} I(E, F)d\overline\mu_\omega(E) \]
	where $\kappa=\min(k, n-k)$.
	
\end{Theorem}
\proof
Considered as a function of $F$, the integral on the right hand side is a Crofton integral as in Proposition \ref{prop:oriented_crofton}. Hence it defines a closed, $2k$-form on $V$ which is $\overline{\Sp(V)}$-invariant. By the fundamental theorem of invariant theory, it is a multiple of $\omega^k$, that is 
\begin{equation}\label{eqn:symplectic_crofton}\int_{\AGr^+_{2n-2k}(V)} I(E, F)d\overline\mu_\omega(E) =C\int_F \omega^k\end{equation}

It remains to find the constant $C$. We will use a compatible Euclidean structure. Let $B_W$ be the unit Euclidean ball in the $\omega$-positively oriented subspace $W\in X^{2k}_+$. We will average the integral over $X^{2k}_+$ with respect to the probability measure $dW$ that is invariant under $\mathfrak {so}(2n)$. For easy computation, we replace the exponent of $\sigma_{\omega}$ in $\overline \mu_{\omega}(E)=\sign \sigma_\omega(E)|\sigma_\omega(E)|^{-2n-1}d E$ with the meromorphic variable $s\in\mathbb C$ and compute for real $s$ which is sufficiently large so that all integrands are continuous.

\[A_s:=\int_{  X^{2k}_+}dW\int_{\AGr^+_{2n-2k}(V)} I(E, B_W)\sign \sigma_\omega(E)|\sigma_\omega(E)|^{s}dE
\]
\[=\int_{\Gr^+_{2n-2k}(V)}\sign \sigma_\omega(E)|\sigma_\omega(E)|^{s}d\sigma_1(E)\int_ { X^{2k}_+}dW\int_{V/E}I(E+x, B_W)dx   \]
Here the inner integral is with respect to $dx$, the Euclidean Lebesgue measure under the identification $V/E=E^\perp$. We may write 
\[\int_{V/E}I(E+x, B_W)dx =\sign\sigma_\omega(E)\vol_{2k}(\textrm{Pr}_{E^\perp}(B_W))\]
Hence 
\[A_s= \int_{\Gr^+_{2n-2k}(V)}|\sigma_\omega(E)|^{s}d\sigma_1(E)\int_ { X^{2k}_+} \vol_{2k}(\textrm{Pr}_{E^\perp}(B_W))dW  \]
The inner integral is independent of $E$ and can be computed using the Kubota formula. 
\[\int _{\Gr_{2n-2k}(V)}\vol_{2k}(\textrm{Pr}_{F^\perp} (B^{2k}))d\sigma_{1}(F) =c_0 \vol_{2k}(B^{2k})=c_0 \frac{\pi^k}{k!}  \]
where $B^{2k}$ is any fixed $2k$-dimensional Euclidean ball, and \[c_0={2n \brack 2k}={2n \choose 2k}\frac{\vol_{2n}(B^{2n})}{\vol_{2k}(B^{2k})\vol_{2n-2k}(B^{2n-2k})}={2n \choose 2k}{n\choose k}^{-1}. \]  We get
\[ A_s=c_0 \frac{\pi^k}{k!}  2\int_{[0,1]^\kappa} (\lambda_1\dots\lambda_\kappa)^{s}d\lambda_1\dots d\lambda_\kappa = 2c_0 \frac{\pi^k}{k!}  \frac{1}{(s+1)^\kappa}  \]
and taking $s=-(2n+1)$ we conclude 
$A=2c_0(-1)^\kappa \frac{\pi^k}{k!}  \frac{1}{(2n)^\kappa}$.
Averaging the right hand side of eq. \eqref{eqn:symplectic_crofton}, we get
\[C \int_{X_{2k}^+}  \vol(B^{2k})|\sigma_\omega(W)|dW
  = C\frac{\pi^k}{k!}\int_{[0,1]^\kappa}\lambda_1\dots\lambda_{\kappa} d\lambda_1\dots d\lambda_{\kappa} =C\frac{\pi^k}{k!2^\kappa}.  \]
Summing up, $ 2c_0 \frac{\pi^k}{k!} (-1)^\kappa \frac{1}{(2n)^\kappa} =\frac{\pi^k}{k!2^\kappa}C$ so that
\[C^{-1}=(-1)^\kappa\frac{(2n)^\kappa}{2^{\kappa+1} c_0}=\frac{(-1)^\kappa}{2}{n\choose k}{2n \choose 2k}^{-1}n^\kappa  \]

\endproof

\begin{appendix}
	
	\section{Oriented valuation theory}\label{app:oriented_valuations}
	In this appendix we draw a common thread between valuation theory and the much simpler theory of closed differential forms and linking integrals.
	
		Let $X^n$ be an oriented manifold. We will think of the closed $k$-forms, denoted $\mathcal Z_k(X)$, as smooth oriented valuations of degree $k$, and consider them as functions on $k$-dimensional oriented submanifolds of $X$ with boundary, given by integration: $\omega(A)=\int_A\omega$. The form is clearly determined by its value on all submanifolds,  which is analogous to the Klain embedding.
		Moreover, $\omega(A)$ only depends on $\partial A$.
		
		The wedge product on $\mathcal Z(X)$ turns it into an algebra. When $X$ is compact, we have Poincar\'e duality, namely that $\mathcal Z_k(X)\otimes \mathcal Z_{n-k}(X)\to \mathcal Z_n(X)=\Omega^n(X)\to \R$ is non-degenerate, where the last arrow is given by $\int_X$.
		
		We also consider the closed currents, which are analogous to the generalized valuations. We will denote them by $\mathcal Z^{-\infty}(X)$. We will sometimes write $\mathcal Z^\infty$ instead of $\mathcal Z$.
		
		Now assume $X=V=\R^n$. The translation-invariant (smooth or generalized) oriented valuations $\mathcal Z^{\pm\infty}(V)^{tr}$ are just $\wedge^\bullet V^*$. The following construction appears in \cite{gelfand_smirnov}.
		
		For oriented manifolds with boundary $A,B\subset V$ of complementary dimensions, at least one of which is compact, let $I(A, B)$ denote their intersection index. It is well-defined when $A, B$ are in general position. 
		Note that for a closed (as a subset) $(n-k)$-dimensional submanifold $E\subset V$, $I(\bullet, E)\in \mathcal Z^{-\infty}_{k}(V)$.
		
		Given a distribution $\mu\in\mathcal M^{-\infty}(\AGr^+_{n-k}(V))$, we define $\Cr(\mu)\in \mathcal Z^{-\infty}_k(V)$ by $\Cr(\mu)=\int_{\AGr^+_{n-k}(V)} I(\bullet, E)d\mu(E)$. Clearly the even measures (with respect to the orientation-reversing map) lie in the kernel of $\Cr$, so we restrict our attention to odd measures.
		
		\begin{Proposition}\label{prop:oriented_crofton}
			The map $\Cr: \mathcal M^{\pm\infty}(\AGr^+_{n-k}(V))^{tr}\to \mathcal Z_k(V)^{tr}$ is surjective.
		\end{Proposition}
		\proof
		The $\GL(V)$-module $Z_k(V)^{tr}=\wedge ^k V^*$ is irreducible. Since the Crofton map is $\GL(V)$-equivariant, it suffices to show $\Cr$ is non-zero. This is not hard to see, for instance $\Cr(\delta_{E}-\delta_{-E})(A)=2I(A, E)$, where $\langle \delta_E, f\rangle=\int_{V/E}f(x+E)dx$ on a compactly supported test function $f$. For a smooth example, one could convolve with an approximate identity on $\GL(V)$.
		\endproof
		
		The analogues of the the Alesker-Poincar\'e and Alesker-Fourier dualities coincide in this setting:
		The Alesker-Poincar\'e pairing is the wedge product $\wedge^kV^*\otimes\wedge^{n-k}V^*\to \wedge^{n}V^*$.
		The Alesker-Fourier duality operation is given by the Hodge star: $\ast: \wedge^kV^*\to (\wedge^{n-k}V^*)^*\otimes \wedge^n V^*$. 
		
		The following easy statement is the analogue for oriented valuations of the principal kinematic formula. Due to the finite-dimensionality of the space of translation-invariant forms, one can average over translations alone.
		
		Let a top form $\vol_n\in \wedge ^n V^*$ be fixed on $V$. Let $\Omega\in \wedge^\bullet V^*\otimes  \wedge^\bullet V^*=(\wedge^\bullet V\otimes  \wedge^\bullet V)^*$ be the wedge product, that is for $\xi, \eta\in \wedge ^\bullet V$, $\langle \xi\wedge\eta, \vol_n\rangle =\Omega(\xi,\eta)$.
		
		For oriented compact manifolds with boundary $A, B$ we write $\Omega(A,B):=\int_{A\times B}\Omega$, which can be written more explicitly by representing $\Omega=\sum \omega_i\otimes \omega'_{i}$, then $\Omega(A, B)=\sum \int _A \omega_i\int_B\omega'_{i}$ (the integrals with mismatched dimension vanish by definition). Consider the kinematic operator
		\[K_V(A,B):=\int_V I(A, B+x)d\vol_n(x)\]
		which is well-defined, since the integrand is compactly supported.
		\begin{Proposition}
			Let $A^{n-k}, B^k\subset V$ be compact, oriented submanifolds with boundary, of complementary dimension. Then $K_V(A, B)=\Omega(A, B)$
		\end{Proposition}
		\proof
		Let $[[B]]$ be the current defined by $B$. Then $[B]_V:=\int_V [[B+x]]d\vol_n(x)$ is translation-invariant, that is $[B]_V\in\wedge^kV$. Now $I(A, B)=[[A]]\cap [[B]]$, so that $K_V(A, B)=[[A]]\cap[B]_V=\vol_n([A]_V\wedge[B]_V)=\Omega(A, B)$ as claimed.

		\endproof 
		We remark that this formula is also reminiscent of the Bezout formula in complex algebraic geometry.

	\section{Invariant sections}\label{app:invariants}
	We will need two technical lemmas concerning invariants of group actions.
	
	The first goes back to Kolk and Varadarajan \cite{kolk_varadarajan}, and appeared in a form best suitable for our needs in \cite{bernig_faifman}. Let us quote the result in its simplest sufficient form.

	Take a Lie group $G$ acting on a manifold $X$ with finitely many orbits, all locally closed submanifolds. Let $E$ be a $G$-vector bundle over $X$. Define for integer $\alpha\geq 0$ and a submanifold $Y\subset X$ the $G$-bundle $F^\alpha_Y$ over $Y$ by \begin{displaymath}
	F^\alpha_Y|_y=E|_y \otimes \Dens^*(N_yY) \otimes\Sym^\alpha(N_yY) .                                                                                                                                                                                                                                                                                                                                                                                                                                                                                                                                                                                    \end{displaymath}
	A generalized section $s\in\Gamma^{-\infty}_{\overline Y}(X,E)$ has a certain transversal order $\alpha$ along $Y$, and a transversal principal symbol $\sigma(s)\in \Gamma^{-\infty}(Y, F^\alpha_Y)$. For details, see e.g. \cite[section 4.4.]{alesker_faifman}.
	\begin{Lemma}	\label{lem:LocallyClosedOrbits}
		Let $Z\subset X$ be a closed $G$-invariant subset.
		Decompose into $G$-orbits: $Z=\bigcup_{j=1}^J Y_j$. Then
		\begin{displaymath}
		\dim \Gamma^{-\infty}_Z(X,E)^G \leq \sum_{\alpha=0}^\infty \sum _{j=1}^J \dim \Gamma^\infty(F^\alpha_{Y_j})^{G}.
		\end{displaymath}
	\end{Lemma}
	\vspace{5pt}
	The second statement is surely well-known, but we were not able to locate a reference. For completeness, we include the proof.
	\begin{Proposition}\label{prop:infinite_dimensional_generalized}
	Let $G$ be a (possibly infinite dimensional) Lie group acting on a manifold $M^n$ transitively, and let $\mathcal E$ be an infinite-dimensional $G$-bundle of Fr\'echet spaces over $X$. Then the space of $G$-invariants of $\Gamma_c^\infty(M, \mathcal E)^*$ is naturally isomorphic to $\Gamma^\infty(M, \mathcal E^*\otimes |\omega_M|)^G$
\end{Proposition}
\begin{Remark} In the infinite dimensional case, we assume that every $T_x M$ admits a basis of infinitesimal generators of 1-parametric subgroups of $G$. This is certainly the case for $Diff(M)$, $Symp(M)$, $Cont(M)$.
\end{Remark}
First we consider the wavefront set of the average of a distribution along a flow.	
	\begin{Lemma}\label{lem:average_one_direction}
		Let $M$ be a manifold. Assume $\R$ acts on $M$, and let the curve $C\subset M$ be an orbit.
		Let $\phi\in C^{-\infty}(M)$ be a generalized function, and $\mu\in C^\infty_c(\R)$. 
		Define $\phi\ast \mu:=\int_{\R} t^*\phi\cdot \mu(t)dt$. Then for all $p\in C$ we have $\WF(\phi\ast\mu)\cap T^*_pM\subset N_p^*C$.
	\end{Lemma}
\proof
Fix a small neighborhood $V$ of $p\in C$ which is $\R$-equivariantly identified with a neighborhood $U\subset\R^n$ of $p=0$, with $\R$ acting by translations along the $x_1$-axis and $C$ coinciding with the $x_1$-axis.
Choose a partition of unity $\rho_i$ on $M$ and $w_j$ on $\R$. Writing $\phi_i:=\rho_i\cdot \phi$, $\mu_j:=w_j\cdot \mu$ we have $\phi\ast{\mu}=\sum_{i,j} \phi_i\ast{\mu_j}$.  We may assume that each $\phi_i$ has finite order, and that each $\psi=\phi_i\ast {\mu_j}$ whose support contains $p$, is in fact supported inside $V$. Taking such $\psi$, we may write
$\psi=\Phi\ast\pi_1^*\nu$ for some $\nu\in C^\infty_c(\R)$ and $\Phi\in C_c^{-\infty}(\R^n)$ of order $k$, where $\pi_1:\R^n\to \R$ is the projection to the first coordinate.
Then the Fourier transform satisfies $\widehat \psi(\xi) =\widehat \Phi(\xi)\widehat \nu(\xi_1)$.
Now $|\widehat \Phi(\xi)|\leq C(1+|\xi|)^k$ while $|\widehat \nu(\xi_1)|\leq C_N(1+|\xi_1|)^{-N}$ for all $N$. It follows that for all $\epsilon>0$, the cone
$C_\epsilon:=\{\xi_1^2>\epsilon(\xi_2^2+\dots+\xi_n^2)\}$ falls outside the wavefront set of $\psi$, so that $\WF(\psi)\subset\{\xi_1=0\}$.
That concludes the proof.
\endproof
Clearly the same statement holds also for $S^1$-actions, and with $\phi$ a generalized section of an arbitrary $\R$-equivariant vector bundle over $M$.

	\textit{Proof of Proposition \ref{prop:infinite_dimensional_generalized}.}
	Take $s\in\Gamma_c^{\infty}(M, \mathcal E)^\ast$. For $\phi\in\Gamma_c^\infty(M,\mathcal E)$ one can define $s\cdot \phi \in\mathcal M^{-\infty}(M)$ by $\int f\cdot d(s\cdot \phi):=s(f\phi)$ for $f\in C^\infty_c(M)$.
	Moreover, the map $\Gamma_c^{\infty}(M, \mathcal E)\to \mathcal M^{-\infty}(M)$ given by $\phi\mapsto s\cdot\phi$ is $G$-equivariant.
	Fix $\phi$, and let us verify that $s\cdot\phi$ is in fact a smooth measure. 
	Consider a 1-parametric subgroup $H\subset G$ which can be either $\R$ or $S^1$. By the Dixmier-Mallavin theorem, we may find smooth probability measures $\mu_1,\dots,\mu_N\in\mathcal M^\infty_c(H)$ and sections $\psi_1,\dots,\psi_N\in \Gamma_c^\infty(M,\mathcal E)$ such that \[\phi=\sum_j \int_H g\psi_j d\mu_j(g)\Rightarrow s\cdot \phi=\sum_j \int_{H} g(s\cdot \psi_j)d\mu_j(g)\] 
 	It follows by Lemma \ref{lem:average_one_direction} that the wavefront set $\WF(s\cdot \phi)$ lies in the conormal bundle to the orbits of $H$ on $M$. Since $G$ acts transitively, we conclude $\WF(s\cdot \phi)=\emptyset$, that is $s\cdot\phi \in\mathcal M^{\infty}(M)$.
	
	In particular, we can consider the density on every tangent plane, $(s\cdot\phi) (x)\in\Dens(T_xM)$.
	We next claim that $(s\cdot\phi)(x)$ only depends on $\phi(x)$. Indeed, if $\phi(x)=0$ we may represent $\phi=f\cdot\psi$ for some $\psi\in\Gamma_c^\infty(M,\mathcal E)$ and $f_c\in C^\infty(M)$ with $f(0)=0$. Then $(s\cdot \phi)=f\cdot (s\cdot \psi)$, and thus $(s\cdot \phi)(x)=0$. That is, we get an element $s(x)\in \mathcal E^*|_x\otimes \Dens(T_xM)$. By $G$-invariance, $x\mapsto s(x)$ is a smooth section of $\mathcal E^*\otimes |\omega_M|$, which clearly defines the same functional as $s$ on $\Gamma_c^{\infty}(M, \mathcal E)$. \qed

	\end{appendix}

\bibliographystyle{amsalpha}
\bibliography{symplectic}
\end{document}